\documentstyle[11pt,amssymb,amsfonts,twoside]{article}

\begin{titlepage}

\title{Regularity of invariant distributions}

\author{Ulrich Bunke\thanks{Mathematisches Institut, Universit\"at G\"ottingen, Bunsenstr. 3-5, 37073 G\"ottingen,
 GERMANY, E-mail : bunke@uni-math.gwdg.de} and
Martin
Olbrich\thanks{Mathematisches Institut, Universit\"at G\"ottingen, Bunsenstr. 3-5, 37073 G\"ottingen,
 GERMANY, E-mail : olbrich@.uni-math.gwdg.de  }
}

\end{titlepage}

\newcommand{\proof}{{\it Proof.$\:\:\:\:$}}
 \newcommand{\dist}{{\rm dist}}
\newcommand{\kaaa}{{\frak k}}
\newcommand{\paaa}{{\frak p}}
\newcommand{\vp}{{\varphi}}

\newcommand{\R}{{\Bbb R}}

\newcommand{\Z}{{\Bbb Z}}

\newcommand{\Fam}{{\rm Fam}}
\newcommand{\Cusp}{{\rm Cusp}}

\newcommand{\C}{{\Bbb C}}

\newcommand{\gaaa}{{\frak g}}

\newcommand{\aaaa}{{\frak a}}
\newcommand{\naaa}{{\frak n}}

\newcommand{\cE}{{\cal E}}

\newcommand{\cU}{{\cal U}}

\newcommand{\vol}{{\rm vol}}
\newcommand{\cO}{{\cal O}}
\newcommand{\End}{{\mbox{\rm End}}}
\newcommand{\Ext}{{\mbox{\rm Ext}}}
\newcommand{\rk}{{\rm rank}}

\newcommand{\cF}{{\cal F}}

\newcommand{\Ree}{{\rm Re }}

\newcommand{\Imm}{{\rm Im}}

\newcommand{\id}{{\mbox{\rm id}}}
\newcommand{\ord}{{\rm ord}}
\newcommand{\nat}{{\Bbb  N}}
\newcommand{\supp}{{\rm supp}}

\newcommand{\aca}{{\aaaa_\C^\ast}}

\def\imath{{\rm i}}

\newcommand{\cR}{{R}}
\def\hB{\hspace*{\fill}$\Box$ \newline\noindent}

\newcommand{\cS}{{S}}

\def\hB{\hspace*{\fill}$\Box$ \\[0.5cm]\noindent}

\newcommand{\cP}{{\cal P}}

\newtheorem{prop}{Proposition}[section]
\newtheorem{lem}[prop]{Lemma}
\newtheorem{ddd}[prop]{Definition}
\newtheorem{theorem}[prop]{Theorem}
\newtheorem{kor}[prop]{Corollary}

\newtheorem{con}[prop]{Conjecture}

\begin{document}
\setcounter{page}{1}
\maketitle

\begin{abstract}
We consider a discrete subgroup $\Gamma$ of the isometry group $G$ of the real
hyperbolic space $X\cong H^n$ of dimension $n\ge 2$. 
We assume that $\Gamma$ admits a fundamental domain with finitely many totally
geodesic faces (i.e. $\Gamma$  is geometrically finite). This domain may have
finitely many cusps of rank varying in $\{1,\dots,n-1\}$.
$\Gamma$ acts by
conformal transformations on the sphere $\partial X\cong S^{n-1}$ at infinity.
By $\Lambda_\Gamma$ we denote its limit set which is a closed
$\Gamma$-invariant subset of $\partial X$ of Hausdorff dimension
$\dim_H(\Lambda_\Gamma)$. We consider a  $G$-equivariant irreducible
complex vector bundle $V\rightarrow \partial X$.   Let  $\Lambda$  be the
bundle of densities on $\partial X$. Then we define a real number $s(V)$ by
the condition that $V\otimes \Lambda^{s(V)} \cong  V^\sharp\otimes \Lambda$ as
$G$-equivariant bundles, where $V^\sharp$ is the hermitian dual of $V$.
The goal of the paper  is to find lower bounds on  the regularity of
$\Gamma$-invariant distribution sections of $V$ which are strongly supported
on the limit set. In order to quantify this regularity we consider the spaces
of sections of regularity $H^{p,r}$, $p\in (1,\infty)$, $r\in\R$. For
$r\in\nat_0$ the space $H^{p,r}$ is the space of distribution sections
which have all their derivatives up to order $r$ in $L^p$. For $r\in
[0,\infty)$ we extend this scale by complex interpolation, while for negative
$r$ we define these spaces by duality. 
Our main result is the following. \\
{\em Let
\begin{eqnarray*}
r^0&:=&(n-1)\frac{s(V)-1}{2}+\frac{n-1-\dim_H(\Lambda_\Gamma)}{p}\\
r&:=&\min_{} \left( r^0,(n-1)s(V)+\frac{n-1}{p}
-\max_{cusps\:c} \rk(c) \right)\ , \end{eqnarray*}
where the maximum over the empty set is, by convention, equal to $-\infty$.
Moreover, $\dim_H(\Lambda_\Gamma)$ has to be replaced by half of the rank of
the cusp if $\Gamma$ is elementary parabolic. 
If $\phi$ is a $\Gamma$-invariant 
distribution section of $V$ which is strongly supported on the limit set (and a
cusp form), then $\phi\in H^{p,s}$ for all $s<r$ ($s<r^0$).}\\[1cm] If
$\Gamma$ is convex cocompact, then strongly supported on the limit set  just
means supported on the limit set as a distribution (and the cusp form
condition holds true automatically). But if $\Gamma$ has cusps, then this
condition is stronger. Generically, boundary values of cusp forms and residues
of Eisenstein series have this property. 
Actually, in order to prove our
result, we also need some additional technical assumptions.
\end{abstract}

\parskip1ex
\tableofcontents
\parskip3ex
\section{Statement of the result}
\subsection{Introduction and statement of the main theorem}

Let $G$ be one of the groups $Spin(1,n)$, $ SO(1,n)_0$, $2\le n\in\nat$. Then
the symmetric space $X$ of $G$ is the real hyperbolic space of dimension $n$.
Furthermore, let $\Gamma\subset G$ be a discrete subgroup. We assume that
$\Gamma$ is {\em geometrically finite}. One of many equivalent formulations  of
this condition is that $\Gamma$ admits a fundamental domain in $X$ which has
finitely many totally geodesic faces. In the following we assume that $\Gamma$
is {\em torsion-free}, but we shall see later that one can drop this
assumption.

By $\partial X$ we denote the geodesic boundary of $X$. It can be identified with the space
 of all minimal parabolic subgroups of $G$ such that 
the minimal parabolic subgroup $P\subset G$ corresponds to its unique fixed
point $\infty_P\in \partial X$.  Since $G$ acts transitively  on the set of
minimal  parabolic subgroups we can identify $\partial X \cong G/P$ for any
fixed $P$. 

The $G$-homogeneous complex vector bundles on $\partial X$ are in one-to-one
correspondence with finite-dimensional complex representations of $P$. On the
one hand if $(\theta,V_\theta)$ is a finite dimensional complex representation
of $P$, then we can form the associated bundle $V(\theta):=G\times_{P,\theta}
V_\theta$. On the other hand, if $V\rightarrow \partial X$ is a
$G$-homogeneous complex vector bundle, then we let $V_\theta$ be the
representation of $P$ on the fibre of $V$ over $\infty_P$, and we  have a
natural isomorphism $V(\theta)\cong V$.

In order to parametrize the irreducible representations of $P$ we consider 
the exact sequence  
$$0\rightarrow N\rightarrow P\rightarrow L\rightarrow 0\
,$$ where $N$ is the unipotent radical of $P$. Any finite-dimensional
irreducible representation of $P$ factors over the quotient $L$. This group
fits into a natural exact sequence 
\begin{equation}\label{i98}0\rightarrow M\rightarrow L\rightarrow
A\rightarrow 0\ ,\end{equation}
 where $M$ is a compact group isomorphic to
$Spin(n-1)$ or $SO(n-1)$, and $A$ is isomorphic to the multiplicative group
$\R^+$. Let $\aaaa$ denote the Lie algebra of $A$ and $\aca$ be the
complexification of its dual. 

For $\lambda\in\aca$ we form the character $a\mapsto
a^\lambda:=\exp(\lambda(\log(a))$. The group $L$ acts on the Lie algebra
$\naaa$ of $N$. The induced representation on $\Lambda^{\dim(\naaa)}\naaa$
factors over $A$ and corresponds to the character $2\rho\in\aaaa^*$.
Furthermore, we define the element $\alpha\in\aaaa^*$ by
$$\alpha:= \frac{2\rho}{\dim(\naaa)}$$
which plays the role of a unit.

If $(\theta,V_\theta)$ is any irreducible representation of $L$, then its
restriction to $M$ is an irreducible  unitary representation. Furthermore,
the representation $\Lambda^{\dim(V_\theta)} \theta$ factors over $A$ and
corresponds to an element of $\aca$ (except in the case $G=SL(2,\R)$.
In this case $M$ may act non-trivially on $\Lambda^{\dim(V_\theta)} \theta$
and one employs the canonical split of
(\ref{i98}) in order to define the character of $A$ corresponding to
$\theta$).  \begin{ddd} Given $(\sigma, V_\sigma) \in\hat M$ and
$\lambda\in\aca$ we let $(\sigma_\lambda,V_{\sigma_\lambda})$ be the unique
irreducible representation of $P$ such that
$(\sigma_\lambda)_{|M}=\sigma$ and $\Lambda^{\dim V_\sigma}\sigma_\lambda$
corresponds to $\dim(V_\sigma) (\rho-\lambda)\in\aca$. 
\end{ddd}
The bundle of densities $\Lambda$ on $\partial X$ is given in this
parametrization by $V(1_{-\rho})$, where $1\in\hat M$ denotes the trivial
representation.  The bundle $V(\sigma_\lambda)^\sharp\otimes\Lambda$ 
is isomorphic to 
$V(\sigma_{-\bar \lambda})$ as $G$-homogeneous bundle. Note that $V(1_{-\rho})^s\otimes
V(\sigma_\lambda)=V(\sigma_{\lambda-2s\rho})$.
We conclude that 
$$s(V(\sigma_\lambda))=\frac{\Ree(\lambda)}{\rho}\ ,$$ where
$s(V(\sigma_\lambda))$ is the number introduced in the abstract.

Another important element of $\aaaa^*$ is the critical exponent $\delta_\Gamma$
of $\Gamma$ which we now define. We fix any maximal compact subgroup $K\subset
G$. This choice induces an embedding of $L$ into $P$ such that $L$ is stable
under the Cartan involution $\theta_K$ associated with $K$. Furthermore, we
obtain a decomposition $L=MA$, where $M=L^{\theta_K}$ is the set of fixed
points of $\theta_K$, and if $a\in A$, then $\theta_K(a)=a^{-1}$.
In particular we now can consider $A$ as a subgroup of $G$.
 The character $\rho$ fixes a semigroup $A_+:=\{a\in A\:|\: a^\rho\ge 1\}$.
By the Cartan decomposition any element $g\in G$ can be written as $g=k_ga_g
h_g$, where $k_g,h_g\in K$, and $a_g\in A_+$ is uniquely determined.

\begin{ddd} The {\em critical exponent}
$\delta_\Gamma\in\aaaa^*$   of $\Gamma$ is defined as the infimum of
the set $\{\mu\in \aaaa^*\:|\: \sum_{\gamma\in\Gamma}
a_\gamma^{-\mu-\rho}<\infty\}$. \end{ddd}

Let $\partial X=\Omega_\Gamma\cup\Lambda_\Gamma$ be the $\Gamma$-invariant
decomposition into the limit set $\Lambda_\Gamma$ and the ordinary set
$\Omega_\Gamma$. The limit set $\Lambda_\Gamma$ is closed and can be defined as
the set of accumulation points of any $\Gamma$-orbit in $X$. By results of
Patterson \cite{patterson76} in case $n=2$ and Sullivan \cite{sullivan79} for
all $n\ge 2$ the Hausdorff dimension of $\Lambda_\Gamma$ is given by
$\dim_H(\Lambda_\Gamma)=\frac{\delta_\Gamma +\rho}{\alpha}$ unless $\Gamma$ is
elementary parabolic.

Recall that a parabolic subgroup $P\subset G$ is called $\Gamma$-cuspidal if 
$\Gamma\cap P$ is infinite and projects to a precompact subset $l(\Gamma\cap
P)\subset L$, where $l:P\rightarrow L$ is the projection. We
form the compact group $M_\Gamma:=\overline{l(\Gamma\cap P)}\subset M$.
Furthermore, there exists a subgroup $N_\Gamma\subset N$ and a Langlands
decomposition $P=MAN$ (choice of $K$) such that $\Gamma_P:=\Gamma \cap
P$ is a cocompact lattice in the unimodular group $P_\Gamma:=N_\Gamma M_\Gamma$. 

A {\em cusp} of $\Gamma$ is a $\Gamma$-conjugacy class $[P]_\Gamma$ of
$\Gamma$-cuspidal parabolic subgroups. The dimension $\dim(N_\Gamma)$ is
called {\em the rank of the cusp} $\rk([P]_\Gamma)$ defined by $P$. We
further define $\rho_{\Gamma_P}:= \frac12 \rk([P]_\Gamma)\alpha $ and
$\rho^{\Gamma_P}:=\rho-  \rho_{\Gamma_P}$.  If $\dim(N_\Gamma)\in
\{1,\dots,\dim(\naaa)-1\}$, then we say that the cusp has {\em smaller rank}.
If $\dim(N_\Gamma)=\dim(\naaa)$, then the cusp has {\em full rank}.

A finite-dimensional representation $(\vp,V_\vp)$ of $\Gamma$ is called
{\em admissible twist} (shortly twist) if for any $\Gamma$-cuspidal
parabolic subgroup $P$ the representation $\vp$ extends to $AP_\Gamma$ such
that $A$ acts by semisimple endomorphisms. For example, if $\vp$ is a
finite-dimensional representation of $G$, then its restriction to $\Gamma$ is
an admissible twist.

Given a representation $(\theta,V_\theta)$ of $P$ and an admissible twist
$(\vp,V_\vp)$ of $\Gamma$ we form the bundle $V(\theta,\vp)=V(\theta)\otimes
V_\vp$. If $\pi^\theta$ denotes the left-regular representation of $G$ on 
sections of $V(\theta)$, then we can form the representation
$\pi^{\theta,\vp}:=\pi^\theta\otimes \vp$ of $\Gamma$ on sections
of $V(\theta,\vp)$. Note that $\pi^{\theta,\vp}$ extends to $AP_\Gamma$ for
each $\Gamma$-cuspidal parabolic subgroup $P$.

\begin{ddd} The {\em exponent $\delta_\vp$ of the twist $\vp$}
is defined as the infimum of
the set $\{\mu\in \aaaa^*\:|\: \sup_{\gamma\in\Gamma}
\|\vp(\gamma)\| a_\gamma^{-\mu}<\infty\}$,
where $\|.\|$ is any norm on $\End(V_\vp)$.  \end{ddd}

Let $\tilde \sigma$ or $\tilde \vp$ denote the $\C$-linear dual
representations of $\sigma$ or $\vp$.
The spaces $$C^{-\infty}(\partial X,V(\sigma_\lambda,\vp)):=C^\infty(\partial
X,V(\tilde\sigma_{-\lambda},\tilde\vp))^*$$ for varying $\lambda\in\aca$
form locally trivial holomorphic bundles of dual Fr\'echet
spaces so that it  makes sense to speak of germs holomorphic families 
$(\phi_\lambda)_{\lambda}$ such that   $\phi_\lambda\in C^{-\infty}(\partial
X,V(\sigma_\lambda,\vp))$.  Let  ${}^\Gamma C^{-\infty}(\partial
X,V(\sigma_\lambda,\vp))$ be the subspace of $\Gamma$-invariant distributions. 

We call $\phi\in {}^\Gamma
C^{-\infty}(\partial X,V(\sigma_\lambda,\vp))$ {\em deformable}, if there is a
germ at $\lambda$ of a holomorphic family $(\phi_\mu)_{\mu}$  such that $\phi_\lambda=\phi$ and $\phi_\mu\in {}^\Gamma
C^{-\infty}(\partial X,V(\sigma_\mu,\vp))$.
\begin{ddd}
Let $\Fam_\Gamma(\sigma_\lambda,\vp)\subset {}^\Gamma
C^{-\infty}(\partial X,V(\sigma_\lambda,\vp))$ be the subspace of all
deformable $\Gamma$-invariant distributions. 
\end{ddd}

If $\phi\in\Fam_\Gamma(\sigma_\lambda,\vp)$, then
we say that $\phi$ is {\em strongly supported on the limit set} if
$res^\Gamma(\phi)=0$. Here $res^\Gamma$ is the restriction map introduced in
\cite{bunkeolbrich99} (we refer to Subsection \ref{reex} for details). More
precisely, this condition means that there exists a germ at $\lambda$ of a
holomorphic  family $(\phi_\mu)_{\mu}$, $\phi_\mu\in {}^\Gamma
C^{-\infty}(\partial X,V(\sigma_\mu,\vp))$, such that $\phi_\lambda=\phi$
and $(res^\Gamma \phi_\mu)_{|\mu=\lambda}=0$.

If $\phi$ is strongly supported on the limit set, then its support as a
distribution $\supp(\phi)$ is contained in $\Lambda_\Gamma$. (In fact it is
equal to $\Lambda_\Gamma$ if $\sharp \Lambda_\Gamma\not=2$.)  
The converse may be false  if  $\Gamma$ has cusps.

\begin{ddd}
Let $\Fam_\Gamma(\Lambda_\Gamma,\sigma_\lambda,\vp)\subset
\Fam_\Gamma(\sigma_\lambda,\vp)$ be the subspace of all deformable invariant
distributions which are strongly supported on the limit set.
\end{ddd}

Being deformable is a non-trivial condition. For example, if $\Gamma$ has
finite-covolume, then the boundary values of cusp forms are not deformable,
while the boundary values of Eisenstein series are deformable.
We also  have examples of non-deformable invariant distributions in the case
of convex cocompact $\Gamma$ (see \cite{bunkeolbrich990},
\cite{olbrichhabil}, \cite{bunkeolbrich00}). 

Using {\em  embedding} and {\em twisting}  we can formulate a weaker
condition. 
We first describe embedding. For a moment we write $G_n$ for
$Spin(1,n)$ or $SO(1,n)_0$, and we attach the subscript $n$ to all related
objects. For $n\ge m$ have inclusions $G_m\subset G_n$, $\partial X_m\subset
\partial X_n$,  $P_m\subset P_n$, $M_m\subset M_n$. Note that $A_m\cong
A_n=:A$ and thus $\aca_m\cong \aca_n=:\aca$ in a natural way.

Let $\sigma_n\in\hat M_n$, $ \sigma_m\in \hat M_m$ and assume that
there is an $M_m$-invariant inclusion  $T:V_{
\sigma_m}\rightarrow V_{\sigma_n}$. Then $T^*$ induces a projection 
 $V^n((\tilde \sigma_n)_{-\lambda-\rho^m+\rho^n})_{|\partial X^n}\rightarrow
V^m((\tilde \sigma_m)_{-\lambda})$. Let
$$i_{n, \sigma_n,T}^*:C^\infty(\partial X^n, V^n((\tilde
\sigma_n)_{-\lambda-\rho^m+\rho^n},\tilde\vp))\rightarrow C^\infty( \partial
X^m, V^m((\tilde \sigma_m)_{-\lambda},\tilde\vp))$$ be the restriction of
functions and  $$i^{n,\sigma_n,T}_*:C^{-\infty}( \partial X^m, V^m((
\sigma_m)_\lambda,\vp))\rightarrow  C^{-\infty}(\partial X^n, V^n((
\sigma_n)_{\lambda+\rho^m-\rho^n},\vp))$$ be the dual map. Note that
$i^{n,\sigma_n,T}_*$ is injective and $G_m$-equivariant. For fixed $m$ and
$\sigma_m$ the various choices of $n$, $\sigma_n$ and  $T$ are are called
embedding data. We will often  abbreviate $i^{n,\sigma_n,T}_*=i_*$.

Now we explain twisting.
If $(\pi,V_\pi)$ is a finite-dimensional representation of $G$,
then there is an isomorphism of $\Gamma$-equivariant bundles
$R:V(\theta\otimes\pi_{|P},\vp)\stackrel{\sim}{\rightarrow} V(\theta,\vp\otimes
\pi)$ given by $R([g,x\otimes y]\otimes z):=[g,x]\otimes z
\otimes \pi(g)y$.
Fix $\sigma\in\hat M$ and $\lambda\in \aca$. Assume that we are given
$\sigma^\prime\in\hat M$, $\mu\in\aca$, a finite-dimensional representation of
 $(\pi,V_\pi)$ of $G$, and a $P$-equivariant inclusion
$T:V_{\sigma_\lambda}\hookrightarrow V_{\sigma^\prime_\mu}\otimes V_{\pi}$.
Then we obtain a $\Gamma$-equivariant inclusion
$$i_*^{\sigma^\prime,\mu,\pi,T}:C^{-\infty}(\partial X,
V(\sigma_\lambda,\vp))\rightarrow   C^{-\infty}(\partial X,
V(\sigma^\prime_\mu\otimes \pi_{|P},\vp))\stackrel{R}{\rightarrow}  C^{-\infty}(\partial X,
V(\sigma^\prime_\mu,\vp\otimes \pi))\ .$$
Note that $i_*^{\sigma^\prime,\mu,\pi,T}$ is $P_\Gamma$-equivariant for any
$\Gamma$-cuspidal parabolic subgroup $P$.
For simplicity we will from now on assume that $(\pi,V_\pi)$ is
{\em irreducible}. The various choices of $\sigma^\prime$, $\mu$, $\pi$ and $T$
are called twisting data. We often abbreviate
$i_*^{\sigma^\prime,\mu,\pi,T}=i_*$.

We also allow iterated twisting and embedding, and we still denote the
resulting inclusion by $i_*$.
We say that $\phi\in  {}^\Gamma
C^{-\infty}(\partial X,V(\sigma_\lambda,\vp))$ is {\em stably deformable} if
$i_*\phi$ becomes deformable for suitable embedding  or twisting data.
\begin{ddd} Let $\Fam^{st}_\Gamma(\sigma_\lambda,\vp)$ denote the subspace of
stably deformable invariant distributions.
\end{ddd}

Being stably deformable is a weaker condition than being deformable. 
E.g., if $\Gamma$ has finite covolume, then the boundary values of cusp forms
are not deformable, but generically they are stably deformable.
In fact, for most $\lambda\in\aca$ all invariant distributions are stably
deformable. Let $I_\aaaa\subset \aaaa^*$ denote the lattice of half-integral
characters spanned by $\frac{1}{2}\alpha$.
\begin{theorem}\label{gfa}
If  $\lambda\not\in I_\aaaa$, then 
$\Fam^{st}_\Gamma(\sigma_\lambda,\vp)={}^\Gamma
C^{-\infty}(\partial X,V(\sigma_\lambda,\vp))$.
\end{theorem}
A proof will appear in a future publication. One of the main steps is to show
that $res^\Gamma$ is regular (i.e. $res ^\Gamma \phi_\mu$ is holomorphic
for any germ at $\lambda$ of a holomorphic family $(\phi_\mu)_{\mu}$
of $\Gamma$-invariant distributions) for
$\lambda\not\in I_\aaaa$ for any $\Gamma$ which is elementary parabolic, i.e.
$\Gamma=\Gamma_P$ for some $\Gamma$-cuspidal parabolic subgroup.

The following was shown in \cite{bunkeolbrich00}, proof of Cor. 6.12.
\begin{theorem}\label{ccalways}
If $\Gamma$ is convex cocompact, then
$$\Fam^{st}_\Gamma(\sigma_\lambda,\vp)={}^\Gamma C^{-\infty}(\partial
X,V(\sigma_\lambda,\vp))$$ holds for any $\lambda\in\aca$.
\end{theorem}
If $\Gamma$ is convex cocompact and not cocompact, then twisting alone
suffices to make an invariant distribution deformable. 
We do not know any example of an invariant distribution which is not stably
deformable.   
\begin{con}
If $\Gamma$ is geometrically finite, then
$$\Fam^{st}_\Gamma(\sigma_\lambda,\vp)={}^\Gamma
C^{-\infty}(\partial X,V(\sigma_\lambda,\vp))$$
for all $\sigma$, $\lambda$ and $\vp$.
\end{con}
Note that if $\phi\in 
\Fam^{st}_\Gamma(\sigma_\lambda,\vp)$, then we can always find
embedding and twisting data such that
$i_*\phi\in \Fam_\Gamma(1_{\lambda^\prime},\vp^\prime)$.
This will be employed later.

If $\phi\in \Fam^{st}_\Gamma(\sigma_\lambda,\vp)$, then we say that it is
{\em strongly supported on the limit set} if $i_*\phi$ is strongly supported
on the limit set for some  choice of embedding and twisting data with the
property  that $i_*\phi$ is deformable.  
This definition makes sense since being strongly supported on the limit set is
stable under twisting and embedding by Lemma \ref{stablesupp}.

\begin{ddd}
We let $\Fam^{st}_\Gamma(\Lambda_\Gamma,\sigma_\lambda,\vp)$ be the space of
stably deformable invariant distribution vectors which are strongly supported
on the limit set.
\end{ddd}

An alternative way to sharpen the condition "supported on the limit set" 
is to require a "cusp form" condition. If $\lambda\not\in I_\aaaa$, then 
"cusp form" implies "strongly supported on the limit set" (Theorem \ref{cfad}).

Let $P$ be a $\Gamma$-cuspidal parabolic subgroup. 
We normalize the Haar measure $dx$ on the unimodular group $P_\Gamma$ such that
$\vol(\Gamma_P\backslash P_\Gamma)=1$.  Note that $\pi^{\sigma_\lambda,\vp}$
extends to $P_\Gamma$.
\begin{ddd}
If $\phi\in {}^\Gamma
C^{-\infty}(\partial X,V(\sigma_\lambda,\vp))$, then we define its {\em
constant term} $$\phi_P\in {}^{P_\Gamma} C^{-\infty}(\partial
X,V(\sigma_\lambda,\vp))$$  with respect to $P$ by
$$\phi_P:=\int_{\Gamma_P\backslash P_\Gamma} \pi^{\sigma_\lambda,\vp}(x) \phi\:
dx\ .$$  \end{ddd} 
\begin{ddd}
We call $\phi\in {}^\Gamma
C^{-\infty}(\partial X,V(\sigma_\lambda,\vp))$ a {\em cusp form} 
if $\supp(\phi)\subset \Lambda_\Gamma$ and $\phi_P=0$ for all $\Gamma$-cuspidal
parabolic subgroups $P$.
\end{ddd}

If $\Gamma$ has finite covolume and $\phi\in {}^\Gamma
C^{-\infty}(\partial X,V(\sigma_\lambda))$ is a cusp form, then the matrix
coefficients  $G\ni g\mapsto \langle \phi,
\pi^{\tilde\sigma_{-\lambda}}(g) f\rangle$ for $K$-finite $f\in
C^\infty(\partial X,V(\tilde\sigma_{-\lambda}))$ are cusp forms in the
classical sense. 
If $\Gamma$ is convex cocompact, then any $\phi\in {}^\Gamma
C^{-\infty}(\partial X,V(\sigma_\lambda,\vp))$ with
$\supp(\phi)\subset\Lambda_\Gamma$ is a cusp form for trivial reasons.
If $\phi$ is a cusp form, then $i_*\phi$ is also a cusp form for any choice of
embedding or twisting data.

\begin{ddd}
Let $\Cusp_\Gamma(\sigma_\lambda,\vp)$ denote the space of all cusp forms in 
${}^\Gamma C^{-\infty}(\partial X,V(\sigma_\lambda,\vp))$.
\end{ddd}
\begin{theorem}\label{cfad}
If $\lambda\not\in I_\aaaa$, then $\Cusp_\Gamma(\sigma_\lambda,\vp) \subset
\Fam^{st}_\Gamma(\Lambda_\Gamma,\sigma_\lambda,\vp)$. 
\end{theorem}

Note that the proof of this theorem  uses  Thm. \ref{gfa} and the fact that
$res^\Gamma$ is regular outside $I_\aaaa$. Assuming these facts in the present
paper we show that a cusp form is strongly supported on the limit set.

Since $\partial X$ is a closed manifold there are natural Sobolev spaces
$H^{p,r}(\partial X,V(\sigma_\lambda,\vp))$, $p\in (1,\infty)$, $r\in\R$.
For $r\in \nat_0$ we let $H^{p,r}(\partial X,V(\sigma_\lambda,\vp))$ be the
space of those distributions in $C^{-\infty}(\partial
X,V(\sigma_\lambda,\vp))$  which have $L^p$-integrable derivatives up to
order $r$. For $r\in -\nat$ we define $H^{p,r}(\partial
X,V(\sigma_\lambda,\vp))$ to be the dual of $H^{q,-r}(\partial X,V(\tilde
\sigma_{-\lambda},\tilde\vp))$, where $\frac{1}{p}+\frac{1}{q}=1$. For
non-integral $r$ we define $H^{p,r}(\partial X,V(\sigma_\lambda,\vp))$ by
complex interpolation (see Subsection \ref{cint}). We also define
$$H^{p,<r}(\partial X,V(
\sigma_{\lambda},\vp)):=\bigcap_{s<r} H^{p,s}(\partial X,V(
\sigma_{\lambda},\vp))\ .$$

Now we can state our main result.
 
\begin{theorem}\label{main}
Let $r_{p,\lambda}(\Gamma),r^0_{p,\lambda}(\Gamma)\in\R$ be determined by
\begin{eqnarray*}
r^0_{p,\lambda}\alpha& :=&\Ree(\lambda) -\delta_\vp- \rho +
\frac{\rho-\delta_\Gamma}{p}\\
r_{p,\lambda}(\Gamma)\alpha&:=&\min\left(r^0_{p,\lambda}(\Gamma)\alpha,\:\:\:\:
2\Ree(\lambda)+ \frac{2\rho}{p}-\max_{cusps \: [P]_\Gamma} \left[
2\delta_{\vp_{|\Gamma_P}} + \alpha \: \rk([P]_\Gamma) \right]\right)
\end{eqnarray*}  
(the maximum over the empty set is defined by convention as $-\infty$).
Then  
\begin{eqnarray*}
\Fam_\Gamma^{st}(\Lambda_\Gamma,\sigma_\lambda,\vp)&\subset&
H^{p,<r_{p,\lambda}(\Gamma)}(\partial X,V(\sigma_\lambda,\vp))\\
\Cusp_\Gamma(\sigma_\lambda,\vp)&\subset&
H^{p,<r^0_{p,\lambda}(\Gamma)}(\partial X,V(\sigma_\lambda,\vp)) 
\hspace{1cm} \mbox{if $\lambda\not\in I_\aaaa$}\ . \end{eqnarray*}
\end{theorem}  \underline{Remarks :} 
\begin{enumerate}
\item
If $\delta_\Gamma<\rho$, then it is
clear that $\phi\in  {}^\Gamma C^{-\infty}(\partial X, V(\sigma_\lambda,\vp))$
with $\supp(\phi)\subset \Lambda_\Gamma$  cannot belong to $H^{p,0}(\partial
X, V(\sigma_\lambda,\vp))$. Indeed, in this case the limit set
$\Lambda_\Gamma$ has trivial Lebesgue measure and cannot be the support of any
non-trivial $L^p$-function. We conclude
\begin{kor}
If $r_{1,\lambda}(\Gamma)>0$ and $\delta_\Gamma<\rho$, then
$$\Fam_\Gamma^{st}(\Lambda_\Gamma,\sigma_\lambda,\vp)=\{0\}\ .$$
\end{kor}
This is the vanishing result obtained in \cite{bunkeolbrich99}, Prop. 4.22,
resp. \cite{bunkeolbrich011}.
Note that the Corollary  does {\em not follow logically}  from the present
paper because it is {\em used} in the argument in order to deal with
positive  $r_{p,\lambda}(\Gamma)$.
\item
Assume that $\Gamma$ has finite covolume. In this case
$\rho=\delta_\Gamma$. If we take $\sigma=1$, $\vp=1$ and $\lambda=\rho$,
then $r_{p,\lambda}(\Gamma)=0$.  However, the constant function on $\partial X$
is a {\em smooth} $\Gamma$-invariant section of $V(1_\rho)$ which is strongly 
supported on the limit set.    \item
Assume that $\Gamma$ is convex cocompact, $\sigma=1$, $\vp=1$, and
$\lambda=\delta_\Gamma$. Then the {\em Patterson-Sullivan measure}  generates
$\Fam_\Gamma^{st}(\Lambda_\Gamma,1_\lambda,1)$.
A measure belongs to $H^{p,<-\dim(\partial X)\frac{1}{q}}$ for all $p\in
(1,\infty)$. In contrast, our estimate gives
$r_{p,\lambda}(\Gamma)\alpha = -(\dim(X)-\dim_H(\Lambda)) \frac{1}{q}
> -\frac{1}{q} \dim(\partial X)\alpha$ except if $\Gamma$ is elementary
hyperbolic. Thus the Patterson-Sullivan measure is more regular than a general
measure. It has the regularity of a smooth measure on $\Lambda$ if the latter
would be a smooth submanifold.
\item It follows from lower regularity bounds obtained by W. Schmid in special
cases that our estimate is essentially optimal. We refer to Subsection
\ref{ssc} for more details. 
\item
Probably, the condition $\lambda\not\in
I_\aaaa$ in the cusp form case is not necessary.

\end{enumerate}

Let us now explain the modifications in the case that $\Gamma$ has torsion. In
this case there is a torsion-free subgroup $\Gamma^\prime\subset \Gamma$ of
finite index which is still geometrically finite. Moreover, we have
$\Lambda_\Gamma=\Lambda_{\Gamma^\prime}$ and
$\delta_\Gamma=\delta_{\Gamma^\prime}$.  
We define
$$
\Fam_{\Gamma}^{st}(\Lambda_\Gamma,\sigma_\lambda,\vp):=
\Fam^{st}_{\Gamma^\prime}(\Lambda_{\Gamma^\prime},\sigma_\lambda,\vp)\cap
{}^\Gamma C^{-\infty}(\partial X,V(\sigma_\lambda,\vp))$$
(this definition does not depend on the choice of $\Gamma^\prime$). 
Cusp forms are characterized as in the torsion-free case. Then
Theorem \ref{main} holds true for $\Gamma$ with torsion elements.

\subsection{Idea of the proof}

Let us demonstrate the main ideas leading to the proof of Theorem \ref{main}
in the special case that $\delta_\vp=0$, $\sigma=1$,
$\Ree(\lambda)<\min(-\delta_\Gamma,0)$, $\lambda\not\in I_\aaaa$, and that
$\Gamma$ has no cusps of full rank and is not cocompact.

Let $\phi\in \Fam_\Gamma(\Lambda_\Gamma,1_\lambda,\vp)$. We apply the
Knapp-Stein intertwining operator (see Subsection \ref{reex}) and obtain
$\psi :=\hat J^w_{1_\lambda,\vp} \phi \in {}^\Gamma C^{-\infty}(\partial
X,V(1_{-\lambda},\vp))$. The fact that $\phi$ is strongly supported on the
limit set implies that $\psi$ is smooth on $\Omega_\Gamma$ and has a controlled
growth along the cusps. This will be formalized by saying that
$\psi_{|\Omega_\Gamma}$ belongs to the space $B_\Gamma(1_{-\lambda},\vp)$
which will be introduced in Subsection \ref{tpp}.  Using the convergence of the
Poincar\'e series $\sum_{\gamma\in\Gamma}
a_\gamma^{-(\delta_\Gamma+\rho+\epsilon)}$, $0<\epsilon\in\aaaa^*$, the
$\Gamma$-invariance of $\psi_{|\Omega_\Gamma}$, and the controlled growth
along the cusps we show (Theorem \ref{hypsum}) that $\psi\in
H^{p,<r_{p,\lambda}(\Gamma)-2\frac{\Ree(\lambda)}{\alpha}}$. Since
$\lambda\not\in I_\aaaa$ we recover $\phi$ from $\psi$ by applying a multiple
of $\hat J^{w^{-1}}_{1_{-\lambda},\vp}$. This operator decreases the
regularity by $-2\frac{\Ree(\lambda)}{\alpha}$ (see Subsection \ref{zc}) so
that we obtain the desired result $\phi\in H^{p,<r_{p,\lambda}(\Gamma)}$. If
$\phi$ is a cusp form, then we get a stronger estimate of $\psi$ along the
cusps leading to $\psi\in
H^{p,<r^0_{p,\lambda}(\Gamma)-2\frac{\Ree(\lambda)}{\alpha}}$ and $\phi\in
H^{p,<r^0_{p,\lambda}(\Gamma)}$.

\subsection{Specializations}
In this subsection we discuss specializations of the main theorem. 
In our main theorem we quantify the regularity of distributions using the
two-parameter family of spaces $H^{p,r}$.  
Let us indicate another way to measure the regularity of a
distributions.

For $0< s\in\R$ , $s=[s]+s^\prime$, $[s]\in\nat_0$, $s^\prime\in
[0,1)$ let  $C^s(\partial X,V(\sigma_\lambda,\vp))$ be the
space of sections which have  H\"older continuous (with exponent $s^\prime$)
derivatives up to order $[s]$. Furthermore, let $ C^{<s}:=\bigcap_{t<s}C^{t}$.
For negative $s$ we define $C^{s}$ (resp. $C^{<s}$) to be the space of
distributions which can locally be written as $N$-fold derivatives of
elements of $C^{s+N}$ (resp. $C^{<s+N}$), where $N$ is sufficiently large.  For
all $s\in \R$ and $p\in (1,\infty)$  have the embeddings
$H^{p,<\frac{n-1}{p}+s}\hookrightarrow C^{<s}$.
Since we deduce H\"older regularity from Sobolev regularity we cannot omit the
$<$-sign. In some cases a direct approach to H\"older regularity
may show $C^s$ instead of $C^{<s}$-regularity (see Subsection \ref{ssc}).

\subsubsection{$\Gamma$ cocompact}

Assume that $\Gamma$ is cocompact.
Then $\Lambda_\Gamma=\partial X$, $\delta_\Gamma=\rho$. 
Any invariant distribution
is stably deformable  by Theorem  \ref{ccalways} and  strongly
supported on the limit set.

\begin{theorem}
Assume that $\Gamma$ is cocompact, and let 
$$r \alpha:=  \Ree(\lambda) -\delta_\vp- \rho\ .$$ 
Then  
$${}^\Gamma C^{-\infty}(\partial X,V(\sigma_\lambda,\vp))\subset
H^{p,<r}(\partial X,V(\sigma_\lambda,\vp))$$ for any $p\in (1,\infty)$.
\end{theorem}
Using the freedom to choose $p$ arbitrary large we obtain the following
corollary. 
\begin{kor}
If $\Gamma$ is cocompact, 
 then 
$${}^\Gamma C^{-\infty}(\partial
X,V(\sigma_\lambda))\subset C^{<-\frac{n-1}{2}+\frac{\Ree(\lambda)}{\alpha}}\
.$$   \end{kor}

\subsubsection{$\Gamma$ convex cocompact}
Assume now that $\Gamma$ is convex cocompact.
If $0\not=\phi\in {}^\Gamma C^{-\infty}(\partial X,V(\sigma_\lambda,\vp))$
satisfies $\supp(\phi)\subset\Lambda_\Gamma$, then it is deformable by Theorem
\ref{ccalways} and strongly supported on the limit set.

We can thus apply
Theorem \ref{main} to the space 
$${}^\Gamma C^{-\infty}(\Lambda_\Gamma,V(\sigma_\lambda,\vp)):=\{\phi\in
{}^\Gamma C^{-\infty}(\partial X,V(\sigma_\lambda,\vp))\:|\:\supp(\phi)\subset
\Lambda_\Gamma\}\ .$$

\begin{theorem}
Let $$r\alpha:= \Ree(\lambda) -\delta_\vp -\rho + 
\frac{\rho-\delta_\Gamma}{p}\ .$$ Then  
$${}^\Gamma C^{-\infty}(\Lambda_\Gamma,V(\sigma_\lambda,\vp)) \subset
H^{p,<r}(\partial X,V(\sigma_\lambda,\vp))\ .$$
\end{theorem}

Using the freedom to choose $p$ arbitrary large   we obtain
\begin{kor}
If $\Gamma$ is convex cocompact, then  
$${}^\Gamma
C^{-\infty}(\Lambda_\Gamma,V(\sigma_\lambda,\vp))\subset
C^{<-\frac{n-1}{2}+\frac{\Ree(\lambda) -\delta_\vp}{\alpha}}\ .$$ \end{kor}

\subsubsection{$\Gamma$ has finite covolume}

Now assume that $\Gamma$ has finite covolume.
Let $\phi\in {}^\Gamma C^{-\infty}(\partial X,V(\sigma_\lambda))$ be a cusp
form and  
$\lambda\not\in I_\aaaa$.  
Applying Theorem \ref{main} and using that $\delta_\Gamma=\rho$ we obtain the
following result.

\begin{theorem}
Assume that $\Gamma$ has finite covolume and $\lambda\not\in I_\aaaa$, then
 $$\Cusp_\Gamma(\sigma_\lambda) \subset H^{p,<
-\frac{n-1}{2}+\frac{\Ree(\lambda)}{\alpha}}(\partial X,V(\sigma_\lambda))$$
for any $p\in (1,\infty)$. \end{theorem}

\begin{kor}
If $\Gamma$ has finite covolume and $\lambda\not\in I_\aaaa$, then 
 $$\Cusp_\Gamma(\sigma_\lambda) \subset C^{<
-\frac{n-1}{2}+\frac{\Ree(\lambda)}{\alpha}}\ .$$ \end{kor}

More general we have the following estimate.
\begin{theorem}
If $\Gamma$ has finite covolume,  then for $p\in (1,\infty)$ we have
 $$\Fam_\Gamma^{st}(\Lambda_\Gamma,\sigma_\lambda,\vp)\subset
  H^{p,<
 -\frac{n-1}{q}+2\frac{\Ree(\lambda)}{\alpha}}(\partial
X,V(\sigma_\lambda))\ .$$ \end{theorem}

\begin{kor}
If $\Gamma$ has finite covolume, then 
 $$\Fam_\Gamma^{st}(\Lambda_\Gamma,\sigma_\lambda,\vp)\subset
C^{<-(n-1)+2\frac{\Ree(\lambda)}{\alpha}}(\partial X,V(\sigma_\lambda,\vp))\
.$$ \end{kor}

\subsubsection{Residues of Eisenstein series}

For general geometrically finite $\Gamma$  a source of interesting
invariant distributions are boundary values of residues of Eisenstein series.
In the language of geometric scattering theory \cite{bunkeolbrich99} these are
invariant distributions obtained from the singular part of the extension map
$ext^\Gamma$. We refer to Subsection \ref{reex} for more details.

\begin{ddd}\label{singext}
Let $\Ext^{sing}_\Gamma(\sigma_\lambda,\vp)\subset {}^\Gamma
C^{-\infty}(\partial X,V(\sigma_\lambda,\vp))$ be the subspace of invariant
distributions generated by the singular part of $ext^\Gamma$.
\end{ddd}

Essentially by definition the elements of
$\Ext^{sing}_\Gamma(\sigma_\lambda,\vp)$ belong to ${}^\Gamma
C^{-\infty}(\partial X,V(\sigma_\lambda,\vp))$ and are deformable. Moreover,
the relation $res^\Gamma\circ ext^\Gamma=\id$ implies that the elements of
$\Ext^{sing}_\Gamma(\sigma_\lambda,\vp)$ are strongly supported on the limit
set. We thus have $\Ext^{sing}_\Gamma(\sigma_\lambda,\vp)\subset
\Fam_\Gamma(\Lambda_\Gamma,\sigma_\lambda,\vp)$. We show in Section \ref{reex}
that these spaces are in fact equal.

Note that  for $\Ree(\lambda)>\delta_\Gamma+\delta_\vp$ the
extension $ext^\Gamma$ is regular
and thus $\Ext^{sing}_\Gamma(\sigma_\lambda,\vp)=0$.
Theorem \ref{main} gives the
following result.

\begin{theorem}
We have
 $$ \Ext^{sing}_\Gamma(\sigma_\lambda,\vp) \subset
H^{p,<r_{p,\lambda}(\Gamma)}(\partial X,V(\sigma_\lambda,\vp))\ .$$
\end{theorem}
Choosing $p$ arbitrary large we obtain the following consequence. 
\begin{kor}
If $\Gamma$ is geometrically finite, then
$$\Ext^{sing}_\Gamma(\sigma_\lambda,\vp) \subset
C^{<r}\ ,$$
where 
$$r:=\min\left(-\frac{n-1}{2}+\frac{\Ree(\lambda)-\delta_\vp}{\alpha},\:
\:2\frac{\Ree(\lambda)}{\alpha}-\max_{cusps\:
[P]_\Gamma}\left[2\frac{\delta_{\vp_{|P_\Gamma}}}{\alpha}+\rk([P]_\Gamma)\right]\right)\
.$$  \end{kor}

\subsection{Related results}

\subsubsection{Results of W. Schmid}\label{ssc}

In \cite{schmid00} W. Schmid considers a finitely generated discrete
subgroup $\Gamma$ of $G=SL(2,\R)$. Note that $G\cong Spin(1,2)$.
If $\phi\in {}^\Gamma C^{-\infty}(\partial X,V(\sigma_\lambda))$, then it
defines a $G$-equivariant map 
$$c_\phi: C^\infty(\partial X,V(\tilde \sigma_{-\lambda}))\rightarrow
C^\infty(\Gamma\backslash G)$$ by
$c_\phi(f)(g):=\langle \phi,\pi^{\tilde\sigma_\lambda}(g) f\rangle$.
Let us call $\phi$ {\em  square integrable}, iff $c_\phi(f)\in
L^2(\Gamma\backslash G)$ for all $f\in C^\infty(\partial X,V(\tilde
\sigma_{-\lambda}))$. Schmid calls a square integrable $\phi$ (or rather the
map $c_\phi$) {\em cuspidal}, if $c_\phi(f)$ is a bounded function
for all $K$-finite $f$, where $K$ is some maximal compact subgroup of $G$.
\begin {theorem}[Schmid]
Let $\phi\in {}^\Gamma C^{-\infty}(\partial X,V(\sigma_\lambda))$ be
square integrable.
\begin{enumerate}
\item
If $\phi$ is cuspidal and $\frac{\lambda}{\alpha}\not\in \Z+\frac12$, then
$\phi$ is of H\"older regularity $C^{\frac{\Ree(\lambda)-\rho}{\alpha}}$.
\item 
If $\frac{\lambda}{\alpha}\in \Z+\frac12$, then $\phi$ is cuspidal
and belongs to $C^{<\frac{\lambda-\rho}{\alpha}}$.
\item
If $\phi$ is not cuspidal, then $\phi\in C^{<2(\lambda-\rho)}$ if $\lambda\in
(0,\rho)$, and $\phi\in C^{<-1}$, if $\lambda\in (-\rho,0)$.
\end{enumerate}
\end{theorem}
Note that in the two-dimensional case  a finitely generated discrete group
is  geometrically finite, so the classes of discrete
subgroups considered by Schmid and in the present paper coincide.
On the one hand, being a square integrable invariant distribution
is a strong restriction on $\phi$. On the other hand there are square
integrable distributions which are not strongly supported on the limit set.
If $\phi$ is cuspidal and $\lambda\not\in I_\aaaa$,  or if $\phi$ is stably
deformable, then our estimate applies to $\phi$ and reproduces Schmid's
estimate. In fact,  in the cuspidal case  the result Schmid is slightly
stronger than ours since he shows $\phi\in C^{r}$, while we show $\phi\in
C^{<r}$ for the appropriate $r$.

In the  case $\lambda\in (0,\rho)$ Schmid shows that his lower bounds of
regularity are optimal in both, the cuspidal and the general case. 
This implies that our result is essentially  optimal among estimates which can
be stated using the same data and are compatible with embedding and twisting.

{\em We thank W. Schmid for showing us a counterexample  to an estimate
claimed in a previous version of this paper.}

\subsubsection{Results of J. Lott}

In \cite{lott98} Lott considers discrete subgroups $\Gamma$ of $SO(1,n)_0$.
Let $\sigma^k:=\Lambda^k\bar \naaa\in\hat M$ and $\lambda_k:=\rho- k
\alpha$. Note that $V(\sigma^k_{\lambda_k})=\Lambda^kT^*\partial X$.
Let $\gaaa=\kaaa\oplus\paaa$ be the Cartan decomposition of $\gaaa$.
Then $\Lambda^k\paaa\in\hat K$ is the minimal $K$-type of
$\pi^{\tilde\sigma^k_{-\lambda_k}}$. Fix $0\not=f_k\in C^{\infty}(\partial
X,V(\tilde\sigma^k_{-\lambda_k}))(\Lambda^k\paaa)$.
The main result of Lott is the following theorem.
\begin{theorem}[Lott]
If $\phi\in C^{-\omega}(\partial X,V(\sigma^k_{\lambda_k}))$ is
exact (as a differential form) and such that $c_\phi(f_k)\in C^\infty(G)$ is
bounded, then $\phi\in H^{2,<{-k}}$. \end{theorem}
Among the variety of results proved in Lott's paper the following can be compared with those of the present paper.
\begin{theorem}[Lott]
Let $\Gamma$ be convex cocompact.
\begin{enumerate}
\item
If $k\in [1,\frac{n-1}{2})$, and $\phi\in {}^\Gamma
C^{-\infty}(\Lambda_\Gamma,V(\sigma^k_{\lambda_k}))$, then $\phi\in
H^{2,<-k}$. \item If $\Gamma$ is not cocompact and $\phi\in {}^\Gamma
C^{-\infty}(\Lambda_\Gamma,V(\sigma^1_{\lambda_1}))$, then $\phi\in
H^{2,-1}$. \end{enumerate} \end{theorem}
This can directly be compared with our result, which gives 
$\phi\in H^{2,<-k+\frac{\rho-\delta_\Gamma}{2}}$.
Thus if $\Gamma$ is not cocompact, then our result improves Lott's estimate.

\subsubsection{The result of J. Bernstein and A. Reznikov}

In \cite{bernsteinreznikov98} J. Bernstein and A. Reznikov
consider $G=SL(2,\R)\cong Spin(1,2)$, and a discrete subgroup
of finite covolume. In Sec. 1.4 they "prove" the following result.
\begin{theorem}[Bernstein-Reznikov]
If $\Ree(\lambda)=0$ and $\phi\in  {}^\Gamma C^{-\infty}(\partial X,V(\sigma_\lambda))$ is square integrable, then $\phi\in (H^{1,1/2+\epsilon})^*$ for any $\epsilon>0$.
\end{theorem}
Our estimate applies if $\lambda\not=0$, because in this case
$\phi$ is a cusp form and $\lambda\not\in I_\aaaa$.
We show 
$\phi\in (H^{p,1/2+\epsilon})^*$ for any $\epsilon>0$ and $p\in (1,\infty)$.
This implies Bernstein-Reznikov's result.

\section{Elements of geometric scattering theory}
\subsection{The space $B_\Gamma(\sigma_\lambda,\vp)$}\label{tpp}

The spaces $B_\Gamma(\sigma_\lambda,\vp)$ were introduced in
\cite{bunkeolbrich99}. In this subsection we recall some of their basic
properties. First we consider the case that all cusps of $\Gamma$ have smaller
rank. 

Recall that $\Gamma$ acts properly discontinuously on the complement
$\Omega_\Gamma=\partial X\setminus\Lambda_\Gamma$ of the limit set
$\Lambda_\Gamma$ with quotient $B_\Gamma$. Furthermore, we have a quotient
bundle $V_{B_\Gamma}(\sigma_\lambda,\vp):=\Gamma\backslash
V(\sigma_\lambda,\vp)_{|\Omega_\Gamma}$. The bundles
$V_{B_\Gamma}(\sigma_\lambda,\vp)$ depend holomorphically on $\lambda\in\aca$.
 We consider the family of
spaces
$(C^{\infty}(B_\Gamma,V_{B_\Gamma}(\sigma_\lambda,\vp)))_{\lambda\in\aca}$ as a
(trivial) bundle of Fr\'echet spaces over $\aca$.

For $\Ree(\lambda)<-\delta_\Gamma-\delta_\vp$ we can define the push-down map
$$\pi^\Gamma_*:C^\infty(\partial X,V(\sigma_\lambda,\vp))\rightarrow
C^\infty(B_\Gamma,V_{B_\Gamma}(\sigma_\lambda,\vp))$$
by $\pi_*^\Gamma(f):=\sum_{\gamma\in\Gamma}\pi^{\sigma_\lambda,\vp}(\gamma)
f_{|\Omega_\Gamma}$. Indeed, this sum converges and defines a smooth
$\Gamma$-invariant section of $V(\sigma_\lambda,\vp)_{|\Omega_\Gamma}$,
and thus an element of $C^\infty(B_\Gamma,V_{B_\Gamma}(\sigma_\lambda,\vp))$.
One of the main results of  \cite{bunkeolbrich99} is that $\pi^\Gamma_*$ has
a meromorphic continuation to all of $\aca$. If $\Gamma$ is not
convex-cocompact, then $\pi_*^\Gamma$ is not surjective.
Roughly speaking, the space $B_\Gamma(\sigma_\lambda,\vp)$ is
defined to be the range of $\pi^\Gamma_*$. 
\begin{ddd}
If $\Gamma$ has no cusps of full rank, then we define
$(B_\Gamma(\sigma_\lambda,\vp))_{\lambda\in\aca}$ to be the minimal holomorphic subbundle of
$(C^\infty(B_\Gamma,V_{B_\Gamma}(\sigma_\lambda,\vp)))_{\lambda\in\aca}$ over
which $\pi^\Gamma_*$ factors.
\end{ddd}
This defines $B_\Gamma(\sigma_\lambda,\vp)$ as a vector space. 
In fact, $f\in B_\Gamma(\sigma_\lambda,\vp)$ iff there exists a meromorphic
family $(F_\mu)_\mu$, $F_\mu\in C^\infty(\partial X,V(\sigma_\mu,\vp))$, such
that $(\pi^\Gamma_* F_\mu)_{\mu=\lambda}=f$.

As a
topological vector space it will be equipped with a topology which is stronger
than  the topology induced from the embedding
$B_\Gamma(\sigma_\lambda,\vp)\rightarrow
C^\infty(B_\Gamma,V_{B_\Gamma}(\sigma_\lambda,\vp))$.
With this stronger topology we cannot show that
$(B_\Gamma(\sigma_\lambda,\vp))_{\lambda\in\aca}$
is a locally trivial bundle of Fr\'echet spaces. As explained in
\cite{bunkeolbrich99} we can consider it as a projective limit (see below)
of locally trivial bundles of Fr\'echet spaces, and it therefore still makes
sense to speak of holomorphic families $(f_\mu)_{\mu\in\aca}$, $f_\mu\in
B_\Gamma(\sigma_\mu,\vp)$.

Note that $C_c^\infty(B_\Gamma,V_{B_\Gamma}(\sigma_\lambda,\vp))\subset
B_\Gamma(\sigma_\lambda,\vp)$, so that the elements of
$B_\Gamma(\sigma_\lambda,\vp)$ are distinguished among all smooth sections by
their behaviour at the ends of $B_\Gamma$. In order to describe this behaviour
in detail we first assume that $\Gamma\subset P$ for some $\Gamma$-cuspidal
parabolic subgroup $P\subset G$. Then we have $\Omega_\Gamma=\partial
X\setminus \infty_P$. We choose a maximal compact subgroup $K$ and a Langlands
decomposition $P=MAN$, $M\subset K$. Furthermore, let $w\in K$ be a
representative of the non-trivial element of the Weyl group $N_K(M)/M$.
Then we can parametrize $\Omega_\Gamma=NwP\subset G/P$.
A function  $f\in C^\infty(\Omega_\Gamma,V(\sigma_\lambda,\vp))$
gives rise to a smooth function on $N$ with values in
$V_{\sigma}\otimes V_\vp$ by
$$N\ni x\rightarrow f(xw)\in V_{\sigma}\otimes V_\vp\ .$$
For $X\in\cU(\naaa)$, $W\subset N$, and $\beta\in\aaaa^*$ we define
\begin{eqnarray}
q^\Gamma_{W,X,\beta}(f)&:=&\sup_{x\in W} \sup_{a\in A_+}
a^{\beta-2(\Ree(\lambda)-\rho^\Gamma)} |\vp(a^{-1})
f(x^a Xw)|\label{schn}\\
q_{W,X}(f)&=&\sup_{x\in W} |f(xX)|\ .
\nonumber \end{eqnarray}
Note that in this context we extend $\vp$ to $AP_\Gamma$ (which is possible
since we assume that $\vp$ is admissible) such that the minimal $A$-weight is
$0$ and the highest $A$-weight is not greater than $2\delta_\vp$.  

The family of seminorms $q_{W,X}$ for all compact $W\subset N$ and
$X\in\cU(\naaa)$ defines the topology of
$C^\infty(\Omega_\Gamma,V(\sigma_\lambda,\vp))$. \begin{ddd}
Let $S_{\Gamma,\beta,\kappa}(\sigma_\lambda,\vp)$ be the closure
of $C_c^\infty(B_\Gamma,V_{B_\Gamma}(\sigma_\lambda,\vp))$ with respect to
Banach space topology defined by the seminorms
$q_{W,X}$, $\deg(X)\le \kappa$, $W\subset \Omega_\Gamma$ compact, 
and $q^\Gamma_{W,X,\beta}$, $\deg(X)\le \kappa$, $W\subset N\setminus N_\Gamma$
compact. Furthermore, we define the Fr\'echet spaces 
$S_{\Gamma,\beta}(\sigma_\lambda,\vp):=\bigcap_{\kappa\in\aaaa^*}
S_{\Gamma,\beta,\kappa}(\sigma_\lambda,\vp)$ and
$S_{\Gamma}(\sigma_\lambda,\vp):=\bigcap_{\beta\in\aaaa^*}S_{\Gamma,\beta}(\sigma_\lambda,\vp)$. 
\end{ddd}
The elements of $B_\Gamma(\sigma_\lambda,\vp)$ will now be characterized by
the property that they have a certain asymptotic expansion near $\infty_P$
with remainder in $S_{\Gamma}(\sigma_\lambda,\vp)$.
In order to describe this asymptotic expansion we introduce the following
space of homogeneous functions.
\begin{ddd}
For $\beta\in\aaaa^*$ we define $A_{P_\Gamma,\beta}(\sigma_\lambda,\vp)\subset
C^\infty(\Omega_\Gamma\setminus N_\Gamma wP, V(\sigma_\lambda,\vp))$ to be the
closed subspace of $P_\Gamma$-invariant sections satisfying
$$\vp(a)^{-1} f(x^aw)=a^{2(\lambda-\rho^\Gamma)-\beta} f(xw)\ .$$
\end{ddd}
In \cite{bunkeolbrich99} we have shown that if $f\in 
B_\Gamma(\sigma_\lambda,\vp)$,
then it has an asymptotic expansion
$$\vp(a^{-1})f(x^aw)\stackrel{a\to\infty}{\sim} \sum_{k\in\nat_0}
a^{2(\lambda-\rho^\Gamma)-k\alpha}f_k(x)\ , $$
where $f_k\in  A_{P_\Gamma,k\alpha}(\sigma_\lambda,\vp)$.
\begin{ddd}
We define $R_{\Gamma}(\sigma_\lambda,\vp)^k\subset
A_{P_\Gamma,k\alpha}(\sigma_\lambda,\vp)$ to be the subspace generated
by the $f_k$ for all $f\in B_\Gamma(\sigma_\lambda,\vp)$.
\end{ddd}
The space $R_\Gamma(\sigma_\lambda,\vp)^k$ is finite-dimensional.
We further set $R_{\Gamma,k}(\sigma_\lambda,\vp):=\bigoplus_{l\le
k}R_\Gamma(\sigma_\lambda,\vp)^l$, 
$B_{\Gamma,k}(\sigma_\lambda,\vp):=S_{\Gamma,k\alpha}(\sigma_\lambda,\vp)+B_{\Gamma}(\sigma_\lambda,\vp)$, such that
we have an exact sequence
$$0\rightarrow S_{\Gamma,k\alpha}(\sigma_\lambda,\vp)\rightarrow
B_{\Gamma,k}(\sigma_\lambda,\vp)\stackrel{AS}{\rightarrow}R_{\Gamma,k}(\sigma_\lambda,\vp)\rightarrow 0\ .$$
We can construct a split $L:R_{\Gamma,k}(\sigma_\lambda,\vp)\rightarrow
B_{\Gamma,k}(\sigma_\lambda,\vp)$
of this exact sequence using a suitable $P_\Gamma$-invariant cut-off function.
So we have
$B_{\Gamma,k}(\sigma_\lambda,\vp)=S_{\Gamma,k\alpha}(\sigma_\lambda,\vp)\oplus
L(R_{\Gamma,k}(\sigma_\lambda,\vp))$, and this decomposition defines a
Fr\'echet space topology on $B_{\Gamma,k}(\sigma_\lambda,\vp)$. Note that
$B_\Gamma(\sigma_\lambda,\vp)=\bigcap_{k\in\nat}
B_{\Gamma,k}(\sigma_\lambda,\vp)$.

We need the following seminorms on $B_{\Gamma}(\sigma_\lambda,\vp)$.

\begin{ddd}\label{pdef}
If $W\subset N\setminus N_{\Gamma_P}$ is compact, $k\in\nat_0$, and $X\in
\cU(\naaa)$ is homogeneous of degree $\deg(X)\in\aaaa^*$, then we define
$p^{\Gamma}_{W,X}(f)$ by $$p^{\Gamma}_{W,X}(f):=\sup_{a\in
A_+}\sup_{x\in W} a^{\deg(X)-2(\Ree(\lambda)-\rho^{\Gamma})}
|\vp(a^{-1}) f(x^aXw)|\ .$$ 
\end{ddd}
\begin{lem}
$p^{\Gamma}_{W,X}$ is a continuous seminorm on $
B_{\Gamma}(\sigma_\lambda,\vp)$.
\end{lem}
\proof
Fix $k\in\nat$ such that $k\alpha>\deg(X)$.
Let $f\in B_{\Gamma,k}(\sigma_\lambda,\vp)$.
We define $f_l$ such
that  $AS(f)=\oplus_{l\le k} f_l$. Furthermore, we set 
$g:=f-L AS(f)\in S_{\Gamma,k\alpha}(\sigma_\lambda,\vp)$. Note
that $f_l$ and $g$ depend continuously on $f$.
The restriction of $p^{\Gamma}_{W,X}$ to the space
$S_{\Gamma,k\alpha}(\sigma_\lambda,\vp)$
belongs to the family of seminorms which define the topology of
$S_{\Gamma,k\alpha}(\sigma_\lambda,\vp)$.
So it remains to show that
$A_{P_\Gamma,k\alpha}(\sigma_\lambda,\vp)\ni f_l\rightarrow
p^{\Gamma}_{W,X}(L f_l)$ is  continuous.

Let $\chi_{\Gamma}$ be the $P_\Gamma$-invariant cut-off
function used to define $ L$. Furthermore, 
let $\Delta(X)=\sum_{j} Y_j\otimes Z_j$ be the coproduct.
Then we can write
\begin{eqnarray*}
\vp(a^{-1})L(f_l)(x^aXw)&=&\vp(a^{-1})\sum_j \chi(x^a Y_jw)
f_l(x^aZ_jw)\\ &=&\vp(a^{-1}) \sum_j \chi( x^a  Y_jw) f_l((x
Z^{a^{-1}}_j)^aw)\\ &=&\sum_j   
a^{2(\lambda-\rho^{\Gamma})-\deg(Z_j)-l\alpha}     \chi(x^aY_jw) f_l(x
Z_jw)\\ &=&a^{2(\lambda-\rho^{\Gamma})-\deg(X)-l\alpha}     \chi(x^aw)
f_l(x Xw) \\&&\hspace{1cm}+ \sum_{j,\deg(Y_j)>0}   \chi( x^a w Y_j) f_l((x
Z^{a^{-1}}_j)^aw)\ . \end{eqnarray*}
The summands for $\deg(Y_j)>0$ vanish for large $a\in A_+$ uniformly for $x\in
W$.
We therefore can majorize
$p^{\Gamma}_{W,X}( L f_l)$ by a continuous seminorm  of
$A_{P_\Gamma,l\alpha}(\sigma_\lambda,\vp)$.
\hB

We can now describe $B_\Gamma(\sigma_\lambda,\vp)$ for a general geometrically
finite torsion-free group $\Gamma$ without cusps of full rank. If $P\subset G$
is a $\Gamma$-cuspidal parabolic subgroup, then there is a representative
$B_P\subset B_\Gamma$ of an end of $B_\Gamma$ which is isomorphic to a
representative of the end of $B_{\Gamma_P}$. Let $\chi_P\in
C^\infty(B_\Gamma)$ be supported on $B_P$ such that it is identically one near
infinity. If $f\in C^\infty(B_\Gamma,V_{B_\Gamma}(\sigma_\lambda,\vp))$, then
we can consider $\chi_P f$ as a section $T_Pf\in
C^\infty(B_{\Gamma_P},V_{B_{\Gamma_P}}(\sigma_\lambda,\vp))$.
Then $f\in B_\Gamma(\sigma_\lambda,\vp)$ iff
$T_Pf\in B_{\Gamma_P}(\sigma_\lambda,\vp)$ for all $\Gamma$-cuspidal parabolic
subgroups $P\subset G$. 
\begin{ddd}
As a topological vector space
$B_\Gamma(\sigma_\lambda,\vp)$ is equipped with the smallest topology such that
the inclusion $B_\Gamma(\sigma_\lambda,\vp)\subset
C^\infty(B_\Gamma,V_{B_\Gamma}(\sigma_\lambda,\vp))$
and all the maps $T_P$ are continuous.
\end{ddd}
In a similar manner we define the topological vector spaces
$B_{\Gamma,k}(\sigma_\lambda,\vp)$ ($S_{\Gamma,\beta}(\sigma_\lambda,\vp)$) to
be the space of all $f\in C^\infty(B_\Gamma,V_{B_\Gamma}(\sigma_\lambda,\vp))$
such that $T_Pf\in B_{\Gamma_P,k}(\sigma_\lambda,\vp)$ ($T_Pf\in
S_{\Gamma_P,\beta}(\sigma_\lambda,\vp)$) for all $\Gamma$-cuspidal parabolic
subgroups $P\subset G$. Furthermore, we set
$R_{\Gamma,k}(\sigma_\lambda,\vp):=\bigoplus_{P\in \tilde \cP}
R_{\Gamma_P,k}(\sigma_\lambda,\vp)$, where
$\tilde\cP$ denotes a set of representatives of the $\Gamma$-conjugacy
classes of $\Gamma$-cuspidal parabolic subgroups ($\tilde \cP$ parametrizes
the ends of $B_\Gamma$). Then we have a split exact sequence
$$0\rightarrow S_{\Gamma,k\alpha}(\sigma_\lambda,\vp)\rightarrow
B_{\Gamma,k}(\sigma_\lambda,\vp)\stackrel{AS}{\rightarrow}R_{\Gamma,k}(\sigma_\lambda,\vp)\rightarrow 0\ ,$$
and we can write
\begin{equation}\label{summn}
B_{\Gamma,k}(\sigma_\lambda,\vp)=S_{\Gamma,k\alpha}(\sigma_\lambda,\vp)+L(R_{\Gamma,k}(\sigma_\lambda,\vp))\ .
\end{equation}

We now describe the necessary modification if $\Gamma$ has cusps of full rank.
In this case the ends of $B_\Gamma$ are parametrized by the set of
$\Gamma$-conjugacy classes of $\Gamma$-cuspidal subgroups $\cP^<$ corresponding
to cusps of smaller rank. The same construction as above gives spaces
$B_{\Gamma}(\sigma_\lambda,\vp)_1$, $B_{\Gamma,k}(\sigma_\lambda,\vp)_1$. 
The space  $B_{\Gamma,k}(\sigma_\lambda,\vp)$ will be the direct  sum of
$B_{\Gamma,k}(\sigma_\lambda,\vp)_1$ and a finite-dimensional vector
space $B_{\Gamma}(\sigma_\lambda,\vp)_2:=\oplus_{P\in\tilde \cP^{max}}
R_{\Gamma_P}(\sigma_\lambda,\vp)$, where $\tilde \cP^{max}$ is a set of
representatives of $\Gamma$-conjugacy classes $\Gamma$-cuspidal parabolic
subgroups corresponding to cusps of full rank.

The space $R_{\Gamma_P}(\sigma_\lambda,\vp)$ is the fibre of a trivial
holomorphic vector bundle as follows.
Let ${}^{\Gamma_P}\cE_{\infty_P}(\sigma,\vp)$ denote
the sheaf of holomorphic families
$f_\nu\in{}^{\Gamma_P}C^{-\infty}(\partial X,V(\sigma_\nu,\vp))$ with
$\supp(f)=\infty_P$. Since ${}^{\Gamma_P}\cE_{\infty_P}(\sigma,\vp)$ is torsion-free it is
the space of sections of a unique holomorphic vector bundle
${}^{\Gamma_P}E_{\infty_P}(\sigma,\vp)$ over $\aca$. By ${}^{\Gamma_P}E_{\infty_P}(\sigma_\lambda,\vp)$
we denote the fibre of ${}^{\Gamma_P}E_{\infty_P}(\sigma,\vp)$ at $\lambda\in\aca$. Then we
define $R_{\Gamma_P}(\sigma_\lambda,\vp):=
{}^{\Gamma_P}E_{\infty_P}(\tilde\sigma_{-\lambda},\tilde\vp)^*$. 

The push-down $\pi^\Gamma_*$ now decomposes as $(\pi^\Gamma_*)_1\oplus
(\pi^\Gamma_*)_2$, where $(\pi^\Gamma_*)_1$ is the average as in the case
without cusps of full rank. The second component further decomposes
as $(\pi^\Gamma_*)_2=\oplus_{P\in \tilde\cP^{max}}[\pi^\Gamma_* ]_{P} $,
where $[\pi^\Gamma_* ]_{P}$ has values in $R_{\Gamma_P}(\sigma_\lambda,\vp)$.
For $\Ree(\lambda)<-\delta_\Gamma-\delta_\vp$ it is defined by the condition
$$\langle \phi,[\pi^\Gamma_*]_P(f)\rangle =
\sum_{[g]\in\Gamma/\Gamma_P}\langle
\pi^{\tilde\sigma_{-\lambda},\tilde\vp}(g)\phi,f \rangle$$ for all $\phi\in
{}^{\Gamma_P}E_{\infty_P}(\tilde\sigma_{-\lambda},\tilde\vp)$, and for general $\lambda$
by meromorphic continuation.

The family of spaces $(B_{\Gamma,k}(\sigma_\lambda,\vp))_{\lambda\in\aca}$
forms a trivial holomorphic bundle of Fr\'echet spaces over $\aca$.
In the present paper we do not need an explicit trivialization, and we
therefore omit its description, which can be found in \cite{bunkeolbrich99}.

\subsection{Restriction, extension, and the scattering matrix}\label{reex}

In the preceding subsection we have  defined the Fr\'echet spaces
$B_\Gamma(\sigma_\lambda,\vp)$.  
\begin{ddd}
We define
$D_\Gamma(\sigma_\lambda,\vp):=B_\Gamma(\tilde\sigma_{-\lambda},\tilde\vp)^*$.
\end{ddd}
Furthermore, we introduce
$$D_{\Gamma,k}(\sigma_\lambda,\vp):=B_{\Gamma,k}(\tilde\sigma_{-\lambda},\tilde\vp)^*\ .$$
The family $(D_{\Gamma,k}(\sigma_\lambda,\vp))_{\lambda\in\aca}$ is a trivial
holomorphic bundle of dual Fr\'echet  spaces. Moreover, we have
$ \displaystyle D_\Gamma(\sigma_\lambda,\vp)=\lim_{\stackrel{\rightarrow}{k}}
D_{\Gamma,k}(\sigma_\lambda,\vp)$.
Therefore we can speak of holomorphic families $(\phi_\mu)_{\mu\in\aca}$,
$\phi_\mu\in D_\Gamma(\sigma_\mu,\vp)$. Locally, $\phi_\mu$ is a holomorphic
family in $D_{\Gamma,k}(\sigma_\mu,\vp)$ for some fixed $k$.

\begin{ddd}
We define the extension map
$ext^\Gamma:D_\Gamma(\sigma_\lambda,\vp)\rightarrow C^{-\infty}(\partial
X,V(\sigma_\lambda,\vp))$ to be the adjoint of $\pi^\Gamma_*$.
\end{ddd}
If $(\phi_\mu)_{\mu\in\aca}$,
$\phi_\mu\in D_\Gamma(\sigma_\mu,\vp)$, is a meromorphic family, then
$ext^\Gamma\phi_\mu\in C^{-\infty}(\partial
X,V(\sigma_\lambda,\vp))$ is a meromorphic family of $\Gamma$-invariant
distributions. 
\begin{ddd} By 
$\Ext_\Gamma(\sigma_\lambda,\vp)$ we denote the subspace
of ${}^\Gamma C^{-\infty}(\partial
X,V(\sigma_\lambda,\vp))$ consisting of evaluations at $\lambda$ of
germs families $(ext^\Gamma \phi_\mu)_\mu$, where $(\phi_\mu)_{\mu}$,
$\phi_\mu\in D_\Gamma(\sigma_\mu,\vp)$, is a germ of a meromorphic family.
By $\Ext_\Gamma^0(\sigma_\lambda,\vp)\subset
\Ext_\Gamma(\sigma_\lambda,\vp)$ we denote the subspace of those
evaluations, where $\phi_\lambda$ is regular.
\end{ddd}
In \cite{bunkeolbrich99} we have shown that for  generic $\lambda\in\aca$ we
have  
$$\Ext_\Gamma^0(\sigma_\lambda,\vp)= 
\Ext_\Gamma(\sigma_\lambda,\vp)={}^\Gamma C^{-\infty}(\partial
X,V(\sigma_\lambda,\vp))\ .$$
\begin{ddd}
If $(\psi_\mu)_\mu$, $\psi_\mu  \in {}^\Gamma C^{-\infty}(\partial
X,V(\sigma_\mu,\vp))$, is a germ of a meromorphic family at $\lambda$, then
we define $(res^\Gamma(\psi_\mu))_\mu$, $res^\Gamma(\psi_\mu)\in 
D_\Gamma(\sigma_\mu,\vp)$, to be the  unique germ of a meromorphic meromorphic
family such that  $ext^\Gamma (res^\Gamma(\psi_\mu))=\psi_\mu$.
\end{ddd}
Note that if $\psi\in \Ext^0_\Gamma(\sigma_\lambda,\vp)$, then we can define
$res^\Gamma \psi\in D_\Gamma(\sigma_\lambda,\vp)$.

Recall Definition \ref{singext} of
$\Ext^{sing}_\Gamma(\sigma_\lambda,\vp)$. 
\begin{lem}\label{eqqu}
We have the equality
$$\Ext^{sing}_\Gamma(\sigma_\lambda,\vp)=\Fam_\Gamma(\Lambda_\Gamma,\sigma_\lambda,\vp)\ .$$
\end{lem}
\proof
The inclusion
$\Ext^{sing}_\Gamma(\sigma_\lambda,\vp)\subset\Fam_\Gamma(\Lambda_\Gamma,\sigma_\lambda,\vp)$
follows from the identity $res^\Gamma\circ ext^\Gamma=\id$, which holds when
both sides are applied to germs of meromorphic families. The inclusion
$\Fam_\Gamma(\Lambda_\Gamma,\sigma_\lambda,\vp)\subset
\Ext^{sing}_\Gamma(\sigma_\lambda,\vp)$ follows from the similar identity
$ext^\Gamma\circ res^\Gamma=\id$, which also holds after application to  germs
of meromorphic families. \hB   

Next we  show that the condition "strongly supported on the limit set" is
stable under twisting and embedding.
\begin{lem}\label{stablesupp}
If $i_*:C^{-\infty}(\partial X,V(\sigma_\lambda,\vp))\rightarrow
C^{-\infty}(\partial X^\prime,V(\sigma^\prime_{\lambda^\prime},\vp^\prime))$ is
induced by some twisting or embedding data, then
$$i_*(\Fam_\Gamma(\Lambda_\Gamma,\sigma_\lambda,\vp))\subset
\Fam_\Gamma(\Lambda_\Gamma,\sigma^\prime_{\lambda^\prime},\vp^\prime)\ .$$
\end{lem}
\proof
We show
$i_*(\Ext^{sing}_\Gamma(\sigma_\lambda,\vp))\subset
\Ext^{sing}_\Gamma(\sigma^\prime_{\lambda^\prime},\vp^\prime)$
and then apply Lemma \ref{eqqu}.
But the latter inclusion follows immediately from the identity (see
\cite{bunkeolbrich99}, Sec. 2.3)  $i_*\circ ext^\Gamma= ext^\Gamma \circ
i_*^\Gamma$, where the maps $ext^\Gamma$ on the two sides of course act on
different spaces,
and $i_*^\Gamma:D_\Gamma(\sigma_\lambda,\vp)\rightarrow
D_\Gamma(\sigma^\prime_{\lambda^\prime},\vp^\prime)$.\hB

We can now define the scattering matrix.
We fix a maximal compact subgroup $K\subset G$ and a parabolic subgroup
$P\subset G$. Let $P=MAN$ be the associated Langlands decomposition and $w\in
N_K(M)/M$ be a representative of the non-trivial element of the Weyl group.
Then we have the meromorphic family of  Knapp-Stein intertwining operators
$\hat J^w_{\sigma_\lambda,\vp}:C^{-\infty}(\partial
X,V(\sigma_\lambda,\vp))\rightarrow C^{-\infty}(\partial
X,V(\sigma^w_{-\lambda},\vp))$, where $\sigma^w(m)=\sigma(m^{w^{-1}})$.
If $\Ree(\lambda)<0$ and $f\in C^{\infty}(\partial
X,V(\sigma_\lambda,\vp))$, then it is given by
$$\hat J^w_{\sigma_\lambda,\vp}(f)(g)=\int_{\bar N} f(gw\bar n) d\bar n\ .$$

Let $(\phi_\mu)_{\mu\in\aca}$, $\phi_\mu\in D_\Gamma(\sigma_\mu,\vp)$, be
a germ of a meromorphic family at $\lambda$. Then
$(J^w_{\sigma_\mu,\vp}(ext^\Gamma\phi_\mu))_\mu$,
$J^w_{\sigma_\mu,\vp}(ext^\Gamma\phi_\mu)\in{}^\Gamma C^{-\infty}(\partial
X,V(\sigma^w_{-\mu},\vp)) $, is a germ of a meromorphic family.
\begin{ddd}\label{sctdeff} We define the germ of a meromorphic family
$(\hat S^w_{\sigma_\mu,\vp}(\phi_\mu))_\mu$ at $\lambda$, $
\hat S^w_{\sigma_\mu,\vp}(\phi_\mu)\in D_\Gamma(\sigma^w_{-\mu },\vp)$, by  
$\hat
S^w_{\sigma_\mu,\vp}(\phi_\mu):=res^\Gamma(\hat
J^w_{\sigma_\mu,\vp}(ext^\Gamma\phi_\mu))$. \end{ddd} Note that for generic
$\lambda$ we have a well-defined continuous map
$\hat S^w_{\sigma_\lambda,\vp}:D_\Gamma(\sigma_\lambda,\vp)\rightarrow
D_\Gamma(\sigma^w_{-\lambda},\vp)$.
The identity 
$$
ext^\Gamma\circ \hat S^w_{\sigma_\mu,\vp}= \hat J^w_{\sigma_\mu,\vp}\circ
ext^\Gamma\ .$$ 
is an immediate consequence of the definitions.
Strictly speaking, this equation holds after applying both sides to a germ of
a meromorphic family.

There is a natural pairing $V_{B_\Gamma}(\sigma_\lambda,\vp)\otimes
V_{B_\Gamma}(\tilde\sigma_{-\lambda},\tilde\vp)\rightarrow
V_{B_\Gamma}(-\rho)$. Since $V_{B_\Gamma}(-\rho)$ is the bundle of densities
on $B_\Gamma$ its sections can be integrated. 
In \cite{bunkeolbrich99} we have shown that if $\max_{P\in\cP^<}(
\delta_{\vp_{|P_\Gamma}} - \rho^{\Gamma_P})   <0$, then we have a
non-degenerate paring $B_\Gamma(\sigma_\lambda,\vp)\otimes
B_\Gamma(\tilde\sigma_{-\lambda},\tilde\vp)\rightarrow \C$ and thus an
inclusion $B_\Gamma(\sigma_\lambda,\vp)\hookrightarrow
D_\Gamma(\sigma_\lambda,\vp)$. Furthermore, we have shown the following result 
\begin{lem}\label{bser}
 If $(f_\mu)_\mu$ is a germ at $\lambda$
of a meromorphic family, $f_\mu\in B_\Gamma(\sigma_\mu,\vp)$, then
$(\hat S^w_{\sigma_\mu,\vp} (f_\mu))_\mu$ is a germ at $\lambda$ of a
meromorphic family, $ \hat S^w_{\sigma_\mu,\vp}
(f_\mu)\in B_\Gamma(\sigma^w_{-\mu},\vp)$. \end{lem}

We now prove Theorem \ref{cfad} assuming that $res^\Gamma$ is regular outside
$I_\aaaa$. Let $\phi\in \Cusp(\sigma_\lambda,\vp)$. By Theorem \ref{gfa} we
know that $\phi$ is stably deformable.
After suitable twisting and embedding we can
assume that $\phi\in \Cusp(\sigma_\lambda,\vp)\cap
\Fam_\Gamma(\sigma_\lambda,\vp)$ and there are no cusps of full rank. Let
$(\phi_\mu)_\mu$ be a germ of a holomorphic family  such that
$\phi_\lambda=\phi$.  Since $\supp(\phi)\subset \Lambda_\Gamma$ we conclude
that $\{res^\Gamma \phi_\mu\}_{|\mu=\lambda}=0$, where $\{\psi\}$ denotes the
restriction of $\psi\in D_\Gamma(\sigma_\lambda,\vp)$ to
$S_\Gamma(\tilde \sigma_{-\lambda},\tilde\vp)\subset B_\Gamma(\tilde
\sigma_{-\lambda},\tilde\vp)$ (here we use that fact that
$C_c^\infty(B_\Gamma,V_{B_\Gamma}(\tilde\sigma_{-\lambda},\tilde\vp)$ is dense
in $S_\Gamma(\tilde \sigma_{-\lambda},\tilde\vp)$).
Using the decomposition \cite{bunkeolbrich99}, (22), 
we can write
$$(res^\Gamma\phi_\mu)_{\mu=\lambda}=\sum_{P\in\tilde \cP} (T_P^*
res^{\Gamma_P} \phi_\mu)_{\mu=\lambda}\ .$$
It suffices to show that
$(T_P^* res^{\Gamma_P} \phi_\mu)_{\mu=\lambda}=0$ for all $P\in \tilde\cP$.
Let $\chi\in C_c^\infty(P_\Gamma)$ be such that
$\sum_{\gamma\in\Gamma_P}\gamma^*\chi\equiv 1$.
Then we can write
$$\phi_P:=\int_{P_\Gamma} \chi(x)\pi^{\sigma_\lambda,\vp}(x)\phi \: dx\ .$$
Since we assume that $\lambda\not\in I_\aaaa$ the push-down $\pi^{\Gamma_P}_*$
is regular and surjective at $\lambda$.
In fact, the poles of $\pi^{\Gamma_P}_*$ are the same as the poles of
$\pi^{P_\Gamma}_*$ (see the definition below). The location of
the latter is determined in \cite{bunkeolbrich99}, Lemma 3.14.  
Furthermore, regularity of $res^\Gamma$ is
equivalent to surjectivity of $\pi^{\Gamma_P}_*$. Thus in order to check that
$T_P^*(res^\Gamma \phi_\mu)_{\mu=\lambda}=0$ it suffices to show
that
$\langle \chi_P res^{\Gamma_P}
\phi_\mu,\pi^{\Gamma_P}_*f_{-\mu}\rangle_{\mu=\lambda}=0$
for all germs $(f_\mu)_\mu$ at $-\lambda$ of  holomorphic families
$f_\mu\in C^\infty(\partial X,V(\tilde\sigma_\mu,\tilde\vp))$.
In fact, since $\{\chi_P  res^{\Gamma_P}\phi_\mu\}_{\mu=\lambda}=0$
and $B_{\Gamma_P}(\tilde\sigma_{\mu},\tilde
\vp)$ is a $C_c^\infty(B_{\Gamma_P})$-module
it suffices to show 
$$\langle   res^{\Gamma_P}
\phi_\mu,\pi^{\Gamma_P}_*f_{-\mu}\rangle_{\mu=\lambda}=0$$
for all families with $\supp(\pi^{\Gamma_P}_*f_\mu)\subset
\{\chi_P=1\}$.
We now consider
$$\pi^{P_\Gamma}_*(f_\mu):=\int_{P_\Gamma} \pi^{\tilde\sigma_{\mu},\tilde
\vp}(x) (f_\mu)_{|\Omega_{\Gamma_P}} dx\ .$$
In \cite{bunkeolbrich99} we have shown that
$\pi^{P_\Gamma}_*(f_\mu)\in B_{\Gamma_P}(\tilde\sigma_{\mu},\tilde
\vp)$, and that $(\pi^{\Gamma_P}(f_\mu)-\pi^{P_\Gamma}_*(f_\mu))_\mu$ is a germ
of a holomorphic family in $S_\Gamma(\tilde\sigma_{\mu},\tilde
\vp)$.  Since $\pi_*^{P_\Gamma}$ factorizes over $\pi_*^{\Gamma_P}$  we also
  have $\supp(\pi^{P_\Gamma}_*f_\mu)\subset \{\chi_P=1\}$. Because of 
$\{\chi_P  res^{\Gamma_P}\phi_\mu\}_{\mu=\lambda}=0$ we get
\begin{equation}\label{gha}\langle   res^{\Gamma_P}
\phi_\mu,\pi^{\Gamma_P}_*f_{-\mu}\rangle_{\mu=\lambda}=\langle   
res^{\Gamma_P} \phi_\mu,\pi^{P_\Gamma}_*f_{-\mu}\rangle_{\mu=\lambda}\
.\end{equation}
 We compute
\begin{eqnarray*}
\langle   res^{\Gamma_P}
\phi_\mu,\pi^{P_\Gamma}_*f_{-\mu}\rangle_{\mu=\lambda}&=&
\langle   res^{\Gamma_P}
\phi_\mu,\pi^{\Gamma_P}_*\int_{P_\Gamma}
\chi(x)\pi^{\sigma_{-\mu},\tilde\vp}(x)f_{-\mu}dx \rangle_{\mu=\lambda}\\
&=&
\langle ext^{\Gamma_P} res^{\Gamma_P}
\phi_\mu, \int_{P_\Gamma}
\chi(x)\pi^{\sigma_{-\mu},\tilde\vp}(x)f_{-\mu}dx \rangle_{\mu=\lambda}\\
&=&\left(\int_{P_\Gamma}
\chi(x) 
\langle  
\phi_\mu,  
\pi^{\sigma_{-\mu},\tilde\vp}(x)f_{-\mu}\rangle dx \right)_{\mu=\lambda}\\
&=&\langle  
(\phi_\mu)_{P}, 
f_{-\mu} \rangle_{\mu=\lambda}\\
&=&0\ .
\end{eqnarray*}
We now apply (\ref{gha}).
\hB

\section{Some functional analytic preparations}
\subsection{Complex interpolation}\label{cint}

For $p\in (1,\infty)$ and $r\in\nat_0$ the space
$H^{p,r}(\partial X,V(\sigma_\lambda,\vp))$ is the usual Sobolev space
of sections of $V(\sigma_\lambda,\vp)$ which have distributional derivatives
in $L^p$ up to order $r$. This space is well-defined, but there is no
prefered norm. There are several slightly different ways to define Sobolev spaces of non-integral order. 
For our purpose a natural choice is complex
interpolation. If $r_0,r_1\in \nat_0$, $\theta\in (0,1)$, and
$r:=r_0+\theta(r_1-r_0)$, then we define (see \cite{reedsimon75}, IX.4,
\cite{berghloefstroem76}, Ch. 6)  $$H^{p,r}(\partial
X,V(\sigma_\lambda,\vp)):=[H^{p,r_0}(\partial
X,V(\sigma_\lambda,\vp)),H^{p,r_1}(\partial X,V(\sigma_\lambda,\vp))]_\theta\
.$$ If $r\in\nat_0$, then this definition coincides with the former. We extend
the scale of Sobolev spaces to negative orders by duality  $$H^{p,-r}(\partial
X,V(\sigma_\lambda,\vp)):=H^{q,r}(\partial
X,V(\tilde\sigma_{-\lambda},\tilde\vp))^*\ ,$$ where $r\ge 0$ and
$\frac{1}{p}+\frac{1}{q}=1$. Then for  $r_0,r,r_1\in\R$ with $r_0<r <r_1$,
$r=r_0+\theta(r_1-r_0)$ and $p_0,p,p_1\in(1,\infty)$,
$\frac{1}{p}=\frac{1}{p_0}+\theta (\frac{1}{p_1}-\frac{1}{p_0})$,  we have
\begin{equation}\label{interpol}H^{p,r}(\partial
X,V(\sigma_\lambda,\vp)):=[H^{p_0,r_0}(\partial
X,V(\sigma_\lambda,\vp)),H^{p_1,r_1}(\partial X,V(\sigma_\lambda,\vp))]_\theta\
.\end{equation} We further will consider the Fr\'echet spaces
$$H^{p,<r}(\partial X,V((\sigma_\lambda,\vp)):=\bigcap_{s<r}  H^{p,s}(\partial
X,V((\sigma_\lambda,\vp))\ .$$

For the convenience of the reader we recall the definition of the
interpolation space (\ref{interpol}). 
The pair of spaces $$\left(H^{p_0,r_0}(\partial X,V(\sigma_\lambda,\vp)),
H^{p_1,r_1}(\partial X,V(\sigma_\lambda,\vp))\right)$$ considered as subspaces
of $C^{-\infty}(\partial X,V(\sigma_\lambda,\vp))$ forms a couple in the sense
of \cite{berghloefstroem76}, Sec. 2.3. There is a natural  Banach norm
on the sum $\Sigma:=H^{p_0,r_0}(\partial
X,V(\sigma_\lambda,\vp))+H^{p_1,r_1}(\partial X,V(\sigma_\lambda,\vp))$.
 Let
$\cF$ be the space of all bounded, continuous functions $f$ on the strip
$\{\Ree(z)\in [0,1]\} $ with values in $\Sigma$ which
 are holomorphic on the interior of this strip, such that for $t\in\R$ we have
that $f(\imath t)\in H^{p_0,r_0}(\partial
X,V(\sigma_\lambda,\vp))$ and $f(\imath t+1)\in H^{p_1,r_1}(\partial
X,V(\sigma_\lambda,\vp))$ are continuous and
$$\sup_{t\in\R } \|f(\imath
t)\|_{H^{p_0,r_0}}<\infty\ ,\quad 
  \sup_{t\in\R } \|f(\imath
t+1)\|_{H^{p_1,r_1}}<\infty\ .$$
Then $\sup_{t\in \R} \{ \|
f(\imath t) \|_{H^{p_0,r_0 }},\|f(\imath t+1)\|_{H^{p_1,r_1 }} \}$ is a norm
making $\cF$ into a Banach space. Let $\cF_\theta$ be the closed subspace of
all $f\in\cF$ with $f(\theta)=0$. Then $$[H^{p_0,r_0}(\partial
X,V(\sigma_\lambda,\vp)),H^{p_1,r_1}(\partial
X,V(\sigma_\lambda,\vp))]_\theta:= \cF/\cF_\theta\ .$$

We employ complex interpolation in the following way.
We will show for a fixed distribution
$f\in C^{-\infty}(\partial X,V(\sigma_\lambda,\vp))$ that  
$f\in H^{p_0,r_0}\cap  H^{p_1,r_1}$. If $\theta\in
(0,1)$, $r=r_0+\theta(r_1-r_0)$ and
$\frac{1}{p}=\frac{1}{p_0}+\theta(\frac{1}{p_1}-\frac{1}{p_0})$, then
we conclude that
$f\in  H^{p,r}$.
In order to see this we consider $f$ as a constant function from
$\{\Ree(z)\in [0,1]\}$ to $\Sigma$ and apply the interpolation result
(\ref{interpol})  above.

\subsection{Regularity of intertwining operators}\label{zc}

In this subsection we consider the Sobolev spaces $H^{p,r}:=H^{p,r}(\R^n)$ for
$p\in (1,\infty)$ and $r\in\R$ and the regularity properties
of the operator $A_s:C^{\infty}_c(\R^n)\rightarrow C^{-\infty}(\R^n)$
which is given by the convolution kernel  $a_s(x):=\|x\|^{-n-2s}$ for $s\in\C$,
where $\|x\|^2:=\sum_{i=1}^n x_i^2$. To be precise, $a_s$ is a regular
distribution for $\Ree(s)<0$, and it is defined by meromorphic continuation for
all $s\in\C$.  Let $\Delta:=\sum_{i=1}^n (\frac{\partial}{\partial x_i})^2$ be
the Laplace operator. Then we have
$\Delta \|x\|^{-n-2s}=(2s+n)(2s+2)    \|x\|^{-n-2-2s}$.
Iterating the equation
\begin{equation}\label{iter} A_{s+1}=\frac{1}{(2s+n)(2s+2)}\Delta
A_s\ .\end{equation} we obtain the required meromorphic continuation.

For $\Ree(s)=0$, $s\not=0$ the operator $A_s$ is a singular integral operator.
It is well-known that $A_s\in B(H^{p,r})$ for $p\in(1,\infty)$, $r\in\R$.  
If $\Ree(s)\in\nat_0$, then we can study the mapping properties of
$A_s$ by reduction to the case $\Ree(s)=0$ using (\ref{iter}).
If $\Ree(s)\ge 0$ is not integral, then we can reduce to the case
$\Ree(s)\in(-1,0)$. In this case $A_s$ is given by convolution by a regular
distribution. But $A_s$ does not act nicely on $H^{p,r}$ since this
distribution grows at infinity. So we must restrict the domain and extend
the range.

Let $H_c^{p,r}\subset H^{p,r}$ be the subspace of all elements of compact
support, and let $H_{loc}^{p,r} \subset C^{-\infty}(\R^n)$ be the space
of all distributions $f$ such that $\chi f\in H^{p,r}$ for all $\chi\in
C_c^\infty(\R^n)$. We define
$H^{p,<r}_{loc}:=\bigcap_{s<r} H^{p,s}_{loc}$.

If $m\in\nat_0$, then let $A_m=(s-m)^{-1} B_{m,-1}+B_{m,0}+\dots$
be the Laurent expansion of $A_s$ at $s=m$.
Our main result is the following theorem.
\begin{theorem}\label{reg}
If $\Ree(s_0)\ge 0$ and $s_0\not\in\nat_0$, then we have 
$$A_s:H_c^{p,r}\rightarrow H_{loc}^{p,r-2\Ree(s_0)}\ .$$
For $m\in\nat_0$ we have
\begin{eqnarray*}
B_{m,-1}&:&H^{p,r}\rightarrow H^{p,r-2m}\ ,\\
B_{m,k}&:&H_c^{p,r}\rightarrow H_{loc}^{p,<r-2m} ,\quad k\in\nat_0\ .
\end{eqnarray*}
\end{theorem}
\proof
The space $H^{p,r}_c$ is a limit of Banach spaces, while
$H^{p,r}_{loc}$ is a Fr\'echet space in a natural way.
For each compact $K\subset \R^n$ let 
$H_K^{p,r}:=\{f\in H^{p,r}\:|\:\supp(f)\subset K\}$. Then $H^{p,r}_K$ is a
closed subspace of $H^{p,r}$, and $H_c^{p,r}=\lim_{\stackrel{\rightarrow}{K}}
H^{p,r}$. For each $\chi\in C^\infty_c(\R^n)$ we define the seminorm 
$q_{\chi,p,r}(f):=\|\chi f\|_{H^{p,r}}$ on $H^{p,r}_{loc}$. The family of
these seminorms defines the topology of $H^{p,r}_{loc}$.
If $A:H_c^{p,r}\rightarrow H_{loc}^{p,s}$ is continuous, then for each
$K$ and $\chi$, we define
$$q(A)_{K,\chi,p,r,s}:=\sup_{f\in H^{p,r}_K,\|f\|_{H^{p,r}}=1 }
q_{\chi,p,s}(Af)\ .$$
First we assume that
$\Ree(s)\in[-1,0]$. In order to fix the pole at $s=0$ we set $\hat A:= s A_s$.
\begin{prop}\label{rr}
For each compact subset $K\subset \R^n$ and $\chi\in
C_c^\infty(\R^n)$ there exists $C<\infty$ such that for all $s\in \C$ with
$\Ree(s)\in [-1,0]$ we have
$$q_{K,\chi,p,0,0}(\hat A_s)\le C
(1+|s|)^3\ .$$
\end{prop}
\proof
Fix a compact $K\subset \R^n$ and $\chi\in C_c^\infty(\R^n)$.
Let 
\begin{equation}\label{rdef}
R:=\sup_{x\in K,y\in\supp(\chi)}\dist(x,y)\ .\end{equation}
Let $\kappa\in C^\infty_c(\R^n)$ be a cut-off function such that
$\kappa(x)=1$ for $\|x\|\le R$ and $\kappa(x)=0$ for $\|x\|\ge 2R$. We define
$a_s^\kappa:=\kappa a_s$ and let $A^\kappa_s$ be the convolution operator
defined by $a_s^\kappa$. If $\supp(f)\subset K$, then we have $\chi A_s f=\chi
A_s^\kappa(f)$. It suffices to show
$$\| s A_s^\kappa\|_{B(L^p)}\le C(1+|s|)^3\ ,$$
where $C$ is independent of $s \in \{\Ree(s)\in[-1,0]\}$.
If $\Ree(s)<0$, then $a_s^\kappa\in L^1(\R^n)$, so that by Yang's
inequality $$\|A_s^\kappa\|_{B(L^p)}\le \|a_s^\kappa\|_{L^1}\ .$$
There is a constant $C<\infty$ such that 
$ \|a_s^\kappa\|_{L^1}\le C$ uniformly for $\{\Ree(s)\in [-1,-1/8]\}$,
and hence 
\begin{equation}\label{yang} \|A_s^\kappa\|_{B(L^p)}\le C\ .\end{equation}
We cannot use this estimate close to the imaginary axis, since in that
region $\|a^\kappa_s\|_{L^1}$ explodes like
$\Ree(s)^{-1}$, and we need a uniform estimate. The point is that we must take
into account the oscillating behaviour of $a_s^\kappa(x)$ near $x=0$. Then
the divergence above only occurs at $s=0$, and this is fixed because we
consider $\hat A_s^\kappa=s A^\kappa_s$.

\begin{lem}\label{l2}
There exists $C<\infty$ such that on $\{\Ree(s)\in[-1/4,0]\}$
$$\|\hat A_s^\kappa\|_{B(L^2)}\le C(1+|s|)\ .$$
\end{lem}
\proof
We first estimate $\|A_s^\kappa\|_{B(L^2)}$ using the Fourier transform.
The Fourier transform $\cF(a_s^\kappa)$ 
is the convolution of $\cF(a_s)$ and $\cF(\kappa)\in \cS(\R^n)$.
We  have 
$$\cF(a_s)(\xi)= c \: 2^{-2s} \frac{\Gamma(-s)}{\Gamma(\frac12 n+s)}\|\xi\|^{2s}\ ,$$ where $c$ is independent of $s$.
Note that (see \cite{ryshikgradstein57}, 6.328)
$$\lim_{|y|\to\infty} |\Gamma(x+\imath y)| e^{\frac{\pi}{2}|y|}
|y|^{\frac12-x}=\sqrt{2\pi}\ ,$$
locally uniformly for $x\in \R$.
Using this property of the $\Gamma$-function
we see that there is a constant $C<\infty$ such that for $\Ree(s)\in
[-1/4,0]$  we have  $|s\frac{\Gamma(-s)}{\Gamma(\frac12
n+s)}| \le C(1+ \Imm(s))^{-2\Ree(s)-n/2+1}$.
Let $\chi_1\in C_c^\infty(\R^n)$ be a cut-off function such that
$\chi_1(\xi)=1$ for $\|\xi\|\le 1$. Then uniformly on 
$\{\Ree(s)\in [-1/4,0]\}$ we have $\|(1-\chi_1) \cF(\hat a_s)\|_{L^\infty}\le
C(1+|s|)$. Furthermore, on  $\{\Ree(s)\in[-1/4,0]\}$ we obtain a uniform
estimate $\|\chi_1 \cF(\hat a_s)\|_{L^1}\le C (1+|s|)$.
Writing
$\cF(\hat a_s^\kappa)=\cF(\kappa)*\chi_1 \cF(\hat a_s)+\cF(\kappa)*(1-\chi_1)
\cF(\hat a_s)$ we can estimate
$$\|\cF(\hat a_s^\kappa)\|_{L^\infty}\le
\|\cF(\kappa)\|_{L^\infty} \|\chi_1 \cF(\hat
a_s)\|_{L^1}+\|\cF(\kappa)\|_{L^1}\|(1-\chi_1) \cF(\hat a_s)\|_{L^\infty}\le
C(1+|s|)$$
uniformly on $\{\Ree(s)\in[1/4,0]\}$. This gives
$$\|\hat A_s^\kappa\|_{B(L^2)}\le C(1+|s|)$$ uniformly on
$\{\Ree(s)\in[1/4,0]\}$.
\hB

In order to extend the assertion of Lemma \ref{l2} from $L^2$ to $L^p$ we need
the following estimate. 
\begin{lem}\label{smooth}
There exists a constant $C<\infty$ such that if $\{\Ree(s)\in[-1/4,0]\}$, then
\begin{equation}\label{bg5}\int_{\dist(y,z)\ge 4 \dist(x,z)} |\hat
a_s^\kappa(y-x)-\hat a_s^\kappa(y-z)| dy\le C(1+|s|^3)\end{equation}
for all $x,z\in\R^n$.
\end{lem}
\proof
On the domain of integration we have
$$\frac12 \dist(z,y)\le \dist(x,y)\le 2\dist(z,y)\ .$$
Recall the definition (\ref{rdef}) of $R$.
We decompose the domain of integration in (\ref{bg5}) into the regions
$\{\dist(z,y)\le R/2\}$, $\{R/2\le \dist(z,y)\le 2R\}$, $\{\dist(z,y) \ge
2R\}$, and let $I_1,I_2,I_3$ denote the corresponding integrals.
Since the integrand vanishes on the third region we have $I_3=0$.
Next we consider $I_1$. In this region we can replace $\hat a_s^\kappa$ by
$\hat a_s$.
For $0\not=x\in\R^n$ let $x^0:=\frac{x}{\|x\|}$.
We compute
\begin{eqnarray}
 \lefteqn{|\hat a_s(y-z)-\hat a_s(y-x)|}&&\\&=& |s | \left|\|y-z\|^{-n-2s} -
\|y-x\|^{-n-2s}\right|\nonumber \\ 
&=&|s| \|y-z\|^{-n-2\Ree(s)} \left|1 -
\frac{\|y-z+z-x\|^{-n-2s}}{\|y-z\|^{-n-2s}}\right|\nonumber\\  &=&|s|
\|y-z\|^{-n-2\Ree(s)} \left|1 -
\|(y-z)^0+\frac{z-x}{\|y-z\|}\|^{-n-2s}\right|\nonumber \\ &=&|s|
\|y-z\|^{-n-2\Ree(s)} \left|1 - \left(1+2\left\langle (y-z)^0,
\frac{z-x}{\|y-z\|}\right\rangle +
\frac{\|z-x\|^2}{\|y-z\|^2}\right)^{\frac{-n-2s}{2}}\right|\label{thiseq}\ .
\end{eqnarray} 
Note that $b:= 2\langle (y-z)^0, \frac{z-x}{\|y-z\|}\rangle +
\frac{\|z-x\|^2}{\|y-z\|^2}$ satisfies $|b| \le   1/2 + 1/16< 3/4 $. For
$b\in\R$, $|b|<3/4$ we write $(1+b)^{-n-2s}=\exp((-n-2s) \log(1+b))$.
Furthermore, $\log(1+b)=b(1+r_1)$, where $\sup_{|b|\le 3/4} | r_1|
<\infty$. Then we expand
$\exp((-n-2s) b(1+r_1))=1+(-n-2s) b(1+r_1)+ b r_2$,
where $|r_2|\le 
\frac12 |b|[(-n-2s) (1+r_1)]^2
\sup_{\xi\in [0,1]} |\exp(\xi (-n-2s) b(1+r_1))|$.
We conclude that 
$|r_2| < C (1+|s|)^2$, where $C$ is independent of $s$ and $b$.
We further get
$(1-(1+b)^{-n-2s})=(-n-2s)b + b r_3$,
where $|r_3|\le C (1+|s|)^{2}$ uniformly in $s$.
Inserting this estimate into (\ref{thiseq}) we obtain
\begin{eqnarray*}
\lefteqn{|\hat a_s(y,z)-\hat a_s(y,x)|}&&\\&=& |s|  \|y-z\|^{-n-2\Ree(s)}  
\left|(2\left\langle (y-z)^0, \frac{z-x}{\|y-z\|}\right\rangle +
\frac{\|z-x\|^2}{\|y-z\|^2}) ((-n-2s)  +  r_3)\right|\\ &=&|s|
\|y-z\|^{-n-2\Ree(s)} \frac{\|z-x\|}{\|y-z\| } \left|(2\left\langle (y-z)^0,
(z-x)^0\right\rangle + \frac{\|z-x\|}{\|y-z\|}) ((-n-2s)  +  r_3)\right|\\
&\le& (1+|s|)^3 C \|y-z\|^{-n} \frac{\|z-x\|}{\|y-z\| }\ .
\end{eqnarray*}
We now estimate
\begin{eqnarray*}
I_1&\le&(1+|s|)^3 C \|z-x\|
\int_{R/2\ge  \dist(y,z)\ge 4 \dist(x,z)}  \|y-z\|^{-n-1} dy\\
&\le& (1+|s|)^3 C^\prime \ ,
\end{eqnarray*}
where $C^\prime$ can be chosen uniformly for $\Ree(s)\in [-1/4,0]$.
Finally, on the compact region $\{R/2\le \dist(z,y)\le 2R\}$ and for
$\Ree(s)\in [-1/4,0]$
the integrand is uniformly bounded with respect to $x,y,z$. 
Thus 
$$I_2\le C$$  uniformly on  $\{\Ree(s)\in [-1/4,0]\}$.
\hB

\begin{lem}
For $p\in (1,\infty)$
there exists $C<\infty$ such that on $\{\Ree(s)\in[-1/4,0)\}$
$$\|\hat A_s^\kappa\|_{B(L^p)}\le C(1+|s|)^3\ .$$
\end{lem}
\proof
We follow the argument of \cite{coiffmanweiss71}, Thm. III.2.4.
There the following fact was shown.
Let $K$ be an integral operator  given by an integral kernel $k(x,y)$.
Assume that $K$ extends to a bounded operator on $L^2$ and
$\|K\|_{B(L^2)}\le C_1$. Furthermore, assume that there are constants
$C_2,C_3$ such that for all $z,x\in\R^n$
$$\int_{\dist(y,z)\ge C_2 \dist(x,z)} |k(y,x)-
k(y,z)| dy \le  C_3\ .$$
Then $K$ extends to a bounded operator on $L^p$, $p\in (1,2]$, and
$\|K\|_{B(L^p)}\le C(p,n,C_2) (1+C_1+C_3)$.

For $\Ree(s)<0$ the operator $\hat A^\kappa_s$ is given by an integral kernel.
We set $C_2=4$. Then 
we obtain the assertion of the lemma by combining
Lemma \ref{smooth} (which estimates $C_3$) and Lemma \ref{l2}
(which provides an estimate for $C_1$). 
For $p\in[2,\infty)$ we argue by duality using 
$(\hat A_s^\kappa)^*=\hat A_{\bar{s}}^\kappa$.
\hB

We now finish the proof of Proposition \ref{rr}. We must
discuss the case $\Ree(s)=0$. Let $f\in\cS(\R^n)$.
Then we know that in $\cS(\R^n)$ and thus in $L^p(\R^n)$
$$\lim_{\epsilon\downarrow 0} \hat A_{s-\epsilon}^\kappa f = \hat A_s^\kappa f\
.$$ Since we have the bound
$\|\hat A_{s-\epsilon}^\kappa\|_{B(L^p)}\le C(1+|s|)^3$ uniformly for
$\epsilon\in (0,1)$ we conclude that  
$\|\hat A_s^\kappa\|_{B(L^p)}\le C(1+|s|)^3$ uniformly if $\Ree(s)\in [-1,0]$.
\hB

\begin{lem}\label{caser}
Let $r\in \R$.
For each compact subset $K\subset \R^n$ and $\chi\in
C_c^\infty(\R^n)$ there exists $C<\infty$ such that for all $s\in \C$ with
$\Ree(s)\in [-1,0]$ we have
$$q_{K,\chi,p,r,r}(\hat A_s)\le C
(1+|s|)^3\ .$$
\end{lem}
\proof
It again suffices to show that
$$\|\hat A_s^\kappa\|_{B(H^{p,r})}\le C(1+|s|)^3\ .$$
Let $J_r$ be the operator on $\cS^\prime(\R^n)$ given by
$\cF(J_r)f(\xi)=(1+|\xi|^2)^{r/2} \cF(f)(\xi)$.
Then it is well-known that
$H^{p,r}=\{f\in \cS^\prime(\R^n)\:|\: J_rf\in L^p\}$, and
$\|f\|_{H^{p,r}}:=\|J_r f\|_{L^p}$
is a norm on the Banach space $H^{p,r}$. Since $\hat A_s^\kappa$ commutes with
translation, it also commutes with $J_r$ (when restricted to the Schwartz
space). So we can write $\hat A_s^\kappa=J_{-r} \hat A_s^\kappa J_r$. The right
hand side obviously acts on $H^{p,r}$, and we have the estimate
$$\|\hat A_s^\kappa\|_{B(H^{p,r})}\le \|\hat A_s^\kappa\|_{B(L^p)}\le
C(1+|s|)^3\ .$$ This proves the lemma. \hB

We now consider the case that $\Ree(s)\in [m-1,m]$, $m\in \nat$.
In order to fix the poles of $A_s$ at $s=m-1,m$ we consider
$\tilde A_ {m,s}:= s(s+1) A_{s+m}$, $\Ree(s)\in [-1,0]$.
\begin{lem}\label{wel} 
Let $r\in \R$.
For each compact subset $K\subset \R^n$ and $\chi\in
C_c^\infty(\R^n)$ there exists $C<\infty$ such that for all $s\in \C$ with
$\Ree(s)\in [-1,0]$ we have
$$q_{K,\chi,p,r,r-2m}(\tilde A_s)\le C
(1+|s|)^2\ .$$
\end{lem}
\proof
We use (\ref{iter}) in order to reduce to the case $\Ree(s)\in [-1,0]$.
We have
$$\tilde A_{m,s}=  q_m(s)^{-1} \Delta^m \hat A_s\ ,$$
where $q_m$ is some polynomial  of degree $2m-1$ which does not vanish if
$\Ree(s)\in[-1,0]$.
Let $\chi^\prime\in C_c^\infty$ be such that $\chi^\prime\chi=\chi$.
Since $\|\chi \Delta^m\|_{B(H^{p,r},H^{p,r-2m})}\le C$ 
we get for $f\in H_K^{p,r}$ with $\|f\|_{H^{p,r}=1}$
\begin{eqnarray*}\|\chi \tilde A_{m,s}(f)\|_{H^{p,r-2m}}&=&
|q_m(s)^{-1}|\|\chi  \Delta^m \hat A_s(f)\|_{H^{p,r-2m}}\\
&=& |q_m(s)^{-1}|\|\chi  \Delta^m \chi^\prime \hat A_s(f)\|_{H^{p,r-2m}}\\
&\le&|q_m(s)^{-1}|C q_{K,\chi^\prime,p,r,r}(\hat A_s)\\
&\le& C(1+|s|)^2
\end{eqnarray*} 
using Lemma \ref{caser} in the last step. \hB

We now prove the theorem.
First we consider $s_0\in \C\setminus \nat_0$ such that $\Ree(s_0)\ge 0$.
Then there exists $m\in\nat$ such that $\Ree(s_0)\in [m-1,m]$.
Let $K\subset \R^n$ be compact and $\chi\in C^\infty_c(\R^n)$.
Let $f\in H^{p,r}_K$ be given. We must show that
$\chi A_s f\in H^{p,r-2\Ree(s_0)}$.

Assume first that $f\in C_c^\infty(\R^n)\cap H^{p,r}_K$.
We consider the  holomorphic family $F_s$ with values in $\cS(\R^n)$
given by $F_s:=(s+2)^{-3}f$ for $\Ree(s)\in [-1,0]$.
Then $\chi \tilde A_{m,s} F_s$ is a holomorphic family in $C^\infty_c(\R^n)$.
Thus it is holomorphic  with values in $H^{p,r-2m}$ and
$\sup_{\Ree(s)=[-1,0]}   \|\chi \tilde A_{m,s} F_s\|_{H^{p,r-2m}}\le C
\|f\|_{H^{p,r}}$ (see  Lemma \ref{wel}). Furthermore, we know that if
$\Ree(s)=0$, then $\chi \tilde A_{m-1,s} F_{s-1}\in H^{p,r-m+2}$ is uniformly
bounded by $\|f\|_{H^{p,r}}$.  Thus
$\sup_{\Ree(s)=-1}   \|\chi \tilde
A_{m,s} F_s\|_{H^{p,r-2m+2}}\le C \|f\|_{H^{p,r}}$. 
 By complex interpolation we now conclude that
\begin{eqnarray*}\chi \tilde A_{m,s_0-m} F_{s_0}&\in&
H^{p,r-2\Ree(s_0)}\\
\|\chi \tilde A_{m,s_0-m} F_{s_0}\|_{H^{p,r-2\Ree(s_0)}}&\le& C\|f\|_{H^{p,r}}
\end{eqnarray*}

Note that any element of $f\in H^{p,r}_K$ can be approximated by smooth
compactly supported functions. Since $\tilde A_{m,s_0-m}$ and $F_{s_0}$ are
multiples of  $A_{s_0}$ and $f$ we conclude that
$$\chi A_{s_0} f\in H^{p,r-2\Ree(s_0)}\ .$$ 

It remains to consider the integral points $m\in\nat_0$. First of all note
that $B_{m,-1}$ is a multiple of $\Delta^m$.
For $k\in \nat_0$ and $f\in \cS^\prime(\R^n)$ we write
$$B_{m,k}f=\frac{1}{2\pi\imath}\int_{S(m,\epsilon)} \frac{1}{(s-m)^{k+1}} A_s f
ds\ ,$$ where $S(m,\epsilon)$ is the circle with radius $\epsilon>0$ centered
at $m$, and the integral is taken in $\cS^\prime(\R^n)$. If $f\in H^{p,r}$ and
$\chi\in C_c^\infty(\R^n)$, then we have $\| \chi A_s
f\|_{H^{p,r-2m-2\epsilon}}\le C$ uniformly on $S(m,\epsilon)$. For any
$\phi\in C_c^\infty(\R^n)$ we can estimate
\begin{eqnarray*} |\langle \chi B_{m,k}
f,\phi\rangle|&=&|\frac{1}{2\pi\imath}\int_{S(m,\epsilon)} 
\frac{1}{(s-m)^{k+1}} \langle \chi A_s f,\phi\rangle  ds|\\ &\le& 
\frac{1}{2\pi}\int_{S(m,\epsilon)} | \frac{1}{(s-m)^{k+1}}\langle \chi  A_s
f,\phi\rangle|  ds\\ &\le&\frac{1}{\epsilon^{k}} \sup_{s\in S(m,\epsilon)}
\| \chi  A_s f\|_{H^{p,r-2m-2\epsilon}}
\|\phi\|_{H^{q,-r+2m+2\epsilon}}  ds\\ &\le &
C\|\phi\|_{H^{q,-r+2m+2\epsilon}}\ . \end{eqnarray*} We conclude that $\chi
B_{m,k}f\in (H^{q,-r+2m+2\epsilon})^*=H^{p,r-2m-2\epsilon}$. Since
$\epsilon>0$ can be chosen arbitrary small we have $\chi B_{m,k}f\in
H^{p,<(r-2m)}$ as asserted. \hB

We now apply this result to the spherical intertwining operator $\hat
J^w_{1_\mu,\vp}$. Let $\lambda\in\aca$ be such that $\Ree(\lambda)\ge 0$
and $$\hat J^w_{1_\mu,\vp}=\sum_{l}(\frac{\mu-\lambda}{\alpha})^l B_l$$ be the 
Laurent expansion of the meromorphic family $(\hat
J^w_{1_\mu,\vp})_\mu$ at $\lambda$. The coefficients have the following
mapping property. 
\begin{kor}\label{mapint}
For all  $p\in (1,\infty)$ and $r\in \R$
$$B_l: H^{p,<r}\rightarrow H^{p,<r-2\Ree(\lambda)}\ .$$
\end{kor}
\proof
We fix a parabolic subgroup $P$, a maximal compact subgroup $K$, a Langlands
decomposition $P=MAN$ with $M\subset K$, and a representative $w\in N_K(M)$
of the non-trivial element of the Weyl group $N_K(M)/M$. 
Let $G\setminus wP=\bar NMAN$, $g=\bar n(g)m(g)\alpha(g)\tilde n(g)$,
be the Bruhat decomposition. Furthermore, we define $F:N\setminus
\{1\}\rightarrow \bar N\setminus\{1\}$ such that $F(x)P=xwP$.

If
$\Ree(\lambda)<0$ and  $f\in C^\infty(\partial X,V(1_\lambda,\vp))$, then by
definition $$\hat J^w_{1_\lambda,\vp}f(g)=\int_{\bar N} f(gw\bar x) d\bar x\
.$$ We now assume that $\supp(f)\subset \bar NP$. We compute
\begin{eqnarray*}
(\hat J^w_{1_\lambda,\vp}f)(1)&=&\int_{\bar N} f( w \bar x ) d\bar x\\
&=&\int_{\bar N} f( \bar x^w w ) d\bar x\\
&=&\int_{N} f( x w ) dx\\
&=&\int_{ N\setminus\{1\}} \alpha( xw)^{\lambda-\rho}  f(F(x)) dx
 . \end{eqnarray*}
We now identify the group $N$ with its Lie algebra $\naaa$ using the
exponential map, and further we identify $\naaa$ with $\R^{n-1}$ such that
the action of $M$ is orthogonal (this fixes the identification up to
scale) and such that the measure $dx$ identifies with the standard Lebesgue
measure (this fixes the scale).
The function $\alpha(xw)^\alpha$ is $M$-invariant and satisfies
$\alpha(x^aw)^\alpha=a^{2\alpha}\alpha(xw)$, $a\in A$. Thus in the coordinates
it is given by $\alpha(xw)^\alpha=c_1\|x\|^2$ for some constant $c_1>0$.
We identify $\bar N$ with $\bar\naaa$ using the exponential map, and
furthermore, $\bar\naaa$ with $\R^{n-1}$ using the isomorphism
$\bar\naaa\ni \bar X\mapsto X^w\in\naaa$.
The map $F:N\setminus \{0\}\rightarrow \bar N\setminus \{0\}$
is $M$-equivariant and satisfies $F(x^a)=F(x)^a$. Therefore, in our coordinates
$F(x)=c_2\frac{x}{\|x\|^2}$ for some $0\not=c_2\in\R$.
So we can further compute
\begin{eqnarray*}
\int_{ N\setminus\{1\}} \alpha( xw)^{\lambda-\rho}  f(F(x))
dx&=&c_1^{2\lambda/\alpha-(n-1)}  \int_{\R^{n-1}}
\|x\|^{2\lambda/\alpha-(n-1)}  f(c_2\frac{x}{\|x\|^2}) dx\\
&=&(c_1c_2)^{2\lambda/\alpha-(n-1)}  \int_{\R^{n-1}}
\|\bar y\|^{-2\lambda/\alpha-(n-1)}   f(\bar y) d\bar y \end{eqnarray*}
Let now $\bar z\in \bar N$.
Then we define $f_{\bar z}(g)=f(\bar zg)$.
Since $\hat J^w_{1_\lambda,\vp}$ is $G$-equivariant, we have
$(\hat J^w_{1_\lambda,\vp} f)(\bar z)=(\hat J^w_{1_\lambda,\vp} f_{\bar
z})(1)$. In our coordinates we thus have
\begin{eqnarray*}
(\hat J^w_{1_\lambda,\vp} f)(\bar
z)&=&(c_1c_2)^{2\lambda/\alpha-(n-1)}  \int_{\R^{n-1}}
\|\bar y\|^{-2\lambda/\alpha-(n-1)}   f(\bar z+\bar y) d\bar y\\
&=&(c_1c_2)^{2\lambda/\alpha-(n-1)}  \int_{\R^{n-1}}
\|\bar z-\bar y\|^{-2\lambda/\alpha-(n-1)}   f(\bar y) d\bar y\ .
\end{eqnarray*}
Thus up to the harmless factor $(c_1c_2)^{2\lambda/\alpha-(n-1)}$ the
intertwining operator coincides with $A_s$ introduced at the beginning of the
present subsection, where $s=\lambda/\alpha$.
The assertion of the corollary now follows from Theorem \ref{reg}
since we can cover $\partial X$ by finitely many charts of the form
$\bar N P$ (with varying $P$). \hB

\section{Proof of Theorem \ref{main}}
\subsection{Compatibility with twisting and embedding}

We fix $\sigma\in\hat M$, $\lambda\in \aca$, and a twist $\vp$.
If $\sigma^\prime$, $\mu$, $\pi$ and  $T$ is some twisting data, then we have
the map 
$$i_*:C^{-\infty}(\partial
X,V(\sigma_\lambda,\vp))\rightarrow C^{-\infty}(\partial
X,V(\sigma^\prime_\mu,\vp\otimes \pi))\ .$$ 
\begin{lem}\label{tw}
Let $r\in \R$ and $p\in (1,\infty)$.   If $\phi\in
C^{-\infty}(\partial X,V(\sigma_\lambda,\vp))$ satisfies\\ $i_*\phi\in
H^{p,r}(\partial X,V(\sigma^\prime_\mu,\vp\otimes\pi))$, 
then $\phi\in H^{p,r}(\partial X,V(\sigma_\lambda,\vp))$.
\end{lem} \proof
The map $i_*$ is induced by an inclusion $V(\sigma_\lambda,\vp)\hookrightarrow
V(\sigma^\prime_\mu,\vp\otimes\pi)$ of vector bundles. Since
we can find a complementary bundle of $V(\sigma_\lambda,\vp)$ inside
$V(\sigma^\prime_\mu,\vp\otimes\pi)$ there is also a projection from
 $V(\sigma^\prime_\mu,\vp\otimes\pi)$ to $V(\sigma_\lambda,\vp)$.
In local trivializations this projection is given by a matrix with smooth
entries. Since smooth functions act continuously on the Sobolev
space $H^{p,r}$ the projection provides a left-inverse of $i_*$. 
Now the lemma follows.
\hB

Now we discuss embedding. Let $G=G_m$, $\sigma_m\in\hat M_m$ and $\lambda\in
\aca$. If $n>m$, $\sigma_n$ and $T$ is some embedding data, then we have the
map
$$i_*:C^{-\infty}(\partial X^m,V^m((\sigma_m)_\lambda,\vp))\rightarrow
C^{-\infty}(\partial X^n,V^n((\sigma_n)_{\lambda+\rho^m-\rho^n},\vp))\ .$$
\begin{lem}\label{embb}
Let $r\in (-\infty,0]$,  $p\in (1,\infty)$, and $q$ be given by
$\frac{1}{p}+\frac{1}{q}=1$. \\ If $\phi\in C^{-\infty}(\partial
X^m,V^m((\sigma_m)_\lambda,\vp))$  satisfies $i_*\phi\in
H^{p,<r-\frac{n-m}{q}}(\partial
X^n,V^n((\sigma_n)_{\lambda+\rho^m-\rho^n},\vp))$, then $\phi\in
H^{p,<r}(\partial X^m,V^m((\sigma_m)_\lambda,\vp))$.  \end{lem} \proof
The map $T$ induces an inclusion of bundles
$$V^m((\sigma_m)_\lambda,\vp)\hookrightarrow 
V^n((\sigma_n)_{\lambda+\rho^m-\rho^n},\vp)_{|\partial X^m}\ .$$ Let
$E\rightarrow U$ be an  extension of $V^m((\sigma_m)_\lambda,\vp)$ to a 
tubular neighbourhood $U$ of  $\partial X^m$.  It is now clear that
$i_*=i^2_*\circ i^1_*$ is a composition of the push-forward $i^1_*$ of
distributions induced by the inclusion $\partial X^m\hookrightarrow \partial
X^n$ of a regular submanifold, and the inclusion $i^2_*$ induced by the
inclusion of bundles. If  $i_*\phi\in
H^{p,<r-\frac{n-m}{q}}(\partial
X^n,V^n((\sigma_n)_{\lambda+\rho^m-\rho^n},\vp))$, then
$i^1_*(\phi)\in H^{p,<r -\frac{n-m}{q}}(U,E)$ by Lemma \ref{tw}. It
remains to conclude that $\phi\in H^{p,<r}(\partial
X^m,V^m((\sigma_m)_\lambda,\vp))$.

Fix $\epsilon>0$.
We first assume that $n=m+1$.
Let $l:=[-r+\epsilon+\frac{1}{q}]$ be the integral part of
$-r+\epsilon+\frac{1}{q}$.  We define a restriction
operator  $\cR: C^\infty(U,E^+)\rightarrow \bigoplus_{i=0}^l C^\infty(\partial
X^m,E^+_{|\partial
X^m})$ by
$$\cR(f):=(f_{|\partial
X^m}, (\partial_n f)_{|\partial X^m},\dots,(\partial^l_n f)_{|\partial X^m})\
,$$ where $\partial_n$ denotes the normal derivative,
and $E^+:=E^*\otimes \Lambda$, $\Lambda$ being the bundle of densities.
We have an exact sequence
$$0\rightarrow H_0^{q,-r+\epsilon+\frac{1}{q}}(U,E^+)\rightarrow
H^{q,-r+\epsilon+\frac{1}{q}}(U,E^+)\stackrel{\cR}{\rightarrow}
\bigoplus_{i=0}^l B^{-r+\epsilon-i}_{q,q}(\partial
X^m,E^+_{|\partial
X^m})\rightarrow 0\ ,$$
(see \cite{triebel78} for the definition of the Besov spaces and
\cite{triebel78}, 4.7.1,  for the characterization of the range and the kernel
of $\cR$). In particular, the space of smooth sections of $E^+$ vanishing to
infinite order along $\partial X^m$ is dense in $\ker \cR$.

Since $i^1_*(\phi)\in H^{p,r-\epsilon-\frac{1}{q}}(U,E)$ we know that
$i^1_*\phi$ extends continuously to $H^{q,-r+\epsilon+\frac{1}{q}}(U,E^+)$.
Furthermore, since $i^1_*\phi$ vanishes on the space of smooth sections of
$E^+$ vanishing to infinite order along $\partial X^m$, we conclude that
$i^1_*\phi$ vanishes on $H_0^{q,-r+\epsilon+\frac{1}{q}}(U,E^+)$ and thus
factors over the quotient $\bigoplus_{i=0}^l
B^{-r+\epsilon-i}_{q,q}(\partial
X^m,E^+_{|\partial
X^m})$.
This means that $\phi$ extends continuously to
$B^{-r+\epsilon}_{q,q}(\partial X^m,E^+_{|\partial
X^m})$.
We now use the continuous embedding (\cite{triebel78}, 4.6.1)
$$H^{q,-r+2\epsilon}\hookrightarrow
B^{-r+2\epsilon}_{q,\max(2,q)}\hookrightarrow
B^{-r+\epsilon}_{q,q}$$
in order to conclude that $\phi$ restricts to a continuous functional on
$H^{q,-r+2\epsilon}$.
Since this holds for all $\epsilon>0$ we have $\phi\in H^{p,<r}(\partial
X^m,V^m((\sigma_m)_\lambda,\vp))$. 
If $n-m>1$, then we argue by induction.
\hB

We now state the following stable version of our main theorem.
\begin{theorem}\label{stmain}
\begin{enumerate}
\item
If $\phi\in \Fam_\Gamma^{st}(\Lambda_\Gamma,\sigma_\lambda,\vp)$, then
for $p$ in a generic subset of $(1,\infty)$
 there
exist suitable twisting and embedding data such that
$i_*\phi\in H^{p,<r_{p,\mu}(\Gamma)^\prime}(\partial
X^\prime,V(\sigma_\mu^\prime,\vp^\prime))$, where we attach a "$\:{}^\prime\:$"
to objects which are associated to $\partial X^\prime$ and $\vp^\prime$.
\item If $\phi\in \Cusp_\Gamma(\sigma_\lambda,\vp)$ and $\lambda\not\in
I_\aaaa$, then for all $p\in (1,\infty)$ there exist suitable twisting and
embedding data such that $i_*\phi\in H^{p,<r^0_{p,\mu}(\Gamma)^\prime}(\partial
X^\prime,V(\sigma_\mu^\prime,\vp^\prime))$.
\end{enumerate}
\end{theorem}
 
The contents of this theorem is that given a stably deformable invariant
distribution which is strongly supported on the limit set, or given a cusp
form and $\lambda\not\in I_\aaaa$, then after suitable embedding and twisting
it belongs to the Sobolev space as asserted by  Theorem \ref{main} (for
generic $p$).

\begin{prop}\label{impl}
Theorem \ref{stmain} implies Theorem  \ref{main}.
\end{prop}
\proof
Assume that $i_*$ is induced by some embedding data. 
In this case $\vp=\vp^\prime$, $\mu=\lambda+\rho-\rho^\prime$,
$\delta_\Gamma^\prime=\delta_\Gamma+\rho-\rho^\prime$.
We observe that 
\begin{eqnarray*}
r_{p,\lambda}(\Gamma)^\prime&=&r_{p,\mu}(\Gamma)  -  
\frac{\dim(\partial X^\prime)-\dim(\partial X)}{q}\\
(\:\:\: r^0_{p,\lambda}(\Gamma)^\prime&=&r^0_{p,\mu}(\Gamma)  -  
\frac{\dim(\partial X^\prime)-\dim(\partial X)}{q}\:\:\:)\ .
\end{eqnarray*}
If $r_{p,\lambda}(\Gamma)>0$ ($\:\:$
$r^0_{p,\lambda}(\Gamma)>0$ $\:\:$), then we employ the vanishing result
\cite{bunkeolbrich011} in order to conclude that $\phi=0$. Else we combine the
regularity of $i_*\phi$ asserted in Theorem \ref{stmain} with Lemma \ref{embb}
in order to see that $\phi$ has the regularity claimed in Theorem \ref{main}
provided $p$ is generic or $\phi$ is a cusp form and
$\lambda\not\in I_\aaaa$.

Now assume that $i_*$ is induced by some twisting data.
In this case $\partial X=\partial X^\prime$, $\rho=\rho^\prime$,
$\delta_\Gamma=\delta_\Gamma^\prime$. Moreover,
$\vp^\prime=\vp\otimes\pi$. Let $\nu$ be the highest $A$-weight of $\pi$.
Then $\mu-\nu=\lambda$ and $\delta_\vp+\nu=\delta_{\vp^\prime}$.
It is at this point where we use the assumption that $\pi$ is irreducible.
We observe that 
\begin{eqnarray*}
r_{p,\lambda}(\Gamma)^\prime&=&r_{p,\mu}(\Gamma)   \\
(\:\:\:r^0_{p,\lambda}(\Gamma)^\prime&=&r^0_{p,\mu}(\Gamma)\:\:\:)   \ .
\end{eqnarray*}
Combining the regularity of $i_*\phi$ asserted in Theorem \ref{stmain} with
Lemma \ref{tw} we see that $\phi$ has the regularity claimed in Theorem
\ref{main} again provided that $p$ is generic or $\phi$ is a cusp form and
$\lambda\not\in I_\aaaa$.

In the general case, where $i_*$ is the composition of several embeddings and
twistings, we argue by induction. 

In order to drop the assumption that $p$ is generic we argue as follows.
Let $s<r_{p,\lambda}(\Gamma)$. Then we find $\tilde p>p$
such that Theorem \ref{stmain} applies and $s<r_{\tilde p,\lambda}(\Gamma)$.
The argument above gives $\phi\in H^{\tilde p,<r_{\tilde p,\lambda}(\Gamma)}$
By the embedding theorem $H^{\tilde p,<r_{\tilde p,\lambda}(\Gamma)}\subset 
H^{p,s}$. Since $s<r_{p,\lambda}(\Gamma)$ was arbitrary,
we conclude that $\phi\in H^{p,<r_{p,\lambda}(\Gamma)}$.
\hB

\subsection{Sobolev regularity for cusps}\label{pur}

In the present subsection we assume that $\Gamma=\Gamma_P$ for some (uniquely
determined) $\Gamma$-cuspidal parabolic subgroup $P$. We assume that the cusp
defined by $\Gamma$ has smaller rank. We have a natural embedding
$$S_{\Gamma}(\sigma_\lambda,\vp)\subset  
B_{\Gamma }(\sigma_\lambda,\vp) \subset  {}^\Gamma
C^\infty(\Omega_\Gamma,V(\sigma_\lambda,\vp))\ .$$ Under certain conditions
$f\in B_{\Gamma }(\sigma_\lambda,\vp)$ ($\:\:$ $f\in 
S_{\Gamma  }(\sigma_\lambda,\vp)$ $\:\:$) defines a regular distribution
$\hat f\in C^{-\infty}(\partial X,V(\sigma_\lambda,\vp))$. Our goal in the
present subsection is to determine these conditions and to estimate the
regularity of $\hat f$.

Let $N\subset P$ be the unipotent radical of $P$. Then we have the 
 connected subgroup $N_\Gamma\subset N$ which is
$AP_\Gamma$-invariant and such that $\Gamma\backslash N_\Gamma$ is compact.
Here $P$ acts on $N$ by $(p,x)\mapsto p.x$ such that $pxwP=p.xwP$.
Let $2\rho_\Gamma\in\aaaa^*$ denote the
character of the action of $A$ on $\Lambda^{max} \naaa_\Gamma$, where
$\naaa_\Gamma$ denotes the Lie-algebra of $N_\Gamma$. 
This is consistent with the definition given earlier in the present paper.
Furthermore,  recall that
$\rho^\Gamma:=\rho-\rho_\Gamma$. Note that 
$\delta_\Gamma=-\rho^\Gamma$.

Let $S(N/N_\Gamma)$ be the unit sphere with respect to
some choice of an euclidean $M_\Gamma$-invariant  metric on the vector space
$N/N_\Gamma$. Furthermore, we let $U$ be the preimage of $S(N/N_\Gamma)^{A_+}$
under the projection $N\rightarrow N/N_\Gamma$. Then $Uw\subset \Omega_\Gamma$
is a $\Gamma$-invariant subset which projects onto a neighbourhood $B_\Gamma^+$
of the cusp of $B_\Gamma$.

We fix a maximal compact subgroup $K$ and a corresponding Langlands
decomposition $P=MAN$. Let $w\in K$ be a representative
of the non-trivial element of the Weyl group $N_K(M)/M$. Then we define the
group $\bar N:=N^w$ with Lie algebra $\bar \naaa$. For $r\in\nat_0$ let  $(\bar
X_j)_j$ be a base of $\cU(\bar\naaa)^{\ge -r\alpha}$ (note that $X\in
\bar\naaa$ has degree $-\alpha$). If $\supp(f)\subset B_\Gamma^+$, then $\hat
f$ is supported near $\infty_P$.  It therefore  makes sense to consider the
norm $${}_{\bar N}\|\hat f\|_{H^{p,r}}^p:=\sum_{j} \int_{\bar N}
|\hat f(\bar x \bar X_j)|^p d\bar x\ ,$$ where we have fixed a norm on
$V_\vp$.

\begin{theorem}\label{pc}
 Let $f\in   B_{\Gamma}(\sigma_\lambda,\vp)$
($\:\:$ $f\in S_{\Gamma}(\sigma_\lambda,\vp)$ $\:\:$)
 be supported in
$B_\Gamma^+$.
\begin{enumerate}
\item
 If $\rho^\Gamma>\delta_\vp$ ($\:\:$ no condition $\:\:$ ) and
$\Ree(\lambda)>-\rho^\Gamma+\delta_\vp$, then $f$ determines a
regular distribution $\hat f$. 
\item  If $p\in
(1,\infty)$ and $r,r^0\in\nat_0$ satisfy 
 $$r\alpha < \min\left( \Ree(\lambda)-\rho -\delta_\vp
+\frac{\rho+\rho^\Gamma}{p}\:,\: -2\delta_\vp -2\rho_\Gamma
+2\frac{\rho}{p}\right)$$
$$( \:\:\:\: r^0\alpha <   \Ree(\lambda)-\rho
-\delta_\vp +\frac{\rho+\rho^\Gamma}{p}\:\:\:\:) \ ,$$  
then we have $\hat f\in H^{p,r}(\partial X,V(\sigma_\lambda,\vp))$ ($\:\:$
$\hat f\in H^{p,r^0}(\partial X,V(\sigma_\lambda,\vp))$ $\:\:$).
\end{enumerate} 
\end{theorem} \proof Note that $\hat f$ is smooth on $\Omega_\Gamma$.
In order to prove 1.) me must show that $\hat f$ is
locally integrable near $\infty_P$.
Since $ \bar N\ni \bar x \mapsto \bar x P$ defines coordinates of $\partial X$
near $\infty_P$  it therefore suffices to show that 
\begin{equation}\label{fine1} {}_{\bar N}\|\hat f\|_{H^{1,0}}:=\int_{\bar
N\setminus \{1\}} |f(\bar {x})| d\bar{x}<\infty\ .\end{equation} 
 Let $F:N\setminus \{1\}\rightarrow \bar{N}\setminus \{1\}$ be the
diffeomorphism  such that $F(x)m(xw)\alpha(xw)\tilde n(xw)=xw$ with $m(xw)\alpha(xw)\tilde n(xw)\in
MAN$. We have $$f(x w)=f(F(x)m(xw)\alpha(xw)\tilde n(xw))= \sigma(m(xw))^{-1}
\alpha(xw)^{\lambda-\rho} f(F(x))\ . $$ Thus $|f(F(x))|=|f(xw)|
\alpha(xw)^{\rho-\Ree(\lambda)}$.
Furthermore $F^* d\bar{x}=|\det  DF(x)| dx$.
Note that $F(x^a)=F(x)^a$. From this we conclude that 
\begin{equation}\label{detest}|\det
DF(x^a)|=|\det DF(x)| a^{-4\rho}\ . \end{equation}
Now (\ref{fine1}) is equivalent to
$$\int_N  |f(xw)|
\alpha(xw)^{\rho-\Ree(\lambda)} |\det  DF(x)| dx <\infty\ .$$
Taking into account that $\supp(f)\in B_\Gamma^+$ we can write this integral
in double polar coordinates as
 \begin{equation}
\label{finte1}\int_{A_+} \int_{S(N/N_\Gamma)} \int_{A}\int_{S(N_\Gamma)}  
|f(\eta^b\xi^a w)|
\alpha(\eta^b\xi^a w)^{\rho-\Ree(\lambda)} |\det  DF(\eta^b \xi^a)| 
a^{2\rho^\Gamma} b^{2\rho_\Gamma}d\eta db d\xi  da\ .\end{equation}
We choose some compact preimage $\hat S$  of $S(N/N_\Gamma)$ inside $N$.

Recall that $\vp$ extends to a semisimple representation of $A$.
We can and will assume that $V_\vp$
decomposes into eigenspaces $V_{\vp}(\beta)$, $\beta\in\aaaa^*_+\cup \{0\}$,
w.r.t. the action of $A$ so that $a\in A$ acts on $V_{\vp}(\beta)$ by
$a^\beta$, and such that $\max_{V_\vp(\beta)\not=\{0\}}\beta =2\delta_\vp$.  
Let $f_\beta$ be the component of $f$ in
$V_\vp(\beta)$. Furthermore, let $\vp_{\beta_1,\beta_2}$ be the component of
$\vp$ mapping $V_{\vp}(\beta_2)$ to $V_{\vp}(\beta_1)$.
Note that $\vp_{\beta_1,\beta_2}(p)=0$ if $\beta_2>\beta_1$ and $p\in
AP_\Gamma$.  
Since $f$ is $\Gamma$-invariant we have
$$ f(\gamma x w)=\vp(\gamma)f(x w)=\sum_{\beta_1 \ge \beta_2}
\vp_{\beta_1,\beta_2}(\gamma) f_{\beta_2}(x w)\ .$$
We obtain the estimate
$$\sup_{\xi\in\hat S,\eta\in S(N_\Gamma)}|f(\eta^b v^a w)|\le C
\sum_{\beta_1\ge\beta_2} \sup_{v\in V}
|f_{\beta_2}(\xi^a w)|\sup_{\eta\in
S(N_\Gamma)} \|\vp_{\beta_1,\beta_2}(\eta^b)\|\ , $$ where $V\subset N\setminus
N_\Gamma$ is a sufficiently large compact subset, and $C$ is independent of
$a,b$. Furthermore, using  \ref{pdef}, $\alpha(x^aw)=a^2
\alpha(xw)$, and (\ref{detest}),  we get  \begin{eqnarray} 
\sup_{\eta\in S(N_\Gamma)} \|\vp_{\beta_1,\beta_2}(\eta^b)\|&\le& C
b^{\beta_1-\beta_2} \nonumber\\
\sup_{v\in V} |f_{\beta_2}(v^a w)|& \le& p^{\Gamma}_{V,1}(f)
a^{2(\Ree(\lambda)-\rho^\Gamma)+\beta_2} \label{tc3} \\
\sup_{\eta\in
S(N_\Gamma),\xi\in\hat S}  |\alpha(\eta^b\xi^a w)^{\rho-\Ree(\lambda)}|&\le &
C\max(a,b)^ {2(\rho-\Ree(\lambda))}\nonumber \\ \sup_{\eta\in
S(N_\Gamma),\xi\in\hat S}   |\det  DF(\eta^b \xi^a)| &\le&C \max(a,b)^{-4\rho}
 \ .\nonumber \end{eqnarray}  Thus for  $a\le b$ we have uniformly in
$(\eta,\xi)$ \begin{eqnarray*}\lefteqn{ |f(\eta^b\xi^a w)| \alpha(\eta^b\xi^a
w)^{\rho-\Ree(\lambda)} |\det  DF(\eta^b \xi^a)|  a^{2\rho^\Gamma}
b^{2\rho_\Gamma}}&&\\&\le&C \sum_{\beta_1\ge\beta_2}
b^{\beta_1-\beta_2}a^{2(\Ree(\lambda)-\rho^\Gamma)+\beta_2}
b^{2(\rho-\Ree(\lambda))} b^{-4\rho}a^{2\rho^\Gamma} b^{2\rho_\Gamma}\ .
\end{eqnarray*} Similarly, for $b\le a$ we have  uniformly in $(\eta,\xi)$
\begin{eqnarray*}\lefteqn{ |f(\eta^b\xi^a w)| \alpha(\eta^b\xi^a
w)^{\rho-\Ree(\lambda)} |\det  DF(\eta^b \xi^a)|  a^{2\rho^\Gamma}
b^{2\rho_\Gamma}}&&\\&\le&C \sum_{\beta_1\ge\beta_2}
b^{\beta_1-\beta_2}a^{2(\Ree(\lambda)-\rho^\Gamma)+\beta_2}
a^{2(\rho-\Ree(\lambda))} a^{-4\rho}a^{2\rho^\Gamma} b^{2\rho_\Gamma}
\end{eqnarray*} We now easily conclude that (\ref{finte1}) converges if 
$\rho^\Gamma>\delta_\vp$ and 
$\Ree(\lambda)>-\rho^\Gamma+ \delta_\vp$.

If $f\in  S_{\Gamma}(\sigma_\lambda,\vp)$, then
(\ref{tc3}) can be improved to
$$\sup_{v\in V} |f_{\beta_2}(v^a w)| \le q^{\Gamma}_{W,1,k\alpha}(f)
a^{2(\Ree(\lambda)-\rho^\Gamma)+\beta_2-k\alpha}$$
for any $k\in\nat_0$.
We see that in this case (\ref{finte1}) converges if 
$\Ree(\lambda)>-\rho^\Gamma+ \delta_\vp$.

In order to prove assertions 2.) and 3.) of the theorem we extend the argument
above. Let $\bar X\in\cU(\naaa)$ be homogeneous of degree $-\nu\in\aaaa^*$,  
i.e. $\bar X^a=a^{-\nu}\bar X$ for $a\in A$. We then estimate
\begin{equation}\label{fine2} \int_{\bar N\setminus \{1\}} |f(\bar
{x}\bar{X})|^p d\bar{x}  \ .  \end{equation}   We again transform this
integral into an integral over $N$. We define the function
$X:N\setminus\{1\}\rightarrow \cU(\naaa)$ so that the   the equality 
$F(x X(x))=F(x)\bar X$ of differential operators holds true. Let
$X(x)=\sum_{\mu\le \nu} X_\mu(x)$ be its decomposition into homogeneous
components. Because of  \begin{eqnarray*}
F(x^a X(x^a))=F(x^a)\bar X &=& F(x)^a \bar X = a^{\nu} (F(x) \bar X)^a=a^{\nu}
F(x X(x))^a\\&=&a^{\nu}F(x^aX(x)^a)=\sum_{\mu\le \nu} a^{\nu +\mu}  F(x^a
X_\mu(x)) \end{eqnarray*} 
we conclude that
$$X_\mu(x^a)=a^{\mu+\nu}X_\mu(x)\ .$$
Because of 
$\alpha(xw)^{\rho-\lambda} \sigma(m(xw))  f(x w)=  f(F(x))$
we must apply the Leibniz rule in order to express $f(F(x)\bar X)$. 
We obtain a sum of terms of the form
$$h(x) \alpha(xY_1w)^{\rho-\lambda}
\sigma(m(xY_2w))  f(x Y_3 w)\ ,$$ where
$Y_i\in\cU(\naaa)$ is of degree $d_i$, $d_1+d_2+d_3=:\mu\le \nu$, and
$h(x^a)=a^{\mu+\nu}h(x)$
for $a\in A$. 
We must estimate
$$\int_N | h(x)  \alpha(xY_{1}w)^{\rho-\lambda}
\sigma(m(xY_{2}w))  f(x Y_{3} w)|^p  |\det  DF(x)| dx\ .$$
We again write this integral in the form 
\begin{eqnarray}
&&
\int_{A_+} \int_{S(N/N_\Gamma)} \int_{A}\int_{S(N_\Gamma)} 
 | h(\eta^b\xi^a )  \alpha(\eta^b\xi^a Y_{1}w)^{\rho-\lambda} \sigma(m(\eta^b\xi^a
Y_{2}w))  f(\eta^b\xi^a  Y_{3} w)|^p \nonumber\\
&&\hspace{2cm}|\det  DF(\eta^b
\xi^a)|  a^{2\rho^\Gamma} b^{2\rho_\Gamma}d\eta db  d\xi  da\
.\label{finit3}\end{eqnarray} 
Using the $\Gamma$-equivariance of $f$, $d_3\le \nu\le k\alpha$  we
can estimate   \begin{eqnarray}
|f(\eta^b\xi^a  Y_3 w)|&\le& C  \sum_{\beta_1\ge
\beta_2}\sup_{v\in V}|f_{\beta_2}(v^a Y_3 w)|
\sup_{\eta\in S(N_\Gamma)}\|\vp_{\beta_1,\beta_2}(\eta^b)\| \nonumber \\
\sup_{v\in V} |f_{\beta_2}(v^a Y_3 w)|& \le&  p^{\Gamma}_{V,Y_3}(f)
a^{2(\Ree(\lambda)-\rho^\Gamma)+\beta_2-d_3}\label{unten}\\  |h(\eta^b\xi^a
)|&\le& C \max(a,b)^{\mu+\nu}\nonumber\\ |\alpha(\eta^b\xi^a
Y_{1}w)^{\rho-\lambda}|&\le&   C
 \max(a,b)^{2(\rho-\Ree(\lambda))-d_1}\nonumber\\
 |\sigma(m(\eta^b\xi^a Y_{2}w)) |&\le& C\max(a,b)^{-d_2}\nonumber\ .
\end{eqnarray}   

 For $a\le b$ we have uniformly in $(\eta,\xi)$  
\begin{eqnarray*}
\lefteqn{  | h(\eta^b\xi^a )  \alpha(\eta^b\xi^a Y_{1}w)^{\rho-\lambda}
\sigma(m(\eta^b\xi^a Y_{2}w))  f(\eta^b\xi^a  Y_{3} w)|^p 
 |\det  DF(\eta^b \xi^a)| 
a^{2\rho^\Gamma} b^{2\rho_\Gamma}}&&\\
&\le &  C   p^{\Gamma}_{V,Y_3}(f)^p
  \sum_{\beta_1\ge \beta_2} b^{p(\mu+\nu)}
b^{p[2(\rho-\Ree(\lambda))-d_1]} b^{-p d_2}   a^{p
[2(\Ree(\lambda)-\rho^\Gamma)+\beta_2-d_3]} b^{p(\beta_1-\beta_2)}  
b^{-4\rho}a^{2\rho^\Gamma} b^{2\rho_\Gamma}  \ .  \end{eqnarray*}
Similarly for $b \le a$ we have uniformly in $(\eta,\xi)$   
\begin{eqnarray*}
\lefteqn{  | h(\eta^b\xi^a )  \alpha(\eta^b\xi^a Y_{1}w)^{\rho-\lambda}
\sigma(m(\eta^b\xi^a Y_{2}w))  f(\eta^b\xi^a  Y_{3} w)|^p 
 |\det  DF(\eta^b \xi^a)| 
a^{2\rho^\Gamma} b^{2\rho_\Gamma}}&&\\
&\le &  C  
p^{\Gamma}_{V,Y_3}(f)^p  \sum_{\beta_1\ge \beta_2}
a^{p(\mu+\nu)} a^{p[2(\rho-\Ree(\lambda))-d_1]} a^{-p d_2}   a^{p
[2(\Ree(\lambda)-\rho^\Gamma)+\beta_2-d_3]} b^{p(\beta_1-\beta_2)}  
a^{-4\rho}a^{2\rho^\Gamma} b^{2\rho_\Gamma}  \ .  \end{eqnarray*}
>From these two estimates we conclude that (\ref{finit3}) converges for
$\nu <
\min(\Ree(\lambda)-\rho-\delta_\vp+\frac{\rho+\rho^\Gamma}{p},
-2\delta_\vp+2\rho^\Gamma-2\frac{\rho}{q})$.  
This implies assertion 2.) provided that $\hat f$ belongs to
$H^{p,r}$.

If $f\in S_{\Gamma}(\sigma_\lambda,\vp)$,
then (\ref{unten}) can be improved to 
$$\sup_{v\in V} |f_{\beta_2}(v^a Y_3 w)| \le  q^{\Gamma}_{V,Y_3,k\alpha}(f)
a^{2(\Ree(\lambda)-\rho^\Gamma)+\beta_2-k\alpha}$$
for any $k\in\nat_0$.
We see that (\ref{finit3}) now converges for
$\nu < \Ree(\lambda)-\rho-\delta_\vp+\frac{\rho+\rho^\Gamma}{p}$.  
We thus obtain assertion  2.)  in the case that
$f$ belongs to the Schwartz space again provided that $\hat f$ belongs to
$H^{p,r^0}$.

In order to show that $\hat f\in H^{p,r}$ ($\:\:$ resp. $\hat f\in
H^{p,r^0}$ $\:\:$)  we must show that its distributional derivatives up
to order $r$ are indeed regular distributions. We proceed as follows. We
construct a family of cut-off functions  $\kappa_n\in
C_c^\infty(\Omega_\Gamma)$, $n\in\nat$, such that $\kappa_n\rightarrow 1$
locally uniformly on $\Omega_\Gamma$ as $n\to \infty$. We then have
$\lim_{n\to\infty }\kappa_n \hat f=\hat f$ in the sense of distributions. For
$\bar X\in \cU(\bar \naaa)$ of degree $-\nu$, $\nu \le r\alpha$ we consider 
$(\kappa_n \hat f)(\bar x\bar X)=\kappa_n(x)  \hat f(\bar x\bar X) +  r_n(\bar
x)$. Thus $r_n$ subsumes all terms involving  at least one derivative of
$\kappa_n$. The estimate above already shows that $ \kappa_n(x)  \hat f(\bar
x\bar X) $ converges in $H^{p,0}$ to   $\widehat{ f(.\bar X)}$. We will then
show that $r_n$ tends to zero in $H^{p,0}$. This proves that the
distributional derivatives up to order $r$ of $\hat f$ belong to $H^{p,0}$.

Let $\chi\in C^\infty(A)$ be a cut-off function such that $\chi(a)\equiv 1$ for
$a\le 1$ and $\chi(a)\equiv 0$ for $a^\alpha\ge 2$. Furthermore, fix  $1<a_0\in
A$. Then we define $\tilde \kappa_n(x)\in C^\infty(N\setminus \{1\})$ by
$\tilde \kappa_n(x):=\chi(\alpha(xw) a_0^{-n})$. Furthermore, we set  
$\kappa_n(xw) := \tilde \kappa_n(x)$.
In order to estimate $\|r_n\|_{H^{p,r}}$ we proceed as above but
we replace $f(xw)$ by $\tilde\kappa_n(x) f(xw)$. Then
$$(  \kappa_n  f)(\eta^b\xi^a  Y_3 w) =\sum_{k} \tilde 
\kappa_n(\eta^b\xi^a  Z^\prime_k)   f(\eta^b\xi^a  Z_k w)$$ where $Z_k,Z^\prime_k\in\cU(\naaa)$,
$\deg(Z_k)+\deg(Z^\prime_k)=\deg(Y_3)$ are defined by the decomposition of the
coproduct $\Delta(Y_3)=\sum_k Z_k\otimes Z_k^\prime$.  The terms which
contribute to $r_n$ have $\deg(Z^\prime_k)>0$. In this case we have for
$0<\epsilon\in\aaaa^*$ \begin{eqnarray}
|\tilde  \kappa_n(\eta^b\xi^a  Z^\prime_k)| &=&\chi(\alpha(\eta^b\xi^a 
Z^\prime_kw) a_0^{-n })\nonumber\\
&\le &C \max(a,b)^{-\deg(Z^\prime_k)+\epsilon} a_0^{-n \epsilon/2}\label{eeeedd}\\ 
 |f(\eta^b\xi^a  Z_k w)|&\le & C \sum_{\beta_1\ge \beta_2}
a^{2(\Ree(\lambda)-\rho^\Gamma)+\beta_2-\deg(Z_k)} b^{\beta_1-\beta_2}
\nonumber\ . \end{eqnarray}
In order to see (\ref{eeeedd}) note that $ \chi(\alpha(\eta^b\xi^a 
Z^\prime_kw)a_0)\not=0$ implies that $ a_0 \alpha(\eta^b\xi^aw) \in
[1, a_1]$, where $a_1^\alpha=2$. 
We can choose $\epsilon$ so small that
$(\nu+\epsilon)<\min(\Ree(\lambda)-\rho-\delta_\vp+\frac{\rho+\rho^\Gamma}{p},
-2\delta_\vp+2\rho^\Gamma-2\frac{\rho}{q})$.
Then we get a bound
$$\|r_n\|_{H^{p,r}}\le C a_0^{-n \epsilon/2}\ .$$

In the case that $f\in S_{\Gamma}(\sigma_\lambda,\vp)$
the argument is similar.
\hB

\subsection{Sobolev regularity of invariant functions}

In this subsection $\Gamma$ is a general geometrically finite torsion-free
subgroup of $G$. We assume that all cusps have smaller rank. Again, we 
have   inclusions
$$S_{\Gamma}(\sigma_\lambda,\vp)\subset
B_{\Gamma}(\sigma_\lambda,\vp)\subset {}^\Gamma
C^\infty(\Omega_\Gamma,V(\sigma_\lambda,\vp))\ .$$ 
Under certain
conditions $f\in B_\Gamma(\sigma_\lambda,\vp)$ ($\:\:$ $f\in
S_\Gamma(\sigma_\lambda,\vp)$ $\:\:$ ) determines a regular distribution $\hat
f\in C^{-\infty}(\partial X,V(\sigma_\lambda,\vp))$. Our goal ist to determine
sufficient conditions for this to happen and to estimate the  regularity of
$\hat f$.

Let $\cP$ denote the set of $\Gamma$-conjugacy classes of $\Gamma$-cuspidal
parabolic subgroups. It parametrizes the set of cusps of $\Gamma$. For
$P\in \cP$ we define $\Gamma_P:=\Gamma\cap P$. Then the results of Subsection
\ref{pur} can be applied to  $\Gamma_P$.

 The main result of the present subsection is the following theorem.

\begin{theorem}\label{hypsum}
Let $f\in   B_\Gamma(\sigma_\lambda,\vp)$ ($\:\:$ $f\in
S_\Gamma(\sigma_\lambda,\vp)$ $\:\:$). 
\begin{enumerate}
\item If $\max_{[P]_\Gamma \in\cP}
\left(\delta_{\vp_{|\Gamma_P}}-\rho^{\Gamma_P}\right)<0$
( $\:\:$ no condition $\:\:$ )
 and
$\Ree(\lambda)>\delta_\Gamma+\delta_\vp$, then  $f$ determines a regular
distribution $\hat f$.
\item
If $p\in
(1,\infty)$  and $r, r^0\in \nat_0$, $r \le r^0$,
 $$r\alpha < \min\left( \Ree(\lambda)-\rho -\delta_\vp
+\frac{\rho-\delta_\Gamma}{p}\:,\: -2\max_{[P]_\Gamma
\in\cP}\left[\delta_{\vp_{|\Gamma_P}} +\rho_{\Gamma_P}\right]
+2\frac{\rho}{p}\right)$$ $$( \:\:\:\: r^0\alpha <   \Ree(\lambda)-\rho
-\delta_\vp +\frac{\rho-\delta_\Gamma}{p} \:\:\:\:) \ ,$$  
then we have $\hat f\in
H^{p,r}(\partial X,V(\sigma_\lambda,\vp))$ ($\:\:$ $\hat f\in
H^{p,r^0}(\partial X,V(\sigma_\lambda,\vp))$ $\:\:$). 
\end{enumerate} 
\end{theorem}
\proof
We first prove 1.). It suffices to show that
\begin{equation}\label{kloo2}\int_{\Omega_\Gamma M} |f(k)| dk <\infty\
.\end{equation}
Let $Y_\Gamma:=\Gamma \backslash X$ be the locally symmetric space of $\Gamma$
and $\bar Y_\Gamma:=\Gamma\backslash( X\cup \Omega_\Gamma)$ be its geodesic
compactification.  Let $P\in \cP$ and $\bar Y_P$ be a representative of the
end of $\bar Y_\Gamma$ corresponding to the cusp labeled by $P$ such that we
have an embedding $e_P:\bar Y_P \hookrightarrow \bar Y_{\Gamma_P}$. We choose
a maximal compact subgroup $K$ of $G$ such that $wP\in\Omega_\Gamma$. Let
$\cO\in X$ denote the origin of $X$ determined by this choice of $K$. We let
$\Gamma^P\subset \Gamma$ be a set of representatives of  $\Gamma_P\backslash
\Gamma$, such that if $\gamma\in\Gamma^P$, then $\dist_X(\cO,\gamma\cO)\le 
\dist(\cO,\gamma^\prime \cO) $ for all $\gamma^\prime\in [\gamma]$    (comp.
\cite{bunkeolbrich99}, Sec. 4.2).   Let $B_P:=B_\Gamma\cap \bar Y_P$.
We first assume that $\supp(f)\subset B_P$.
Since $\hat f$ is smooth in a neighbourhood of $wP$ in order to prove
(\ref{kloo2}) we must verify that 
$$ \int_{\bar V} |f(\bar x)| d\bar x <\infty \ ,$$
where $\bar V\subset \bar N$ is a sufficiently
large compact subset.

We consider the lift $(e_P^{-1})^*f=:f_P\in
B_{\Gamma_P}(\sigma_\lambda,\vp)$. Then we have 
$$f=\sum_{\gamma\in\Gamma^P} \pi^{\sigma_\lambda,\vp}(\gamma^{-1})
(f_P)_{|\Omega_\Gamma}$$ as smooth functions on $\Omega_\Gamma$, where the sum
is locally finite. 
 Since $\delta_{\vp_{|\Gamma_P}}-\rho^{\Gamma_P}<0$ ($\:\:$ or
$f$ belongs to the Schwartz space $\:\:$) and
$\Ree(\lambda)>\delta_\Gamma+\delta_\vp\ge
\delta_{\vp_{|\Gamma_P}}-\rho^{\Gamma_P}$ we already know by Theorem \ref{pc}
 that $f_P$ determines a regular distribution $\widehat {(f_P)} \in L^1$.  
Let $G\setminus wP \ni g=\bar n(g)m(g)\alpha(g) \tilde n(g)\in \bar NMAN$ be
the Bruhat decomposition of $g$. 
We now compute using $\frac{d\bar n(\gamma^{-1}\bar x)}{d\bar x}
=\alpha(\gamma^{-1}\bar x)^{-2\rho}$ and $\alpha(\gamma\bar n(\gamma^{-1}\bar
x))=\alpha(\gamma^{-1} \bar x)^{-1}$
\begin{eqnarray*} \int_{\bar V} |f(\bar x)|
d\bar x &\le &\sum_{\gamma\in \Gamma^P} \int_{\bar V} |\vp(\gamma)^{-1}  
f_P(\gamma \bar x)| d\bar x\\ &=&\sum_{\gamma\in \Gamma^P} \int_{\bar V}
|\vp(\gamma)^{-1}   \alpha(\gamma \bar x) ^{\lambda-\rho}
\sigma(m(\gamma \bar x)^{-1})f_P(\bar n(\gamma \bar x)) |d\bar x\\ 
&=&\sum_{\gamma\in \Gamma^P} \int_{\bar V} |\vp(\gamma)^{-1}    f_P(\bar
n(\gamma \bar x))|\alpha(\gamma \bar x) ^{\Ree(\lambda)-\rho}d\bar x\\ 
&= &\sum_{\gamma\in
\Gamma^P} \int_{\bar n(\gamma \bar V)} |\vp(\gamma)^{-1}   f_P(\bar x)|
\alpha(\gamma \bar n(\gamma^{-1} \bar x)) ^{\Ree(\lambda)-\rho} \frac{d\bar
n(\gamma^{-1}\bar x)}{d\bar x} d\bar x\\ 
&=& \sum_{\gamma\in \Gamma^P}
\int_{\bar n(\gamma \bar V)} |\vp(\gamma)^{-1}   f_P(\bar x)| 
\alpha(\gamma^{-1}\bar x)^{-\Ree(\lambda)-\rho} d\bar x\\
&\le& \sum_{\gamma\in \Gamma^P} \|\vp(\gamma^{-1}\| \int_{\bar N}
 |f_P(\bar x)|  \alpha(\gamma^{-1}\bar
x)^{-\Ree(\lambda)-\rho}  d\bar x\ .
\end{eqnarray*}
In order to proceed further we employ the following estimate.
Let $g=k_ga_gh_g\in KA_+K$ denote the Cartan decomposition of $g$. 
\begin{lem}\label{geoestP}
There is $c\in A_+$ such that
$c^{-1} a_\gamma \le \alpha(\gamma^{-1} \bar x)\le c a_\gamma$
uniformly for all $\bar x\in\supp(f_P)$ and $\gamma\in \Gamma^P$.
\end{lem}
Assuming the lemma and using
\begin{eqnarray*}
\|\vp(\gamma)^{-1} \| &\le&C_\epsilon a_\gamma^{\delta_\vp+\epsilon}\\
\int_{\bar N}
 |f_P(\bar x)|  d\bar x &<&\infty
\end{eqnarray*}
for arbitrary small $\aaaa^*\ni\epsilon>0$
we obtain
$$\int_{\bar V} |f(\bar x)|
d\bar x \le C \sum_{\gamma\in \Gamma^P}
a_\gamma^{-\Ree(\lambda)-\rho+\delta_\vp+\epsilon}\ .$$
The sum on the right hand side is finite if
$\Ree(\lambda)>\delta_\Gamma+\delta_\vp+\epsilon$.

This argument for all $P\in\cP$ shows that
 $\hat f$ is a regular distribution if
$\supp(f)$ is concentrated near the cusps of $B_\Gamma$.

Assume now that $f\in C_c^\infty(B_\Gamma,V_{B_\Gamma}(\sigma_\lambda,\vp))$.
We choose some cut-off function $\chi \in C_c^\infty(\Omega_\Gamma)$
such that $\sum_{\gamma\in\Gamma} \gamma^*\chi \equiv 1$ on $\supp(f)$.
Then we can write
$$f=\sum_{\gamma\in\Gamma} \pi^{\sigma_\lambda,\vp}(\gamma^{-1}) (\chi
f)\ .$$
We further choose any parabolic subgroup $P\subset G$ and a maximal compact
subgroup $K$ such  that 
\begin{equation}\label{parachoose} wP\in\Omega_\Gamma\setminus
\supp(f)\ .\end{equation} This subgroup will be denoted by $P_1$ below.
In order to prove that $f$ defines a regular distribution $\hat f$
is suffices to show that
$$\int_{\bar V} |f(\bar x)| d\bar x <\infty$$
for a sufficiently large compact subset $\bar V\subset \bar N$.
We proceed in a similar manner as above.
Define $\tilde\chi$ by $\tilde \chi(\bar n(g))=\chi(g)$.
\begin{eqnarray*} \int_{\bar V} |f(\bar x)|
d\bar x &\le &\sum_{\gamma\in \Gamma} \int_{\bar V} |\vp(\gamma)^{-1}  
  \chi(\gamma \bar x)    f(\gamma \bar x)| d\bar x\\ &=&\sum_{\gamma\in
\Gamma} \int_{\bar V}\tilde \chi(\bar n( \gamma  \bar x)) |\vp(\gamma)^{-1}  
\alpha(\gamma \bar x) ^{\lambda-\rho} \sigma(m(\gamma \bar x)^{-1})f(\bar
n(\gamma \bar x))| d\bar x\\  &=&\sum_{\gamma\in \Gamma} \int_{\bar V}\tilde
\chi(\bar n( \gamma  \bar x))  |\vp(\gamma)^{-1}    f(\bar n(\gamma \bar
x))|\alpha(\gamma \bar x) ^{\Ree(\lambda)-\rho}d\bar x\\  &= &\sum_{\gamma\in
\Gamma} \int_{\bar n(\gamma \bar V)}\tilde \chi(  \bar x)  
|\vp(\gamma)^{-1}   f(\bar x)| \alpha(\gamma \bar n(\gamma^{-1} \bar x))
^{\Ree(\lambda)-\rho} \frac{d\bar n(\gamma^{-1}\bar x)}{d\bar x} d\bar x\\ 
&=& \sum_{\gamma\in \Gamma}
\int_{\bar n(\gamma \bar V)} \tilde \chi(  \bar x)  |\vp(\gamma)^{-1}   f(\bar
x)|  \alpha(\gamma^{-1}\bar x)^{-\Ree(\lambda)-\rho} d\bar x\\
&\le& \sum_{\gamma\in \Gamma} \|\vp(\gamma^{-1}\| \int_{\bar N}\tilde \chi(  \bar x)  
 |f(\bar x)|  \alpha(\gamma^{-1}\bar
x)^{-\Ree(\lambda)-\rho}  d\bar x\ .
\end{eqnarray*}
We now use
\begin{lem}\label{geoestH}
There is $c\in A_+$ such that
$c^{-1} a_\gamma \le \alpha(\gamma^{-1} \bar x)\le c a_\gamma$
uniformly for all $\bar x\in\supp(\chi f)$ and $\gamma\in \Gamma$.
\end{lem}
Again assuming the lemma and using
\begin{eqnarray*}
\|\vp(\gamma)^{-1} \| &\le&C_\epsilon a_\gamma^{\delta_\vp+\epsilon}\\
\int_{\bar N}
 \tilde \chi(\bar x) |f(\bar x)|  d\bar x &<&\infty
\end{eqnarray*}
for arbitrary small $\aaaa^*\ni\epsilon>0$
we obtain
$$\int_{\bar V} |f(\bar x)|
d\bar x \le C \sum_{\gamma\in \Gamma}
a_\gamma^{-\Ree(\lambda)-\rho+\delta_\vp+\epsilon}\ .$$
As above the  sum on the right hand side is finite if
$\Ree(\lambda)>\delta_\Gamma+\delta_\vp+\epsilon$.

Any element $f\in   B_\Gamma(\sigma_\lambda,\vp)$ ($\:\:$ $f\in
S_\Gamma(\sigma_\lambda,\vp)$ $\:\:$) can be
decomposed into a sum of functions which are either supported near the cusps of
$B_\Gamma$ or in the interior of $B_\Gamma$. Hence, we have shown the first
part of the theorem asserting that $\hat f$ is a regular distribution.

We will now show $\hat f$ belongs to $H^{p,r}$ 
( $\:\:$ $H^{p,r^0}$ $\:\:$). 
So far we have shown that $\hat f\in H^{1,0}$. 
In order to obtain the stronger estimate we modify
the proof above. We give the details in the case that $f\in \tilde
B_\Gamma(\sigma_\lambda,\vp)$.

Let $P\in \cP$ and assume that $\supp(f)\subset B_P$.
Let $\bar X\in \cU(\bar \naaa)$ be homogeneous of degree $-\nu$ such
that $\nu\le r\alpha$. We now must show
that $$\|f\|_{P,\bar V,\bar X,p}^p:=\int_{\bar V} |f(\bar x \bar X)|^p d\bar x
<\infty$$ Let $\Delta:\cU(\bar\naaa)\rightarrow \cU(\bar \naaa)\otimes\cU(\bar
\naaa)$ be the coproduct. We write $\Delta\otimes \id(\Delta (\bar
X))=\sum_{k} A_k\otimes B_k\otimes C_k$. Then   we compute 
\begin{eqnarray*} \lefteqn{\int_{\bar V} |f(\bar x
\bar X)|^p d\bar x}&&\\ &\le &\sum_{\gamma\in \Gamma^P} \int_{\bar V}
|\vp(\gamma)^{-1}   f_P(\gamma \bar x\bar X)|^p d\bar x\\ &=&\sum_{\gamma\in
\Gamma^P} \int_{\bar V} |\vp(\gamma)^{-1}   \sum_k \alpha(\gamma \bar x A_k)
^{\lambda-\rho} \sigma(m(\gamma \bar x B_k)^{-1})f_P(\bar n(\gamma \bar x
C_k))|^p d\bar x\\  &\le &\sum_{\gamma\in \Gamma^P} \sum_k \int_{\bar V}
\|\vp(\gamma)^{-1}\|^p   | f_P(\bar n(\gamma \bar x C_k))|^p\|\sigma(m(\gamma
\bar x B_k)^{-1})\|^p |\alpha(\gamma \bar x A_k)^{\lambda-\rho}|^p d\bar x\\ 
&= &\sum_{\gamma\in \Gamma^P} \sum_k \int_{\bar n(\gamma \bar
V)}\|\vp(\gamma)^{-1}\|^p \|\sigma(m(\gamma \bar n(\gamma^{-1} \bar x)
B_k)^{-1})\|^p |\alpha(\gamma \bar n(\gamma^{-1} \bar x)
A_k)^{\lambda-\rho}|^p \\&&\hspace{1cm} | f_P(\bar n(\gamma \bar n(\gamma^{-1}\bar x)
C_k))|^p   \alpha(\gamma^{-1}\bar x)^{-2\rho} d\bar x \ . \end{eqnarray*}
We compute 
\begin{eqnarray*}
m(\gamma \bar n(\gamma^{-1} \bar x)
B_k)&=&m(\gamma \gamma^{-1} \bar x \tilde n(\gamma^{-1} \bar x)^{-1}
\alpha(\gamma^{-1} \bar x)^{-1} m(\gamma^{-1} \bar x)^{-1}  B_k)\\
&=&m(\gamma^{-1} \bar x)^{-1} m(B_k^{[\tilde n(\gamma^{-1} \bar
x)^{m(\gamma^{-1} \bar x)}]^{-1} \alpha(\gamma^{-1} \bar x)^{-1} })\\
&=&\alpha(\gamma^{-1} \bar x)^{-\deg(B_k)} m(\gamma^{-1} \bar x)^{-1}
m(B_k^{[\tilde n(\gamma^{-1} \bar x)^{m(\gamma^{-1} \bar x)}]^{-1}})\\
 \alpha(\gamma \bar n(\gamma^{-1}
\bar x) A_k)&=&
\alpha(\gamma \gamma^{-1} \bar x \tilde n(\gamma^{-1} \bar x)^{-1}
\alpha(\gamma^{-1} \bar x)^{-1} m(\gamma^{-1} \bar x)^{-1}  A_k)\\
&=&\alpha(\gamma^{-1} \bar x)^{-\deg(A_k)}  \alpha(A_k^{[\tilde n(\gamma^{-1}
\bar x)^{m(\gamma^{-1} \bar x)}]^{-1}})\\
 \bar n(\gamma \bar n(\gamma^{-1}\bar x)
C_k)&=&\bar n(\gamma \gamma^{-1} \bar x \tilde n(\gamma^{-1} \bar x)^{-1}
\alpha(\gamma^{-1} \bar x)^{-1} m(\gamma^{-1} \bar x)^{-1}  C_k)\\
&=&\alpha(\gamma^{-1}\bar x)^{-\deg(C_k)} \bar x \bar n(C_k^{[m(\gamma^{-1}\bar x
)\tilde n(\gamma^{-1}
\bar x)]^{-1}})\ .
\end{eqnarray*}

We now employ the following lemma.
\begin{lem}\label{geoestS1}
There exists a compact subset $U\subset N$ such that 
$\tilde n(\gamma^{-1}\bar x ) \in U$
for all $\bar x\in\supp(f_P)$ and $\gamma\in\Gamma^P$.
\end{lem}

Assuming the lemma for a moment we obtain uniformly for
$\bar x\in\supp(f_P)$ and $\gamma\in\Gamma^P$
\begin{eqnarray*}
\|\sigma(m(\gamma \bar n(\gamma^{-1} \bar x)
B_k)^{-1})\|&\le& C \alpha(\gamma^{-1}\bar x)^{-\deg(B_k)}\\
|\alpha(\gamma \bar n(\gamma^{-1} \bar x)
A_k)^{\lambda-\rho}|&\le&C   
\alpha(\gamma^{-1}\bar x )^{-\deg(A_k)} \alpha(\gamma \bar
x)^{-\Ree(\lambda)+\rho}\ . \end{eqnarray*}
Furthermore, using  Theorem \ref{pc} we have
uniformly for $\gamma\in\Gamma^P$
\begin{eqnarray*}
\int_{\bar N} 
 \lefteqn{ | f_P(\bar n(\gamma \bar n(\gamma^{-1}\bar x)
C_k))|^p  d\bar x}&&\\ &\le&C   \sup_{\bar x\in\supp(f_P)}
\alpha(\gamma^{-1}\bar x)^{-p\deg(C_k)}   {}_{\bar N}
\|\widehat{f_P}\|_{H^{p,\frac{\deg(C_k)}{\alpha}}}^p \\ &\le
&C^\prime
 \sup_{\bar x\in\supp(f_P)} \alpha(\gamma^{-1}\bar
x)^{-p\deg(C_k)}  \
. \end{eqnarray*}  
 
Using in addition Lemma \ref{geoestP} we find
$$\|f\|_{P,\bar V,X,p} 
 \le C_\epsilon
  \left(\sum_{\gamma\in \Gamma^P} \sum_k
a_\gamma^{p(-\deg(C_k)-\deg(B_k)
-\deg(A_k)-\Ree(\lambda)+\rho+\delta_\vp+\epsilon)-2\rho}
\right)^{\frac{1}{p}}$$  for arbitrary small  $\epsilon>0$.
Since $-\deg(C_k)-\deg(B_k)-\deg(A_k)=-\mu\le r \alpha$ the sum over
$\Gamma^P$ on the right hand side is finite if  $$r \alpha <  \Ree(\lambda)
-\rho -\delta_\vp +\frac{\rho-\delta_\Gamma}{p}\ .$$
We conclude that $\hat f\in H^{p,r}$ if
$f\in \tilde B_\Gamma(\sigma_\lambda,\vp)$ is supported near the cusps.

In the case that $f$ belongs to the Schwartz space the obvious modification
of the argument above yields 
$\hat f\in H^{p,r^0}$.  

Assume now that $f\in C_c^\infty(B_\Gamma,V(\sigma_\lambda,\vp))$.
We  choose the parabolic subgroup $P\subset G$ and the  maximal compact
subgroup $K$ such  that (\ref{parachoose}) holds.
We must show that 
$$\|\hat f\|_{P,V,\bar X,p}^p:=\int_{\bar V} |f(\bar x \bar X)|^p d\bar x
<\infty\ .$$ We compute  
  \begin{eqnarray*} \lefteqn{\int_{\bar V} |f(\bar x
\bar X)|^p d\bar x}&&\\ &\le &\sum_{\gamma\in \Gamma} \int_{\bar V} 
|\vp(\gamma)^{-1}   (\chi f)(\gamma \bar x\bar X)|^p d\bar x\\
&=&\sum_{\gamma\in \Gamma} \int_{\bar V} |  \sum_k \vp(\gamma)^{-1}  
\alpha(\gamma \bar x A_k) ^{\lambda-\rho} \sigma(m(\gamma \bar x
B_k)^{-1}) (\tilde \chi f)(\bar n(\gamma \bar x C_k))|^p d\bar x\\  &\le
&\sum_{\gamma\in \Gamma} \sum_k \int_{\bar V} \|\vp(\gamma)^{-1}\|^p   |
(\tilde \chi f)(\bar n(\gamma \bar x C_k))|^p\|\sigma(m(\gamma \bar x
B_k)^{-1})\|^p |\alpha(\gamma \bar x A_k)^{\lambda-\rho}|^p d\bar x\\  &=
&\sum_{\gamma\in \Gamma} \sum_k \int_{\bar n(\gamma \bar
V)}\|\vp(\gamma)^{-1}\|^p \|\sigma(m(\gamma \bar n(\gamma^{-1} \bar x
)B_k)^{-1})\|^p |\alpha(\gamma \bar n(\gamma^{-1} \bar x
)A_k)^{\lambda-\rho}|^p \\&&\hspace{1cm}  | (\tilde \chi f)(\bar n(\gamma \bar n(\gamma^{-1}\bar x)
C_k))|^p   \alpha(\gamma^{-1}\bar x)^{-2\rho} d\bar x  \ .\end{eqnarray*}

\begin{lem}\label{geoestS2}
There exists a compact subset $U\subset N$ such that 
$\tilde n(\gamma^{-1}\bar x )\in U$
for all $\bar x\in\supp(\chi f)$ and $\gamma\in\Gamma$.
\end{lem}

Using this lemma   we  conclude that for arbitrary small  $\epsilon>0$
$$\|f\|_{P,\bar V,\bar X,p}^p  \le  C_\epsilon
  \sum_{\gamma\in \Gamma} \sum_k
a_\gamma^{p(-\deg(C_k)-\deg(B_k)-\deg(A_k)-\Ree(\lambda)+\rho+\delta_\vp+\epsilon)-
2\rho}$$ for arbitrary small  $\epsilon>0$.  Note that
$\deg(C_k)+\deg(B_k)+\deg(A_k)=\deg(X)$.
We conclude that
$\hat f\in H^{p,r^0}$.  

Again, since any  element $f\in   B_\Gamma(\sigma_\lambda,\vp)$ ($\:\:$
$f\in S_\Gamma(\sigma_\lambda,\vp)$ $\:\:$) can be
decomposed into a sum of functions which are either supported near the cusps of
$B_\Gamma$ or in the interior of $B_\Gamma$  we have shown  that $\hat f\in
H^{p,r}$ ($\:\:$ 
$\hat f\in H^{p,r^0}$ $\:\:$), if $r\in\nat_0$ ($\:\:$
$r^0\in\nat_0$ $\:\:$) satisfy the assumptions of the theorem.

\proof[Lemma \ref{geoestP}]
For $g\in G$ let $g=\kappa(g)a(g)n(g)\in KAN$ be the Iwasawa decomposition.
Furthermore, if  $g\in G\setminus wP$, then we have the Bruhat decomposition 
$g=\bar n(g)m(g)\alpha(g)\tilde n(g)\in \bar NMAN$. 
In particular, 
\begin{eqnarray*}
\kappa(g)a(g)n(g)&=&\bar n(g)m(g)\alpha(g)\tilde n(g)\\
&=&\kappa(\bar n(g))a(\bar n(g)) n(\bar n(g))m(g)\alpha(g)\tilde n(g)\\
&=& \kappa(\bar n(g))m(g) a(\bar n(g))\alpha(g) n(\bar
n(g))^{\alpha(g)^{-1}m(g)^{-1}}\tilde n(g)\ ,
\end{eqnarray*}
so that $a(g)=a(\bar n(g))\alpha(g)$. Since $a(\bar n)\ge 1$
we conclude that $\alpha(g)\le a(g)$.

Note that there is at most one $\gamma_0\in\Gamma^P$ such that $wP\in
\gamma^{-1} \supp(f_P)$. If this element exists, then we define
$\Gamma_0^P:=\Gamma^P\setminus \{\gamma_0\}$. If not, then we set
$\Gamma_0^P:=\Gamma^P$.

There is a compact subset $V\subset \bar N$ such that $\bar n(\gamma^{-1} \bar
x)\in V$ for all $\bar x\in\supp(f_P)$ and $\gamma\in\Gamma_0^P$. 
In particular, we have $a(\bar n(\gamma^{-1}\bar x))\le c_1$ for some $c_1\in
A$. We conclude 
that 
$$c_1^{-1} a(\gamma^{-1} 
\bar x )\le \alpha(\bar n(\gamma^{-1}\bar
x))\le   a(\gamma^{-1} 
\bar x)\ .$$
Furthermore note that $a(\gamma^{-1}\bar x)=a(\gamma^{-1}\kappa(\bar x))
a(\bar x)$, and $1\le  a(\bar x)\le c_2$ for all $\bar x\in\supp(f_P)$.
Thus
\begin{equation}\label{nr1} c_1^{-1} a(\gamma^{-1} \kappa(
\bar x ))\le \alpha(\bar n(\gamma^{-1}\bar
x))\le  c_2 a(\gamma^{-1} 
\kappa(\bar x))\ .\
\end{equation}

We now want to apply  Lemma \cite{bunkeolbrich00}, Lemma 2.3 asserting that
there is $c_3\in A_+$ such that for all $\gamma\in \Gamma^P$ and $\bar
x\in\supp(f_P)$
 \begin{equation}\label{nr2}c_3^{-1} a_\gamma \le
a(\gamma^{-1}\kappa(\bar x)) \le a_\gamma\ .\end{equation}
In order to verify the assumption of this Lemma we must show that
$WA_+K\cap\Gamma^P$ is finite, where $W$ is some neighbourhood of
$\kappa(\supp(f_P))M$. Note that $\supp f_P\subset \Gamma_P e_P(B_P)$. We
choose $W$ such that $\kappa(\supp(f_P))M\subset W\subset \kappa(\Gamma_P
e_P(B_P)) M$.

We argue by contradiction. So assume that $WA_+K\cap\Gamma^P$ is 
infinite. Then we can find a sequence $\gamma_i\in WA_+K\cap\Gamma^P$ such
that $\gamma_i\rightarrow \Lambda_{\Gamma^P}$, where $\Lambda_{\Gamma^P}$
denotes the set of accumulation points in $\partial X$ of $\Gamma^P\cO$. Let
$D(\cO,\Gamma_P)$ be the Dirichlet domain of $\Gamma_P$ with respect to the
origin. By construction of $\Gamma^P$ we have $\Lambda_{\Gamma^P}\subset
\overline{D(\cO,\Gamma_P)}\setminus \Gamma_P e_P(\bar Y_P)$. Since $W\subset
\kappa(\Gamma_P e_P(B_P))M$ we have $\gamma_i\not\in WA_+K$
for $i>>0$. This is the contradiction.  

Combining (\ref{nr1}) and (\ref{nr2}) we obtain the assertion of Lemma
\ref{geoestP}.\hB  

\proof[Lemma \ref{geoestH}]
Note that there is at most one $\gamma_0\in\Gamma$ such that $wP\in
\gamma^{-1} \supp(\chi f)$. If this element exists then we define
$\Gamma_0:=\Gamma\setminus \{\gamma_0\}$. If not, then we set
$\Gamma_0:=\Gamma$.

There is a compact subset $V\subset \bar N$ such that $\bar n(\gamma^{-1} \bar
x)\in V$ for all $\bar x\in\supp(\chi f)$ and $\gamma\in\Gamma_0$.
We argue as in the proof of Lemma \ref{geoestP} that 
there is $c_1\in A_+$ such that  
\begin{equation}\label{nr3}c_1^{-1} a(\gamma^{-1}\kappa(\bar
x))\le \alpha(\bar n(\gamma^{-1}\bar
x))\le  c_1 a(\gamma^{-1}\kappa(\bar
x))\ .\end{equation}
We again want to apply  Lemma \cite{bunkeolbrich00}, Lemma 2.3, which gives
$c_2\in A_+$ such that for all $\gamma\in\Gamma$ and $\bar x\in\supp (\chi f)$
\begin{equation}\label{nr4}c_2^{-1} a_\gamma \le a(\gamma^{-1}\kappa(\bar
x)) \le a_\gamma\ .\end{equation}
In order to verify the assumption of this Lemma we must show that
$WA_+K\cap\Gamma$ is finite, where $W$ is some neighbourhood of
$\kappa(\supp(\chi f))M$. Note that $\supp (\chi f)\subset \Omega_\Gamma$.
We choose a compact $W$ such that $\kappa(\supp(\chi f))M\subset W\subset
\kappa (\Omega_\Gamma)M$.

We again argue by contradiction. So assume that $WA_+K\cap\Gamma$ is 
infinite. Then we can find a sequence $\gamma_i\in WA_+K\cap\Gamma$ such
that $\gamma_i\rightarrow \Lambda_{\Gamma}$.
 Since $W\subset
\kappa(\Omega_\Gamma)M$ we have $\gamma_i\not\in WA_+K$
for $i>>0$. This is the contradiction.  

Combining (\ref{nr3}) and (\ref{nr4}) we obtain the assertion of Lemma
\ref{geoestH}.\hB  

\proof[Lemma \ref{geoestS1}]
We have 
\begin{eqnarray*}
(\gamma^{-1} \bar x)^{-1}  wP&=&\tilde n(\gamma^{-1} \bar x)^{-1}
\alpha(\gamma^{-1} \bar x)^{-1}m(\gamma^{-1} \bar x)^{-1} \bar n(\gamma^{-1} \bar x)^{-1} wP\\&=&\tilde n(\gamma^{-1} \bar x)^{-1} wP
\end{eqnarray*}
Assume that there are sequences $\bar x_i\in\supp(f_P)$ and 
$\gamma_i\in\Gamma^P$ such that $\tilde n(\gamma^{-1}_i \bar x_i)\rightarrow
\infty$. Then  $\lim_{i\to\infty}\bar x_i^{-1} \gamma_i wP= \infty_P$.
Taking a subsequence we can assume that $\bar x_i$ converges to $\bar x_0\in\supp(f_P)$.
Then $\lim_{i\to\infty} \gamma_i wP = \bar x_0 \infty_P\in \supp(f_P)$.
This is impossible since $\Lambda_{\Gamma^P}\cap\supp(f_P)=\emptyset$.
\hB

\proof[Lemma \ref{geoestS2}]
The proof is similar to that of Lemma \ref{geoestS1}.
Assume that there are sequences $\bar x_i\in\supp(\chi f)$ and 
$\gamma_i\in\Gamma$ such that $\tilde n(\gamma^{-1}_i \bar x_i)\rightarrow
\infty$. Then  $\lim_{i\to\infty}\bar x_i^{-1} \gamma_i wP= \infty_P$.
Taking a subsequence we can assume that $\bar x_i$ converges to $\bar
x_0\in\supp(\chi f)$. Then $\lim_{i\to\infty} \gamma_i wP = \bar x_0
\infty_P\in \supp(f_P)$. This is impossible since
$\Lambda_{\Gamma}\cap\supp(\chi f)=\emptyset$. \hB

We now use the freedom to choose the exponent $p$ in the assertion of
Theorem \ref{hypsum} in order to extend the result to fractional
Sobolev order. 
\begin{kor}\label{hypsum1}
1.$\:$ 
 Let $f\in   B_\Gamma(\sigma_\lambda,\vp)$
and assume that $\max_{[P]_\Gamma \in\cP}
\left(\delta_{\vp_{|\Gamma_P}}-\rho^{\Gamma_P}\right)<0$
and
$\Ree(\lambda)>\delta_\Gamma+\delta_\vp$.
In addition we assume that for $p\in (1,\infty)$ we have
that $$\epsilon:=\left[\Ree(\lambda)-\rho -\delta_\vp
+\frac{\rho-\delta_\Gamma}{p}\right]-\left[-2\max_{[P]_\Gamma
\in\cP}\left[\delta_{\vp_{|\Gamma_P}} +\rho_{\Gamma_P}
\right]+2\frac{\rho}{p}\right]$$ satisfies
$|\frac{\epsilon}{\alpha}|>\frac{\rho+\delta_\Gamma}{\rho-\delta_\Gamma}$,
$\frac{\rho-\delta^\Gamma}{p} > \alpha$, $\frac{\rho-\delta^\Gamma}{q} >
\alpha$, and there
is an integer $k\in\nat$ such that
\begin{eqnarray}
&&0<\left[\Ree(\lambda)-\rho -\delta_\vp
+\frac{\rho-\delta_\Gamma}{p}\right]< k\alpha \le \left[-2\max_{[P]_\Gamma
\in\cP}\left[\delta_{\vp_{|\Gamma_P}} +\rho_{\Gamma_P}
\right]+2\frac{\rho}{p}\right] \quad  \mbox{if}\:\:\epsilon<0 \label{k1}\\
&&0<\left[-2\max_{[P]_\Gamma
\in\cP}\left[\delta_{\vp_{|\Gamma_P}} +\rho_{\Gamma_P}
\right]+2\frac{\rho}{p}\right]<k\alpha\le \left[\Ree(\lambda)-\rho
-\delta_\vp +\frac{\rho-\delta_\Gamma}{p}\right] \quad \mbox{if}\:\:
\epsilon>0\ .\label{k2}\end{eqnarray} Then $\hat f\in H^{p,r}$ provided $r\in
\R$ satisfies  $$r\alpha <\min\left( \Ree(\lambda)-\rho -\delta_\vp
+\frac{\rho-\delta_\Gamma}{p}\:,\: -2\max_{[P]_\Gamma
\in\cP}\left[\delta_{\vp_{|\Gamma_P}} +\rho_{\Gamma_P}\right]
+2\frac{\rho}{p}\right)\ .$$ \\
2.$\:\:$
Let $f\in   S_\Gamma(\sigma_\lambda,\vp)$
and assume that  
$\Ree(\lambda)>\delta_\Gamma+\delta_\vp$.
In addition we assume that for $p\in (1,\infty)$ we have
 $\frac{\rho-\delta_\Gamma}{p}>\alpha$, $\frac{\rho-\delta^\Gamma}{q} >
\alpha$, and $\Ree(\lambda)-\rho
-\delta_\vp +\frac{\rho-\delta_\Gamma}{p}>0$. 
Then $\hat f\in H^{p,r^0}$
provided $r^0\in\R$ satisfies    $$ r^0\alpha <   \Ree(\lambda)-\rho
-\delta_\vp +\frac{\rho-\delta_\Gamma}{p} \ .$$  
\end{kor}
\proof
We begin with the first assertion in case $\epsilon>0$.
Let $k\in \nat$ be the minimal number such that (\ref{k2}) holds true.
For sufficiently small $\delta>0$
we define $1<p_0<p<p_1<\infty$ such that
\begin{eqnarray*}
(k-1)\alpha+\delta&=&\left[-2\max_{[P]_\Gamma
\in\cP}\left[\delta_{\vp_{|\Gamma_P}} +\rho_{\Gamma_P}
\right]+2\frac{\rho}{p_1}\right]\\
 k\alpha+\delta&=&\left[-2\max_{[P]_\Gamma
\in\cP}\left[\delta_{\vp_{|\Gamma_P}} +\rho_{\Gamma_P}
\right]+2\frac{\rho}{p_0}\right]\ .
\end{eqnarray*}
If $\delta$ is sufficiently small, then we have 
\begin{eqnarray*}
(k-1)\alpha&<&\min\left[\left[\Ree(\lambda)-\rho
-\delta_\vp +\frac{\rho-\delta_\Gamma}{p_1}\right],\left[-2\max_{[P]_\Gamma
\in\cP}\left[\delta_{\vp_{|\Gamma_P}} +\rho_{\Gamma_P}
\right]+2\frac{\rho}{p_1}\right]\right]\\
 k\alpha&<&\min\left[\left[\Ree(\lambda)-\rho
-\delta_\vp +\frac{\rho-\delta_\Gamma}{p_0}\right],\left[-2\max_{[P]_\Gamma
\in\cP}\left[\delta_{\vp_{|\Gamma_P}} +\rho_{\Gamma_P}
\right]+2\frac{\rho}{p_0}\right]\right]\ .
\end{eqnarray*}
We now apply Theorem \ref{hypsum} in order to conclude that
$\hat f\in H^{p_0,k}\cap H^{p_1,k-1}$.
We define $\theta\in (0,1)$ such that
$\frac{1}{p}=\frac{1}{p_1}+\theta(\frac{1}{p_0}-\frac{1}{p_1})$.
Then interpolation gives $\hat f\in H^{p,s}$ for $s:=k-1+\theta$.
Note that $s$ depends on $\delta$.
If we choose $\delta>0$ sufficiently small, then $s>r$, and the assertion of
the corollary follows.

Let now $\epsilon<0$. 
Let $k\in \nat$ be the minimal number such that (\ref{k1}) holds true.
For sufficiently small $\delta>0$
we define $1<p_0<p<p_1<\infty$ such that
\begin{eqnarray*}
(k-1)\alpha+\delta&=&\left[ \Ree(\lambda)-\rho
-\delta_\vp +\frac{\rho-\delta_\Gamma}{p_1}\right]\\
 k\alpha+\delta&=&\left[ \Ree(\lambda)-\rho
-\delta_\vp +\frac{\rho-\delta_\Gamma}{p_0}\right]\ .
\end{eqnarray*}
If $\delta$ is sufficiently small, then we have 
\begin{eqnarray*}
 (k-1)\alpha&<&\min\left[\left[\Ree(\lambda)-\rho
-\delta_\vp +\frac{\rho-\delta_\Gamma}{p_1}\right],\left[-2\max_{[P]_\Gamma
\in\cP}\left[\delta_{\vp_{|\Gamma_P}} +\rho_{\Gamma_P}
\right]+2\frac{\rho}{p_1}\right]\right] \\
 k\alpha&<&\min\left[\left[\Ree(\lambda)-\rho
-\delta_\vp +\frac{\rho-\delta_\Gamma}{p_0}\right],\left[-2\max_{[P]_\Gamma
\in\cP}\left[\delta_{\vp_{|\Gamma_P}} +\rho_{\Gamma_P}
\right]+2\frac{\rho}{p_0}\right]\right]\ . 
\end{eqnarray*}
  We now apply Theorem \ref{hypsum} in order to conclude that $\hat f\in
H^{p_0,k}\cap H^{p_1,k-1}$. We define $\theta\in (0,1)$ such that
$\frac{1}{p}=\frac{1}{p_1}+\theta(\frac{1}{p_0}-\frac{1}{p_1})$.
Then interpolation gives $\hat f\in H^{p,s}$ for $s:=k-1+\theta$.
If we choose $\delta>0$ sufficiently small, then $s>r$, and the assertion of
the corollary follows again.

Now we consider the case that $f\in S_\Gamma(\sigma_\lambda,\vp)$.
For sufficiently  small $\delta>0$
we define $1<p_0<p<p_1<\infty$ such that
\begin{eqnarray*}
(k-1)\alpha+\delta&=&\left[ \Ree(\lambda)-\rho
-\delta_\vp +\frac{\rho-\delta_\Gamma}{p_1}\right]\\
 k\alpha+\delta&=&\left[ \Ree(\lambda)-\rho
-\delta_\vp +\frac{\rho-\delta_\Gamma}{p_0}\right]\ ,
\end{eqnarray*}
where $k\in\nat$ is uniquely determined.
We now
apply Theorem \ref{hypsum} in order to conclude that $\hat f\in
H^{p_0,k}\cap H^{p_1,k-1}$. We define $\theta\in (0,1)$ such that
$\frac{1}{p}=\frac{1}{p_1}+\theta(\frac{1}{p_0}-\frac{1}{p_1})$.
Then interpolation gives $\hat f\in H^{p,s}$ for $s:=k-1+\theta$.
If we choose $\delta>0$ sufficiently small, tean $s>r$, and the assertion of
the corollary follows in this case, too.
\hB

We now deduce the following consequence of Corollary \ref{hypsum1}.
\begin{kor}\label{laurent}
Let  $p\in (1,\infty)$ and $\lambda\in \aca$ satisfy the
assumption of Corollary \ref{hypsum1}, 1. 
Furthermore,  let $(f_\mu)_\mu$, $f_\mu \in B_\Gamma(\sigma_\mu,\vp)$, 
be a germ at $\lambda$ of a meromorphic family  and $\hat f_\mu=\sum_{l}
(\mu-\lambda)^l  \hat f_l$ be the Laurent expansion at $\lambda$ of $(\hat
f_\mu)_\mu$. Then $\hat f_l \in H^{p,<r}$,
where $$r\alpha=\min\left( \Ree(\lambda)-\rho -\delta_\vp
+\frac{\rho-\delta_\Gamma}{p}\:,\: -2\max_{[P]_\Gamma
\in\cP}\left[\delta_{\vp_{|\Gamma_P}} +\rho_{\Gamma_P}\right]
+2\frac{\rho}{p}\right)\ .$$
  \end{kor}
\proof
Choose $\epsilon>0$ sufficiently small such that
$F_\mu:=(\mu-\lambda)^{-\ord_{\mu=\lambda} f_\mu} f_\mu$
is defined and holomorphic on $\{|\mu-\lambda|\le \epsilon\}$,
and $r -2\epsilon > 0$.
The proofs of Theorem \ref{hypsum} and  Corollary \ref{hypsum1} show that
$\|F_\mu\|_{H^{p,r-2\epsilon}}$ is locally bounded.
Since $\hat F_\mu=ext ^\Gamma F_\mu$ is a germ of a holomorphic family
of distributions we conclude that $F_\mu$ is weakly continuous in
$H^{p,r-2\epsilon}$.
Since  the embedding
$H^{p,r -2\epsilon}\hookrightarrow
H^{p,r -3\epsilon}$ is compact the family $F_\mu$ is holomorphic
with values in $H^{p,r -3\epsilon}$.
The coefficients $\hat f_l$ can be expressed as
$$\hat f_l=\frac{1}{2\pi\imath}\int_{|\mu-\lambda|=\epsilon}
(\mu-\lambda)^{\ord_{\mu=\lambda} f_\mu-l-1} F_\mu\ ,$$
and therefore $\hat f_l\in H^{p,r -3\epsilon}$.
Since $\epsilon>0$ can be chosen arbitrarily small the assertion follows.\hB

 \subsection{Proof of Theorem \ref{stmain}}

We begin with the first assertion of the theorem.
For a generic set of $p\in (1,\infty)$ we must show that if $\phi\in
\Fam^{st}_\Gamma(\Lambda_\Gamma,\sigma_{\lambda},\vp)$, then  
there exists twisting and embedding data such that $i_*\phi\in
H^{p,<r_{p,\lambda}(\Gamma)^\prime}(\partial
X,V(\sigma^\prime_{\lambda^\prime},\vp^\prime))$.  

From now on $p$ is called generic if $p$ is irrational
and $$\epsilon:=\left[\Ree(\lambda)-\rho -\delta_\vp
+\frac{\rho-\delta_\Gamma}{p}\right]-\left[-2\max_{[P]_\Gamma
\in\cP}\left[\delta_{\vp_{|\Gamma_P}} +\rho_{\Gamma_P}
\right]+2\frac{\rho}{p}\right]\not=0\ .$$
It is easy to see that the set of generic $p$ is a dense subset of
$(1,\infty)$ since either $\rho-\delta_\Gamma\not=2\rho$
or $\Gamma$ is elementary hyperbolic and has no cusps.

Since $\phi$ is stably deformable we can 
choose twisting and embedding data such that
$$(i_{(1)})_*\phi\in {}^\Gamma
C^{-\infty}(\partial
X^{(1)},V^{(1)}(\sigma^{(1)}_{\lambda^{(1)}},\vp^{(1)}))$$ is
deformable.  By twisting further if necessary we can assume that
$\sigma^{(1)}$ is the trivial representation $1\in\hat M^{(1)}$
(see \cite{bunkeolbrich00}, (33)).
Choosing further embedding data we can in addition assume that
$$(i_{(2)})_*(i_{(1)})_*\phi\in {}^\Gamma
C^{-\infty}(\partial
X^{(2)},V^{(2)}(1_{\lambda^{(2)}},\vp^{(1)}))\ ,$$ such that
\begin{eqnarray}
\Ree(\lambda^{(2)})&<&0\label{poi1}\\
\max_{[P]_\Gamma\in
\cP^{(2)}}\left(\delta_{\vp^{(1)}_{|\Gamma_P}}-(\rho^{(2)})^{\Gamma_P}\right)
&<&0\label{poi2} \\
r^{(2)}_{p,\lambda^{(2)}}(\Gamma)&<&0 \label{poi3}\\
 \frac{\rho^{(2)}+\delta^{(2)}_\Gamma}{\rho^{(2)}-\delta^{(2)}_\Gamma}
&<&|\frac{\epsilon}{\alpha}|\label{poo4}\\ 
(\rho^{(2)}-\delta^{(2)}_\Gamma)\min\left(\frac{1}{p}  \:,\:
 \frac{1}{q} \right) &>&\alpha
\label{poo5} \ \ ,\end{eqnarray} and   \begin{equation} \label{poi33}
\min\left( - \Ree(\lambda^{(2)})-\rho^{(2)} -\delta_{\vp^{(2)}}
+\frac{\rho^{(2)}-\delta_\Gamma^{(2)}}{p}\:,\: -2\max_{[P]_\Gamma
\in\cP^{(2)}} \left[ \delta_{\vp^{(2)}_{|\Gamma_P}}
+\rho^{(2)}_{\Gamma_P} \right]+ 2\frac{\rho^{(2)}}{p}   \right)>0 \ .
\end{equation}
Let $d$ denote the codimension of the
embedding. Then $\lambda^{(2)}=\lambda^{(1)}-\frac12 d\alpha$,
$\rho^{(2)}=\rho^{(1)}+\frac12 d\alpha$, 
$(\rho^{(2)})^{\Gamma_P}=(\rho^{(1)})^{\Gamma_P}+\frac12 d\alpha$, 
and  $\delta_\Gamma^{(2)}=\delta_\Gamma^{(1)}-\frac12 d\alpha$.
So we can get (\ref{poi1}) for sufficiently large $d$.
Moreover, we have the relations
$$\max_{[P]_\Gamma \in
\cP^{(2)}}\left(\delta_{\vp^{(1)}_{|\Gamma_P}}-(\rho^{(2)})^{\Gamma_P}\right)
=\max_{[P]_\Gamma\in
\cP^{(1)}}\left(\delta_{\vp^{(1)}_{|\Gamma_P}}-(\rho^{(1)})^{\Gamma_P}\right)-\frac
12 d\alpha\ .$$ Thus we can obtain (\ref{poi2}) for sufficiently large $d$.
Since $r^{(2)}_{p,\lambda^{(2)}}(\Gamma)=  r^{(1)}_{p,\lambda^{(1)}}(\Gamma) -
\frac{d}{q}$ we can also get (\ref{poi3}).
In order to obtain (\ref{poi33}) note that embedding in codimension $d$ 
increases the left hand side by $d\alpha$. For (\ref{poo4}) and
\ref{poo5} note that $\epsilon$ does not change under embedding and twisting,
but $\frac{\rho^{(2)}+\delta^{(2)}_\Gamma}{\rho^{(2)}-\delta^{(2)}_\Gamma}\to
0$ and $ \frac{\alpha}{\rho^{(2)}-\delta_\Gamma^{(2)}}\to 0$ as $d\to\infty$.
If  $\epsilon>0$, then we consider the non-empty interval
$$I^{(1)}:=\left(\left[-2\max_{[P]_\Gamma
\in\cP^{(1)}}\left[\delta_{\vp^{(1)}_{|\Gamma_P}} +\rho^{(1)}_{\Gamma_P}
\right]+2\frac{\rho^{(1)}}{p}\right], \left[\Ree(\lambda^{(1)})-\rho^{(1)} -\delta^{(1)}_\vp
+\frac{\rho^{(1)}-\delta^{(1)}_\Gamma}{p}\right]\right)\ .$$
If $\epsilon<0$, then we set 
$$I^{(1)}:=\left( \left[\Ree(\lambda^{(1)})-\rho^{(1)} -\delta^{(1)}_\vp
+\frac{\rho^{(1)}-\delta^{(1)}_\Gamma}{p}\right],\left[-2\max_{[P]_\Gamma
\in\cP^{(1)}}\left[\delta_{\vp^{(1)}_{|\Gamma_P}} +\rho^{(1)}_{\Gamma_P}
\right]+2\frac{\rho^{(1)}}{p}\right]\right)\ .$$
Then $I^{(2)}=I^{(1)}+\frac{d\alpha}{p}$.
Since $p$ is irrational there are arbitrary large $d$ such that $I^{(2)}$
contains an integer. We now fix $d$ such that the inequalities above are
satisfied and $I^{(2)}$ contains an integer $k\in\nat$.

In order to simplify the notation we now replace
$\phi$ by $(i_{(2)})_*(i_{(1)})_*\phi$, and we omit the superscripts
$()^{(2)}$ everywhere.
Thus we can assume that  
\begin{eqnarray*}
\phi&\in&\Fam_\Gamma(\Lambda_\Gamma,1_\lambda,\vp)\nonumber \\
\Ree(\lambda)&<&0\nonumber \\ 
\max_{[P]_\Gamma \in
\cP}\left(\delta_{\vp_{|\Gamma_P}}-\rho^{\Gamma_P}\right)
&<&0\\
r_{p,\lambda}(\Gamma)&<&0 \\
\min\left( - \Ree(\lambda)-\rho -\delta_{\vp}
+\frac{\rho-\delta_\Gamma}{p}\:,\: -2\max_{[P]_\Gamma
\in\cP} \left[ \delta_{\vp_{|\Gamma_P}}
+\rho_{\Gamma_P} \right]+ 2\frac{\rho}{p}   \right)&>&0
\\\frac{\rho+\delta_\Gamma}{\rho-\delta_\Gamma}
&<&|\frac{\epsilon}{\alpha}|\\ 
(\rho-\delta_\Gamma)\min\left(\frac{1}{p}  \:,\:
 \frac{1}{q} \right) &>&\alpha ,\end{eqnarray*}
and $I$ contains an integer $k$.

We must show that
$\phi\in H^{p,<r_{p,\lambda}(\Gamma)}$.  
Since $\phi$ is deformable and strongly supported on the limit set we can
find a germ at $\lambda$ of a holomorphic family  $(\phi_\mu)_{\mu}$, $
\phi_\mu\in {}^\Gamma C^{-\infty}(\partial X,V(1_\mu,\vp))$,
 such that $\phi_\lambda=\phi$ and
$(res^\Gamma (\phi_\mu))_{|\mu=\lambda}=0$. By Lemma \ref{eqqu} we have  
that $\phi\in \Ext^{sing}_\Gamma(1_\lambda,\vp)$. Since $ext^\Gamma$ has at
most finite-dimensional singularities and $B_\Gamma(1_\lambda,\vp)\subset
D_\Gamma(1_\lambda,\vp)$ is dense  we can find a germ of a holomorphic
family $(\psi_\mu)_{\mu}$, $\psi_\mu \in B_\Gamma(1_\mu,\vp)$, such
that $\psi_\lambda=0$ and $(ext^\Gamma\psi_\mu)_{|\mu=\lambda}=\phi$.
We now define the germ at $-\lambda$ of a meromorphic family $(f_{\mu})_{\mu}$,
$f_\mu\in B_\Gamma(1_\mu,\vp)$, by
$f_{-\mu}:=\hat S^{w,\Gamma}_{1_\mu,\vp} \psi_\mu$. Here
$\hat S^{w,\Gamma}_{1_\mu,\vp}:D_\Gamma(1_\mu,\vp)\rightarrow 
D_\Gamma(1_{-\mu},\vp)$ (note that $1^w=1$) is the scattering matrix defined in
\ref{sctdeff}, and we employ Lemma \ref{bser} in order to see that $f_\mu\in
B_\Gamma(1_\mu,\vp)$. Let $p_1(\lambda)$ be the Plancherel density. It is a
meromorphic function on $\aca$. We then have the following identities:
\begin{eqnarray*}
ext^\Gamma_{-\mu}\circ \hat S^{w,\Gamma}_{1_\mu,\vp} & =&
\hat J^{w}_{1_\mu,\vp}\circ ext^\Gamma_\mu\\
p_1(-\mu) \hat J^{w^{-1}}_{1_{-\mu},\vp}\circ
\hat J^{w}_{1_\mu,\vp}&=&\id \ .\end{eqnarray*}
Thus we can write 
$$\phi=\left( \frac{1}{p_1(-\mu)} \hat
J^{w^{-1}}_{1_{-\mu},\vp}(ext^\Gamma f_{-\mu})\right)_{|\mu=\lambda} \
.$$
Let $\hat f_\mu:=\sum_{k}(\frac{\mu+\lambda}{\alpha})^k \hat f_k$ be the
Laurent expansion of the germ
$(ext^\Gamma f_\mu)_{\mu}$ at $-\lambda$. Then by Corollary \ref{laurent} we
have $\hat f_k\in H^{p,<r^\prime}$,
where $r^\prime$ is determined by
$$r^\prime\alpha:=\min\left( -\Ree(\lambda)-\rho -\delta_\vp
+\frac{\rho-\delta_\Gamma}{p}\:,\: -2\max_{[P]_\Gamma
\in\cP}\left[\delta_{\vp_{|\Gamma_P}} +\rho_{\Gamma_P}\right]
+2\frac{\rho}{p}\right)\ .$$
Let 
$$\frac{1}{p_1(\mu)} \hat
J^{w^{-1}}_{1_{\mu},\vp}=\sum_{l}(\frac{\mu+\lambda}{\alpha})^l B_l$$
be the Laurent expansion of $(\frac{1}{p_1(\mu)} \hat
J^{w^{-1}}_{1_{\mu},\vp})_{\mu}$ at $-\lambda$.
Since $\Ree(\lambda)<0$ we can apply  
Cor. \ref{mapint} in order to get
$B_l:H^{p,<r^\prime}\rightarrow H^{p,<r^\prime+2\frac{\Ree(\lambda)}{\alpha}}$.
Since $r^\prime+2\frac{\Ree(\lambda)}{\alpha}= r_{p,\lambda}(\Gamma)$ and 
$\phi=\sum_{l+k=0} B_l (\hat f_k)$
we conclude that $\phi\in H^{p,<r_{p,\lambda}(\Gamma)}$. 

Now we show the second assertion of Theorem \ref{stmain}.
Let $\phi\in \Cusp_\Gamma(\sigma_\lambda,\vp)$ and $\lambda\not\in
I_\aaaa$. Again, after twisting and embedding we can assume that
all cusps of $\Gamma$ have smaller rank, 
\begin{eqnarray*}
\phi&\in&\Cusp_\Gamma(1_\lambda,\vp)\cap
\Fam_\Gamma(\Lambda_\Gamma,1_\lambda,\vp)\nonumber \\
\Ree(\lambda)&<&0\nonumber \\  
 r^0_{p,\lambda}(\Gamma)&<&0\\
-\Ree(\lambda)-\rho-\delta_\vp+\frac{\rho-\delta_\Gamma}{p}&>&0\\
 (\rho-\delta_\Gamma)\min\left(\frac{1}{p}  \:,\:
 \frac{1}{q} \right) &>&\alpha
.\end{eqnarray*}
Let $(f_\mu)_{\mu}$, $f_\mu\in B_\Gamma(1_\mu,\vp)$, be constructed as
above.
Because we assume $\lambda\not\in I_\aaaa$ this family is regular at
$\mu=-\lambda$. 

We show that $f_{-\lambda}\in S_\Gamma(1_{-\lambda},\vp)$.
Let $[P]_\Gamma\in\cP$. We choose a $P$-invariant the cut-off function
$\chi_P$ in order to define the map $T_P$ (see Subsection \ref{tpp}). Then we
can define the constant term  $$(T_P f_{-\lambda})_P:= \int_{\Gamma_P\backslash
P_\Gamma} \pi^{1_{-\lambda},\vp}(x) T_P f_{-\lambda} \: dx\ .$$
Taking the constant term commutes with multiplication by $\chi_P$,  restriction
and the intertwining operators. Therefore we have
\begin{eqnarray*}
 \left(T_P f_{-\lambda}\right)_P &=& \left(\chi_P res^{\Gamma_P} \hat
J^w_{1_{-\lambda},\vp} \phi_{\lambda}\right)_P  \\
&=&\chi_P  res^{\Gamma_P}  \hat J^w_{1_{-\lambda},\vp}
\left(\phi_{\lambda}\right)_P\ .
\end{eqnarray*}
Since $\phi$ is a cusp form we conclude that $\left(T_P
f_{-\lambda}\right)_P=0$. Since $AS(T_P f_{-\lambda})=AS\left(T_P
f_{-\lambda}\right)_P=0$ for all $P\in \cP$ we see that $f_{-\lambda}\in
S_\Gamma(1_{-\lambda},\vp)$. 

Using   
$$\phi=\left( \frac{1}{p_1(-\lambda)} \hat
J^{w^{-1}}_{1_{-\lambda},\vp}(ext^\Gamma
f^{mod}_{-\lambda})\right)$$ and  
 Corollary \ref{hypsum1}, 2.) 
we obtain  
$ext^\Gamma f_{-\lambda}=\widehat f_{-\lambda} \in H^{p,<r^\prime}$ for
$r^\prime\alpha:=-\Ree(\lambda)-\rho-\delta_\vp+ \frac{\rho-\delta_\Gamma}{p}$.
It follows from Corollary \ref{mapint} and
$r^0_{p,\lambda}(\Gamma)=r^\prime+2\frac{\Ree(\lambda)}{\alpha}$ that
$\phi\in H^{p,<r^0_{p,\lambda}(\Gamma)}$.
\hB

\bibliographystyle{plain}

\end{document}